\documentclass[11pt]{amsart}

\usepackage{amssymb}
\usepackage{color}
\usepackage{graphicx}
\usepackage{float}
\usepackage[all,cmtip]{xy}
\usepackage{tikz}
\usetikzlibrary{matrix}
\usepackage{url}
\usepackage{subfig}
\usepackage{hyperref}
\newcommand{\SD}{{\mathcal{D}}}
\newcommand{\SW}{{\mathcal{W}}}

\newcommand{\Maps}{\operatorname{\mathcal{M}}\mathfrak{aps}}
\newcommand{\Cusp}{\operatorname{\mathcal{C}}\mathfrak{usp}}
\newcommand{\Stereo}{\operatorname{\mathcal{S}}tereo}

\mathchardef\mhyphen="2D

\newcommand{\Diff}{\operatorname{Diff}}

\newcommand{\SO}{\operatorname{SO}}

\newcommand{\Ar}{\operatorname{Area}}
\newcommand{\Int}{\operatorname{Int}}

\newcommand{\im}{{\operatorname{Image}}}
\newcommand{\Id}{{\operatorname{Id}}}

\newcommand{\Rot}{\operatorname{Rot}}
\newcommand{\tb}{\operatorname{tb}}

\newcommand{\tw}{\operatorname{tw}}
\newcommand{\Span}{\operatorname{span}}

\newcommand{\Top}{\operatorname{Top}}

\newcommand{\FCont}{\operatorname{FCont}}
\newcommand{\Res}{\operatorname{Res}}
\newcommand{\Conf}{\operatorname{Conf}}

\newcommand{\R}{{\mathbb{R}}}
\newcommand{\C}{{\mathbb{C}}}
\newcommand{\Z}{{\mathbb{Z}}}

\newcommand{\NS}{{\mathbb{S}}}
\newcommand{\D}{{\mathbb{D}}}

\newcommand{\Op}{{\mathcal{O}p}}

\newcommand{\F}{{\mathcal{F}}}
\newcommand{\PMas}{{\mathcal{P}^+}}

\newcommand{\Leg}{{\mathfrak{Leg}}}
\newcommand{\LegImm}{{\mathfrak{LegImm}}}
\newcommand{\SLegImm}{{\mathfrak{LegImm}^\textbf{strict}}}
\newcommand{\SSLegImm}{{\mathfrak{SLegImm}}}
\newcommand{\SpLegImm}{{\mathcal{S}p\mathfrak{LegImm}}}
\newcommand{\SLegImmG}{{\mathfrak{LegImm}^\textbf{strict}_{0}}}
\newcommand{\FLeg}{{\mathfrak{FLeg}}}
\newcommand{\FLegImm}{{\mathfrak{FLegImm}}}
\newcommand{\FLegO}{{\widehat{\mathfrak{FLeg}}}}

\newcommand{\Hor}{{\mathfrak{Hor}}}
\newcommand{\HorImm}{{\mathfrak{HorImm}}}
\newcommand{\FHor}{{\mathfrak{FHor}}}
\newcommand{\FHorImm}{{\mathfrak{FHorImm}}}
\newcommand{\FHorO}{{\widehat{\mathfrak{FHor}}}}
\newcommand{\FELeg}{{\mathfrak{FELeg}}}
\newcommand{\FELegO}{{\widehat{\mathfrak{FELeg}}}}

\newcommand{\kalman}{{K\'{a}lm\'{a}n }}
\newcommand{\kalmans}{{K\'{a}lm\'{a}n's }}
\newcommand{\Imm}{{\mathfrak{Imm}}}
\newcommand{\Emb}{{\mathfrak{Emb}}}
\newcommand{\LEmb}{{\mathfrak{LEmb}}}

\newtheorem{theorem}{Theorem}[subsection]
\newtheorem{lemma}[theorem]{Lemma}
\newtheorem{proposition}[theorem]{Proposition}
\newtheorem{corollary}[theorem]{Corollary}

\theoremstyle{definition}
\newtheorem{definition}[theorem]{Definition}

\newtheorem{remark}[theorem]{Remark}

\begin{document} 
	
	\title{Fundamental groups of formal Legendrian and Horizontal embedding spaces}
	
	\date{\today}
	

	\author{Eduardo Fern\'{a}ndez}
	\address{Universidad Complutense de Madrid, 
		Departamento de \'{A}lgebra, Geometr\'{i}a y Topolog\'{i}a, Facultad de Matem\'{a}ticas, and Instituto de Ciencias Matem\'{a}ticas CSIC-UAM-UC3M-UCM, C. Nicol\'{a}s Cabrera, 13-15, 28049 Madrid, Spain.}
	\email{eduarf01@ucm.es;eduardo.fernandez@icmat.es}
	
	\author{Javier Mart\'{i}nez-Aguinaga}
	\address{Universidad Complutense de Madrid, 
		Departamento de \'{A}lgebra, Geometr\'{i}a y Topolog\'{i}a, Facultad de Matem\'{a}ticas, and Instituto de Ciencias Matem\'{a}ticas CSIC-UAM-UC3M-UCM, C. Nicol\'{a}s Cabrera, 13-15, 28049 Madrid, Spain.}
	\email{javier.martinez.aguinaga@icmat.es}
	
	\author{Francisco Presas}
	\address{Instituto de Ciencias Matem\'{a}ticas CSIC-UAM-UC3M-UCM, C. Nicol\'{a}s Cabrera, 13-15, 28049 Madrid, Spain.}
	\email{fpresas@icmat.es}

	\begin{abstract}
		
		We compute the fundamental group of each connected component of the space of formal Legendrian embeddings in $\mathbb{R}^3$. We use it to show that previous examples in the literature of non trivial loops of Legendrian embeddings are already non trivial at the formal level. Likewise, we compute the fundamental group of the different connected components of the space of formal horizontal embeddings into the standard Engel $\mathbb{R}^4$.  We check that the natural inclusion of the space of horizontal embeddings into the space of formal horizontal embeddings induces an isomorphism at $\pi_1$--level. 
		
	\end{abstract}
	
	\maketitle
	\setcounter{tocdepth}{2} 
	\tableofcontents
	
	\section{Introduction.}\label{section:Introduction}
	
	The computation of the homotopy type of the space of Legendrian embeddings into a contact $3$--fold has a long story. For a while, it was thought that the computation could be made at the formal level. We mean by that that the inclusion of the space of Legendrian embeddings into the space of formal Legendrian embeddings, ie the space of pairs: smooth embedding and formal Legendrian derivative, was a weak homotopy equivalence. This was proven to be wrong in the key article of D. Bennequin \cite{Bennequin}; in which it was shown that the formal space associated to the standard contact $\R^3$ possesses some connected components that are not representable by Legendrian knots. In other words, the restriction of the induced map of the inclusion at $\pi_0$--level was not surjective. This was the first hint of ridigity phenomena in contact topology. Later on, there has been an industry checking how far the inclusion map is at $\pi_0$--level from being injective or surjective, see, eg,  the work of Chekanov \cite{Chekanov}, Ding and Geiges \cite{DingGeiges}, Eliashberg and Fraser \cite{EliashbergFraser}, Etnyre and Honda \cite{EtnyreHonda} or Osv\'ath, Szab\'o and Thurston \cite{Szabo}.
	
	The next step was the study of higher homotopy groups. This was developed by \kalman \cite{kalman} in dimension $3$ using pseudoholomorphic curves invariants and by Sabloff and Sullivan \cite{Sabloff} in dimension $2n+1$, $n>1$, using generating function invariants. However, they just checked that several non trivial loops in the space of Legendrian embeddings were trivial as elements in the fundamental group of the space of smooth embeddings. We show that all \kalmans examples are non trivial in the space of formal Legendrian embeddings, see Section \ref{section:KalmanLoop}. This makes unnecessary the use of sophisticated invariants to compute these examples. In order to do that, we compute the fundamental group of the space of formal Legendrian embeddings. This is the content of Section \ref{FormalLegendrian}. 
	
	Moreover, the analogous problem for horizontal embeddings into Engel manifolds satisfies an $h$--principle at $\pi_0$--level, eg see Adachi \cite{Ada}, Geiges \cite{GeigesLoops} or del Pino and Presas \cite{PinoPresas}. In Section \ref{FormalHorizontal}, we compute the fundamental group of the space of formal horizontal embeddings into the standard Engel $\R^4$. We check that it is $\Z\oplus\Z_2$ where the first component is a rotation invariant that captures the formal immersion class of the loop and the second component captures the formal embedding class of the loop. In Section \ref{AreaInvariantSection}, we provide a way of combinatorially computing this $\Z_2$--invariant.
	
	Finally, in Section \ref{h--Principle}, we check that the formal fundamental group is isomorphic to the fundamental group, ie the $\Z_2$--invariant  and the rotation invariant completely classify the elements of the fundamental group of horizontal embeddings. Independently, Casals and del Pino have shown that there is a full $h$--principle for horizontal embeddings \cite{CasalsdelPino}. Although our proof is more specific, it has the advantage of not using the topology of the space of smooth embeddings. Therefore, we reprove as a corollary that the space of smooth embeddings of circles into $\R^4$ is simply connected, see Budney \cite{Budney}.
	
	The methods developed in this article may allow to compute higher rank homotopy groups of the space of embeddings of the circle in $\R^4$. The strategy that we are pursuing in a forthcoming project is based on the following facts. First, we compute the higher homotopy groups of the space of horizontal embeddings by the geometrical method developed in this article, that is heavily simplified thanks to the work of Igusa \cite{Igusa}. Secondly, we use \cite{CasalsdelPino} in order to state that the previous computations are also computing homotopy groups of the formal horizontal embedding space. Finally, by using obstruction theory, see Hatcher \cite[Chapter 4.3]{HatcherBook}, we isolate the homotopy groups of the space of smooth embeddings.

	\textbf{Acknowledgements:} The authors are extremely grateful to  R. Casals, V. Ginzburg and A. Del Pino for several discussions clarifying the paper. Also, they thank A. Hatcher for pointing out the useful reference \cite{Budney}. Daniel \'Alvarez-Gavela has helped as quite a lot with enjoyable meetings whenever he comes to Madrid, he has pointed out several sharp remarks. The third author is also grateful to the organizers of the \em Engel Structures \em workshop held in April 2017 (American Institute of Mathematics, San Jose, California) for providing a nice environment in which this article was discussed. Last, but not least, we want to thank the excellent job that the referees have done: they have produced three long and deep reports that have helped to improve the readability of the paper a lot, not to speak of its soundness.
	
	The authors are supported by the Spanish Research Projects SEV--2015--0554, MTM2016--79400--P, and MTM2015--72876--EXP. The first author is supported by a Master--Severo Ochoa grant and by Beca de Personal
	Investigador en Formaci\'on UCM.  The second author is funded by Programa Predoctoral de Formaci\'{o}n de Personal Investigador No Doctor del Departamento de Educaci\'{o}n del Gobierno Vasco.
	
	\section{Spaces of embeddings of the circle into euclidean space.}\label{section:EmbeddingsSmooth}
	Denote by $\Emb(N,M)$ the space of embeddings of a manifold $N$ into a manifold $M$ equipped with the $C^r$--topology, $r\geq 5$.
	
	\subsection{The space $\Emb(\NS^1,\R^3)$.}\label{subsection:EmbeddingsR3}
	
	\begin{theorem}[Hatcher, \cite{HatcherSmale} Apendix: equivalence $(15)$]\label{thm:HatcherUnknot}
		The space of parametrized unknotted circles in $\R^3$ has the homotopy type of $\SO(3)$.
	\end{theorem}
	
	The group $\SO(4)$ acts freely on the connected component $\Emb_{p,q}(\NS^1,\NS^3)\subseteq\Emb(\NS^1,\NS^3)$ of the parametrized \em $(p,q)$ torus knots \em as
	\begin{center}
		$\begin{array}{rccl}
		\SO(4)\times\Emb_{p,q}(\NS^1,\NS^3) &\longrightarrow &  \Emb_{p,q}(\NS^1,\NS^3)\\
		(A,\gamma)& \longmapsto & A\cdot\gamma.
		\end{array}$
	\end{center}
	Thus, we have an inclusion $\SO(4)\hookrightarrow\Emb_{p,q}(\NS^1,\NS^3)$. The following result holds:
	
	\begin{theorem}[Hatcher, \cite{hatcher} Theorem  $1$]\label{thm:HatcherTorusKnots}
		The inclusion $\SO(4)\hookrightarrow\Emb_{p,q}(\NS^1,\NS^3)$ is a homotopy equivalence.
	\end{theorem}
	As a consequence of these results we obtain that 
	
	\begin{corollary}\label{cor:FundamentalGroupKnots}
		Let $\Emb_0(\NS^1,\R^3),\Emb_{p,q} (\NS^1,\R^3)\subseteq\Emb(\NS^1,\R^3)$ be the connected component of the parametrized \em unknots \em or of the parametrized \em $(p,q)$ torus knots\em, respectively. The fundamental groups of these spaces are given by 
		\begin{itemize}
			\item  $\pi_1 (\Emb_0 (\NS^1,\R^3))\cong \Z_2$,
			\item $\pi_1 (\Emb_{p,q}(\NS^1,\R^3))\cong G_{p,q}\rtimes\Z_2$,
		\end{itemize}
		where $G_{p,q}$ is the knot group of the $(p,q)$ torus knot.
	\end{corollary}
	\begin{proof}
		The case of the connected component $\Emb_0 (\NS^1,\R^3)$ follows from Theorem \ref{thm:HatcherUnknot}.
		
		We need to study the connected component $\Emb_{p,q} (\NS^1,\R^3)$ to conclude the proof. Consider the following space
		\[   \Stereo_{p,q}=\{(\gamma, x):x\notin\im(\gamma)\}\subseteq\Emb_{p,q} (\NS^1,\NS^3)\times\NS^3.   \]
		We have two natural fibrations associated to the projection maps
		\begin{displaymath}
		\xymatrix{
			\Emb_{p,q} (\NS^1,\R^3)  \ar@{^{(}->}[r]  &  \Stereo_{p,q} \ar[d]& & \NS^3\backslash K_{p,q} \ar@{^{(}->}[r]  &  \Stereo_{p,q} \ar[d]\\
			&   \NS^3& & &\Emb_{p,q}(\NS^1,\NS^3)}
		\end{displaymath}
		where $K_{p,q}$ is the image of the standard $(p,q)$ torus knot in $\NS^3$. From the first fibration we obtain \[ \pi_1 (\Emb_{p,q} (\NS^1,\R^3))\cong\pi_1 (\Stereo_{p,q}).  \]
		
		Moreover, from the second one and the fact that $\Emb_{p,q} (\NS^1,\NS^3)$ has the homotopy type of $\SO(4)$ (Theorem \ref{thm:HatcherTorusKnots}), we obtain that the sequence
		\begin{displaymath}
		\xymatrix@M=10pt{
			0\ar[r] & G_{p,q}\ar[r] & \pi_1 (\Stereo_{p,q}) \ar[r] & \Z_2 \ar[r]& 0}
		\end{displaymath}
		is exact. Now, it is a simple exercise to check that this sequence is right split.
	\end{proof}
	
	\subsection{The space $\Emb(\NS^1,\R^4)$.}\label{subsection:EmbeddingsR4}
	A \em long embedding \em of $\R$ into $\R^4$ is an embedding $\gamma:\R\rightarrow\R^4$ that coincides with the standard inclusion $\R\hookrightarrow\R\times\R^3=\R^4$ outside a compact neighborhood of the origin. Let $\LEmb(\R,\R^4)$ denote the space of long embeddings of $\R$ into $\R^4$.
		
	\begin{lemma}[Budney]\label{lem:HomotopyKnotsR4}
		$\Emb(\NS^1,\R^4)$ is homotopy equivalent to $\NS^3\times\NS^2\times\LEmb(\R,\R^4)$.
	\end{lemma}
	\begin{proof}
		It follows from \cite[Proposition 2.2]{Budney} that $\Emb(\NS^1,\R^4)$ is homotopy equivalent to $(\SO(4)
\times\widetilde{\Stereo})/\SO(3)$, where $\widetilde{\Stereo}=\{(p,f):p\notin\im(f)\}\subseteq\R^4\times\LEmb(\R,\R^4)$. Observe that the quaternion structure in $\R^4(i,j,k)$ induces a homotopy equivalence $\NS^3\times\widetilde{\Stereo}\rightarrow(\SO(4)
\times\widetilde{\Stereo})/\SO(3);(v,(p,f))\mapsto[((v|iv|jv|kv),(p,f))]$. Finally, the natural fibration $\widetilde{\Stereo}\mapsto\LEmb(\R,\R^4)$ has fiber homotopy equivalent to $\NS^2$, since the space $\LEmb(\R,\R^4)$ is connected and $\R^4 \backslash\R\stackrel{h.e.}{\sim}\NS^2$. Moreover the fibration homotopically splits since given any family $\varphi:\NS^k\rightarrow\LEmb(\R,\R^4)$, it admits a lift given by $\varphi_{\varepsilon}:\NS^k\rightarrow\widetilde{\Stereo};z\mapsto(\varphi(z)(0)+\varepsilon\cdot i\varphi(z)'(0),\varphi(z))$ for some $\varepsilon>0$ small enough.
\end{proof}
	Let us geometrically explain the homotopy equivalence stated in the last Lemma. Let $\{e_1,e_2,e_3,e_4\}$ be the canonical basis of $\R^4$ and write $\R^4=\R\langle e_1 \rangle \times \R^3 \langle e_2,e_3,e_4\rangle$. From a long embedding we obtain an embedding of $\NS^1$ into $\R^4$ closing it in the plane $\langle e_1, p\rangle$, where $p\in\NS^2\subseteq \R^3 \langle e_2,e_3,e_4\rangle$. Finally, the $\NS^3$ factor acts by quaternionic multiplication on a fixed embedding. 
	
	It follows that $\pi_2 (\Emb(\NS^1,\R^4))\cong\pi_2(\NS^2)\oplus\pi_2(\LEmb(\R,\R^4))\cong\Z\oplus\pi_2(\LEmb(\R,\R^4))$. We provide an explicit construction of the generator of the $\Z$ factor in this decomposition.  For every point $p\in\NS^2\subseteq\R^3\langle e_2,e_3,e_4\rangle$ take the standard parametrized \em unknot \em in $\Span\{e_1,p\}$ centered at $p$ and tangent to the line generated by $e_1$ at $0$. This gives us a $2$--parametric family of knots, that we denote by $S_U$, whose homotopy class is the generator $(1,0)$ of $\pi_2(\Emb(\NS^1,\R^4))$. Explicitly, $S_U$ is defined as:
	\begin{center}
		$\begin{array}{rccl}
		S_U\colon & \NS^2 & \longrightarrow &  \Emb(\NS^1,\R^4) \\
		& p& \longmapsto & \gamma_p (t)
		\end{array}$
	\end{center}
	where, if $p=(\lambda_1,\lambda_2,\lambda_3)\in\NS^2$,
	\[
	\gamma_p (t) =\sin(2\pi t)e_1 + (1-\cos(2\pi t))(\lambda_1 e_2+\lambda_2 e_3 + \lambda_3 e_4).
	\]
	\begin{figure}[h] 
		\includegraphics[scale=0.75]{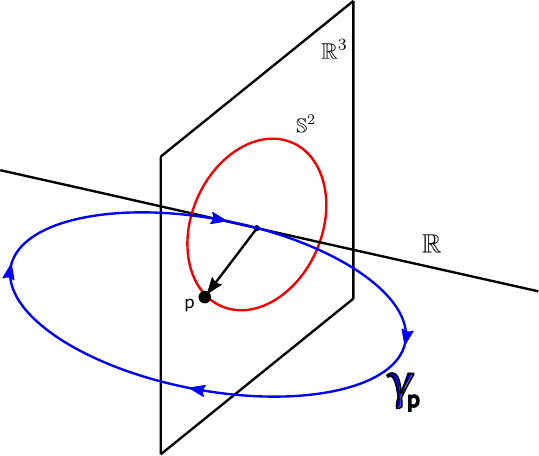}
		\centering
		\caption{Construction of $S_U$ (schematically).}
	\end{figure}
	
	\begin{remark}
		In \cite[Proposition 3.9(4)]{Budney} Budney shows that the space $\Emb(\NS^1,\R^4)$ is simply connected and $\pi_2(\Emb(\NS^1,\R^4))\cong\Z\oplus\Z$. Moreover, he provides an explicit construction of the second generator of the second homotopy group (\cite[Theorem 3.13]{Budney}). 
		
		In section \ref{h--Principle} we provide an alternative proof of the fact that $\Emb(\NS^1,\R^4)$ is simply connected based on the techniques developed in this paper.
	\end{remark}
	
	\section{Formal Legendrian Embeddings in $\R^3$.}\label{FormalLegendrian}
	We denote by $\xi$  the standard contact structure in $\R^3(x,y,z)$ given by $\xi=\ker(dz-ydx)$. Throughout the Section we fix the Legendrian framing $\partial_y$.\footnote{Observe that $\xi$ is naturally (co)oriented by the contact form $dz-ydx$ and, thus, $\partial_y$ determines a unique oriented framing up to homotopy.}
	
	\subsection{Formal Legendrian Embeddings in $\R^3$.}
	
	\begin{definition}
		An immersion $\gamma:\NS^1\rightarrow\R^3$ is said to be \em Legendrian \em if $\gamma'(t)\in\xi_{\gamma(t)}$ for all $t\in\NS^1$. If $\gamma$ is an embedding, we say it is a \em Legendrian embedding\em.
	\end{definition}
	
	\begin{definition}
		\item [(a)]	A \em formal Legendrian immersion \em in $\R^3$ is a pair $(\gamma,F)$ such that: 
		\begin{itemize}
			\item [(i)] $\gamma:\NS^1\rightarrow\R^3$ is a smooth map.
			\item [(ii)] $F:\NS^1\rightarrow \gamma^*(T\R^3\backslash\{0\})$ satisfies $F(t)\in\xi_{\gamma(t)}$.
		\end{itemize}
		
		\item[(b)]	A \em formal Legendrian embedding \em in $\R^3$ is a pair $(\gamma,F_s)$, satisfying: 
		\begin{itemize}
			\item [(i)] $\gamma:\NS^1\rightarrow\R^3$ is an embedding.
			\item [(ii)] $F_s:\NS^1\rightarrow \gamma^*(T\R^3\setminus\lbrace 0\rbrace)$, is a $1$--parametric family, $s\in[0,1]$, such that $F_0=\gamma'$ and $F_1(t)\in\xi_{\gamma(t)}$.
		\end{itemize}
	\end{definition}
	
	Use the framing $\langle \partial_y \rangle $ to trivialize the contact distribution understood as a bundle. This provides a bundle isomorphism $\xi\simeq\R^2$. From now on, we will understand the map $F:\NS^1\rightarrow\NS^1\equiv\NS^2\cap\R^2$ and the family $F_s:\NS^1\rightarrow\NS^2$ with $F_1:\NS^1\rightarrow\NS^1\equiv\NS^2\cap\R^2$. We say that an immersion is \em strict \em if it is a non injective map.
	
	Denote by $\LegImm(\R^3)$ the space of Legendrian immersions in $\R^3$, by $\Leg(\R^3)$ the space of Legendrian embeddings in $\R^3$ and by $\SLegImm(\R^3)=\LegImm(\R^3)\backslash\Leg(\R^3)$ the space of strict Legendrian immersions. Denote also by $\FLegImm(\R^3)$ the space of formal Legendrian immersions and by $\FLeg(\R^3)$ the space of formal Legendrian embeddings. These definitions make sense for immersions and embeddings of the interval. We define $\LegImm([0,1],\R^3)$, $\Leg([0,1],\R^3)$, $\SLegImm([0,1],\R^3)$, $\FLegImm([0,1],\R^3)$ and $\FLeg([0,1],\R^3)$ analogously.
	
	\begin{remark}\label{rem:FrechetImm}
		The space of $C^r$--maps $\Maps^{r}(\NS^1,\R^3)$ endowed with the $C^r$--topology is a Banach vector space. It has an open submanifold that is the space of immersions $\Imm^r(\NS^1,\R^3)$. Denote by $\alpha$ the standard contact form in $\R^3$; ie $\ker(\alpha)=\xi$. There is a smooth submersion $\tilde{\alpha}:\Maps^r(\NS^1,\R^3)\to \Omega^1(\NS^1), \gamma\mapsto\gamma^{*}(\alpha)$. We have that $\LegImm^r(\R^3)$ is a Banach submanifold of $\Imm^r(\NS^1,\R^3)$ (see Lang \cite[Corollary 5.7, page 19]{Lang}). Whenever we speak about families of curves in the space of Legendrians we are implicitly asumming this Banach structure. This means that our maps are $C^r$ instead of smooth but all the definitions in the paper make sense for $C^r$--maps whenever $r\geq 5$. To simplify the notation we will be writing smooth maps instead of $C^5$--maps unless it is stated otherwise. The bound $r\geq 5$ comes from the fact that we are considering $3$--dimensional families of Legendrians. In order to locally classify them we need to consider up to $4$ derivatives of their front projection ($5$ of the original Legendrians).
		
		It is important to note that any $\gamma\in\LegImm^r(\R^3)$ admits a particular type of immersed chart $\varphi_{\gamma}$, called a \em Weinstein chart\em, which identifies via a immersed contactomorphism a tubular neighborhood of the zero section of the jet space $J^1(\NS^1)$ (see Geiges \cite[Example 2.5.11]{GeigesCont}) with a tubular neighborhood of $\gamma(\NS^1)$. 
		This construction provides the local structure produced by the Implicit Function Theorem applied in the last paragraph. We can just check that the exponential map in a $C^r$--neighborhood of $\gamma$ is defined as the following construction: take $f\in \Maps^r(\NS^1,\R)$, a map $\tilde{f}\in \Maps^r(\NS^1,\NS^1)$ such that $\tilde{f}(\theta)=\theta+f(\theta)$. Then, the map is defined as 
		 \begin{equation}\label{eq:ExponentialLegImm}
		 	\begin{array}{rccl}
		 	\mathcal{C}\colon & \mathcal{U}\subseteq \Maps^r(\NS^1, \R) \times \Maps^r(\NS^1, \R) & \longrightarrow &  \mathcal{V}\subseteq\LegImm^r(\R^3) \\
		 	& (f,g)& \longmapsto & \varphi_{\gamma}^{-1}\circ j^1(g\circ\tilde{f}).
		 	\end{array}
		 \end{equation}
	 Note that this is a local chart (ie  a diffeomorphism). However, if we consider $r=\infty$ the previous map is not surjective in the reparametrizations direction since it is well-known that the action of the Lie algebra of $\Diff(\NS^1)$ is not surjective on $\Diff(\NS^1)$. On the other hand, if we quotient out $\LegImm^\infty(\R^3)$ by the action of $\Diff(\NS^1)$, ie we consider non--parametrized Legendrians, the map is a local diffeomorphism. This is an obvious consequence of the characterization of the Legendrian curves in $J^1(\NS^1)$ as graphs of smooth functions. We use $r<\infty$ in order to make sure that the parametrizations are surjective. We do not forget about reparametrizations because they are important in our algebraic topology arguments.
	 
	 From now on we will skip the index $r=5$ in all the statements, unless it is not clear from the context.
	\end{remark}

All the spaces of Legendrians are equipped with the $C^r$--topology. On the other hand, the spaces of formal Legendrians are equipped with the product topology that is the $C^r$--topology for the first factor (the smooth immersion/embedding) and the $C^{r-1}$--topology for the second factor (the formal derivative).

	\begin{remark}
		It is well--known that $h$--principle holds for Legendrian immersions (see, eg, Eliashberg and Mishachev \cite{EliashMisch}). Hence, $\pi_0 (\LegImm(\R^3))\cong\Z$, $\pi_1 (\LegImm(\R^3))\cong\Z$ and $\pi_k (\LegImm(\R^3))=0,$ for all $k\geq2$. The connected components of $\LegImm(\R^3)$ are given by the \em rotation number\em. The rotation number of an Legendrian immersion $\gamma$ is $\Rot(\gamma)=\deg(\gamma':\NS^1\rightarrow\NS^1)$. Let us explain the group $\pi_1 (\LegImm(\R^3))\cong\Z$. Take a loop $\gamma^\theta$ in $\LegImm(\R^3)$, the integer is just $\Rot_L (\gamma^\theta)=\deg(\theta\mapsto(\gamma^\theta)'(0))$, we call this number \em rotation number of the loop\em. These invariants make sense in the formal case and the definitions are the obvious ones. 
	\end{remark}
	
	\subsection{The space $\FLeg(\R^3)$.}
	Consider the space $\FLegO(\R^3)=\{(\gamma,F)| \gamma\in\Emb(\NS^1,\R^3), F\in\Maps(\NS^1,\NS^1)\}$. We have a natural fibration $\FLeg(\R^3)\rightarrow \FLegO(\R^3)$. In order to compute the homotopy groups of $\FLeg(\R^3)$, take $\gamma\in\Leg(\R^3)$ and fix $(\gamma,\gamma')$ as base point. The fiber over this point is $\F=\F_{(\gamma,\gamma')}=\Omega_{\gamma'}(\Maps(\NS^1,\NS^2))$. We have the following exact sequence of homotopy groups associated to the fibration:
	\begin{displaymath}
	\xymatrix@M=10pt{
		& \cdots\ar[r] & \pi_2 (\FLegO(\R^3)) \ar[dll]\\
		\pi_1 (\F) \ar[r] & \pi_1 (\FLeg(\R^3)) \ar[r] & \pi_1 (\FLegO(\R^3)) \ar[dll] \\
		\pi_0 (\F) \ar[r] & \pi_0 (\FLeg(\R^3)) \ar[r] & \pi_0 (\FLegO(\R^3)) \ar[r] & 0 }
	\end{displaymath}
	Notice that $\FLegO(\R^3)$ has the homotopy type of $\Emb(\NS^1,\R^3)\times\NS^1\times\Z$. Hence, $\pi_0 (\FLegO(\R^3))\cong\pi_0 (\Emb(\NS^1,\R^3))\oplus\Z$, where the integer is the rotation number. Moreover, $\pi_1 (\FLegO(\R^3))\cong\pi_1 (\Emb(\NS^1,\R^3))\oplus \Z$ and the $\Z$ factor is given by the rotation number of the loop. Finally, $\pi_k (\FLegO(\R^3))\cong\pi_k (\Emb(\NS^1,\R^3))$ for all $k\geq2$.
	
	The homotopy groups of $\F$ are easily computable. Just observe that there is a fibration $\Maps(\NS^1,\NS^2)\rightarrow\NS^2$ defined via the \em evaluation map\em, with fiber over $p\in\NS^2$ given by $\Omega_p (\NS^2)$. As every element $[f]\in\pi_n (\NS^2)$ can be lifted to an element $[f_n]\in\pi_n (\Maps(\NS^1,\NS^2))$, defined as
	\begin{equation*}
		f_n(p)(t)=f(p), t\in\NS^1, p\in\NS^n,
	\end{equation*}
	all the diagonal maps in the associated exact sequence are zero. This implies that there are short exact sequences $\pi_n (\Omega_p (\NS^2))\rightarrow\pi_n (\Maps(\NS^1,\NS^2))\rightarrow\pi_n(\NS^2)$ for $n\geq1$. In particular, since $\NS^2$ is simply connected, we obtain that \[ \pi_0(\mathcal{F})\cong\pi_1(\Maps(\NS^1,\NS^2))\cong\pi_1 (\Omega_p (\NS^2))\cong\pi_2(\NS^2)\cong\Z.\] Moreover, theses sequences are right split and, thus, split for $n>2$ since the groups involved are abelian. So, we have
	\begin{align*}
		 &\pi_1 (\F)\cong\pi_2(\Maps(\NS^1,\NS^2))\cong\pi_2 (\NS^2)\oplus\pi_2 (\Omega_p (\NS^2))\cong\pi_2 (\NS^2)\oplus\pi_3 (\NS^2)\cong\Z\oplus\Z \\
		\text{and } &\pi_2 (\F)\cong\pi_3(\Maps(\NS^1,\NS^2))\cong\pi_3 (\NS^2)\oplus\pi_3 (\Omega_p (\NS^2))\cong\pi_3 (\NS^2)\oplus\pi_4 (\NS^2)\cong\Z\oplus\Z_2.
	\end{align*}
	
	\begin{lemma}\label{componentsFLeg}
		$\pi_0 (\FLeg(\R^3))\cong\pi_0 (\Emb(\NS^1,\R^3))\oplus\Z\oplus\Z$.
	\end{lemma}
	\begin{proof}
		It is sufficient to show that every element in $\pi_1 (\FLegO(\R^3))\cong\pi_1 (\Emb(\NS^1,\R^3))\oplus\pi_1(\Maps(\NS^1,\NS^1)) \cong\pi_1(\Emb(\NS^1,\R^3))\oplus\Z$ can be lifted to an element in $\pi_1 (\FLeg(\R^3))$. 

		Take a loop $(\gamma^\theta,F_1^{\theta})$ in $\FLegO(\R^3)$. Let $F_0=(\gamma^{\theta})':\NS^1\times\NS^1(\theta,t) \rightarrow\NS^2$ be the derivative $F_0(\theta,t)=(\gamma^\theta)'(t)$, we need to show that $F_0=(\gamma^\theta)'$ is homotopic to the map $F_1:\NS^1\times\NS^1\rightarrow\NS^2$. Observe that the homotopy classes of maps from $\NS^1\times \NS^1$  to $\NS^2$ are classified by the degree and $\deg(F_1)=0$, so we just need to show that $\deg(F_0)=0$ to complete the proof.
		
		Indeed, the map
		\begin{equation}
		G_\varepsilon(\theta,t)=\begin{cases}
		(\gamma^\theta)'(t) & \text{if $\varepsilon=0$,} \\
		\frac{\gamma^\theta(t+\varepsilon)-\gamma^\theta(t)}{||\gamma^\theta(t+\varepsilon)-\gamma^\theta(t)||} & \text{if $0<\varepsilon<1$,}\\
		-(\gamma^\theta)'(t) & \text{if $\varepsilon=1$,}
		\end{cases} \label{eq:homotopyantipodal}
		\end{equation}
		is well--defined, because $\gamma^\theta$, $\theta\in\NS^1$, is an embedding. Thus, $F_0=(\gamma^\theta)'$ is homotopic to $-F_0=-(\gamma^\theta)'$ and $\deg(F_0)=\deg(\gamma^\theta)'=0$.
	\end{proof}
	
	\subsection{Classification of formal Legendrian embeddings in $\R^3$.}
	We have checked that $\pi_0 (\FLeg(\R^3))\cong\pi_0 (\Emb(\NS^1,\R^3))\oplus\Z\oplus\Z$. The first $\Z$ corresponds to the \em rotation number \em and we will show that the second one corresponds to the \em Thurston--Bennequin invariant\em.
	
	Let us refine the definition of formal Legendrian embedding to extend the definition of the \em Thurston--Bennequin invariant \em to the formal case.
	
	\begin{definition}
		A \em formal extended Legendrian embedding \em in $\R^3$ is a pair $(\gamma,G_s)$, satisfying: 
		\begin{itemize}
			\item [(i)] $\gamma:\NS^1\rightarrow\R^3$ is a embedding.
			\item [(ii)] $G_s:\NS^1\rightarrow\SO(3)$, is  a smooth family in the parameter $s\in[0,1]$, such that $G_0=\Id$ and $G_1 (\gamma')\in\xi_{\gamma(t)}$.
		\end{itemize}
	\end{definition}
	
	We denote $\FELeg(\R^3)$ for the space of formal extended Legendrian embeddings in $\R^3$ equipped with the $C^r$--topology in the first factor and the $C^{r-1}$--topology in the second one.
	
	\begin{remark}\label{FELeg=FLeg}
		The natural fibration $t: \FELeg(\R^3)\rightarrow\FLeg(\R^3),(\gamma,G_s)\mapsto(\gamma,G_s(\gamma'))$, has contractible fibers. Thus, this map is a weak homotopy equivalence.
	\end{remark}
	
	Given $(\gamma,G_s)\in\FELeg(\R^3)$ we have a well--defined \em formal contact framing \em $\F_{\FCont}$ of $\nu(G_1(\gamma'))$ given by the Legendrian condition $G_1 (\gamma')\subseteq\xi_{\gamma(t)}$. Then, $G_{1}^{-1}(\F_{\FCont})$ defines a framing of the normal bundle $\nu$ of $\gamma$. On the other hand, we have a \em topological framing \em $\F_{\Top}$ of $\nu$ given by a Seifert surface of $\gamma$.
	
	\begin{definition}
		Let $(\gamma,G_s)\in\FELeg(\R^3)$. The \em Thurston--Bennequin invariant \em is $\tb(\gamma,G_s)=\tw_{\nu}(G_{1}^{-1}(\F_{\FCont}),\F_{\Top})$, ie the twisting of $G_{1}^{-1}(\F_{\FCont})$ with respect to $\F_{\Top}$.
	\end{definition}
	
	The Thurston--Bennequin invariant is defined over $\FELegO(\R^3)=\{(\gamma,G)|\gamma\in\Emb(\NS^1,\R^3),G\in\Maps(\NS^1,\SO(3)),G(\gamma')\in\xi_{\gamma(t)}\}$. Furthermore, since the unique oriented $\NS^1$--bundle over $\NS^1$ is the trivial one, $\pi_0 (\FELegO(\R^3))\cong\pi_0 (\Emb(\NS^1,\R^3))\oplus\Z\oplus\Z$. The first $\Z$ is just the \em rotation number \em and the second one corresponds to the \em Thurston--Bennequin invariant\em.
	
	Now we can state the main result of this Section, which is well--known in the literature.
	
	\begin{theorem}\label{classificationFLeg}
		Formal Legendrian embeddings are classified by their parametrized knot type, rotation number and Thurston--Bennequin invariant.
	\end{theorem}
	The proof of this result follows directly using the fibration  $\hat{F}:\FELeg(\R^3)\rightarrow\FELegO(\R^3)$ and the fact that the map $t: \FELeg(\R^3) \rightarrow \FLeg(\R^3)$ is a weak homotopy equivalence. Note also that the fibration $\hat{F}$ has connected fiber, because its $\pi_0$ is given by $\pi_2(\SO(3))=0$. This completes the proof. However to get a more geometric picture, we will express the isomorphism in more concrete terms. Clearly the isomorphism preserves the rotation invariant, ie the rotation number $(\gamma,G_s)$ is sent to the rotation invariant of $(\gamma,G_s(\gamma'))$. To understand the rest of the isomorphism we fix a base point $(\gamma, F=\gamma')$ in $\FLegO(\R^3)$ with $\Rot(\gamma,\gamma')=0$, i.e we declare the base point to be a Legendrian embedding with zero rotation. Now, given an element of the fiber, ie $(\gamma, F^s)$ with $F^0=F^1=\gamma'$, we claim that for $[(\gamma, 0, k)] \in \pi_0(\FLeg(\R^3))$, the isomorphism $\pi_0(t)$ is given by $\tb(\pi_0(t)^{-1}(\gamma, 0,k))=\tb(\gamma,\gamma')-2k$. In other words, it depends on the choice of base point. This is obvious if we check that given a double stabilization, see Definition \ref{dobleEstab} for a precise definition, 
of the Legendrian knot, the value of the degree invariant in the fiber increases by $1$ and it is a simple computation to check that the $\tb$ decreases by $2$.
	
	\subsection{Fundamental group of formal Legendrian Embeddings in $\R^3$.}
	As a consequence of Lemma \ref{componentsFLeg} we have that the following sequence is exact:
	\begin{displaymath}
	\xymatrix@M=10pt{
		\cdots\ar[r]& \pi_2 (\FLeg(\R^3))\ar[r] & \pi_2 (\FLegO(\R^3)) \ar[dll]\\
		\pi_1  (\F) \ar[r] & \pi_1 (\FLeg(\R^3)) \ar[r] & \pi_1 (\FLegO(\R^3))\ar[r] & 0 }
	\end{displaymath}
	
	Take a $2$--sphere $(\gamma^z, F_{1}^{z})$ in $\FLegO(\R^3)$, the diagonal map $d:\pi_2(\FLegO(\R^3))\rightarrow \pi_1 (\F)$ measures the obstruction to lifting $(\gamma^z, F_{1}^{z})$ to a $2$--sphere $(\gamma^z,F_{s}^{z})$ in $\FLeg(\R^3)$, ie the obstruction to find a homotopy between the derivative map  $F_0=(\gamma^z)':\NS^2\times\NS^1\rightarrow\NS^2$, $(z,t)\mapsto F_0(z,t)=(\gamma^z)'(t)$, and the map $F_1 :\NS^2\times\NS^1\rightarrow \NS^2$. Note that by the Legendrian condition $F_1$ is nullhomotopic, since it is not surjective. The first obstruction to find this homotopy is just the degree of $F_0(z,0)=(\gamma^z)'(0)$ and corresponds to the first $\Z$ factor of $\pi_1 (\F)\cong\pi_2 (\NS^2)\oplus\pi_3(\NS^2)\cong\Z\oplus\Z$. In particular, since $(\gamma^z)'(0)$ is homotopic to $-(\gamma^z)'(0)$ this obstruction vanishes (see equation (\ref{eq:homotopyantipodal})).
	
	\begin{theorem}\label{pi1FLeg}
		The sequence
		\begin{displaymath}
		\xymatrix@M=10pt{
			0\ar[r] & \Z\oplus\Z_m\ar[r] & \pi_1 (\FLeg(\R^3)) \ar[r] & \pi_1 (\Emb(\NS^1,\R^3))\oplus\Z\ar[r] & 0 }
		\end{displaymath}
		is exact, where $m\geq0$. In particular, if we fix the connected component $\Emb' (\NS^1,\R^3)\subseteq \Emb(\NS^1,\R^3)$ of the parametrized \em unknot \em or of the parametrized \em $(p,q)$ torus knot \em  we have that $m=0$ and so
		\begin{displaymath}
		\xymatrix@M=10pt{
			0\ar[r] & \Z\oplus\Z \ar[r] & \pi_1 (\FLeg(\R^3)) \ar[r] & \pi_1 (\Emb'(\NS^1,\R^3))\oplus\Z\ar[r] & 0 }
		\end{displaymath}
		is exact.
	\end{theorem}
	
	\begin{proof}
		We only need to check the particular cases mentioned above. For the connected component $\Emb_0 (\NS^1,\R^3)$ of the parametrized unknot the result follows from Theorem \ref{thm:HatcherUnknot}, since $\pi_2 (\FLegO_0(\R^3))=\pi_2 (\Emb_0(\NS^1,\R^3))=\pi_2 (\SO(3))=0$, where $\FLegO_0 (\R^3)$ stands for a formal Legendrian connected component of the smooth unknot.
		
		On the other hand, fix the connected component $\Emb_{p,q} (\NS^1,\R^3)$ of the parametrized $(p,q)$ torus knot and consider the commutative diagram
		\begin{displaymath} 
		\xymatrix@M=10pt{
			\Emb_{p,q}(\NS^1,\R^3) \ar@{^{(}->}[r]\ar@{^{(}->}[d]& \Emb_{p,q}(\NS^1,\NS^3) \ar@{^{(}->}[d]\\
			\Imm(\NS^1,\R^3) \ar@{^{(}->}[r]& \Imm(\NS^1,\NS^3) }
		\end{displaymath}
		defined by the natural inclusions, where $\Imm(N,M)$ denotes the space of immersions of a manifold $N$ into a manifold $M$ equipped with the $C^r$--topology, $r\geq5$. 
		
		By the Smale--Hirsch Theorem for immersions (see \cite{SmaleHirsch} or \cite{EliashMisch}) we have that $\Imm(\NS^1,\R^3)$ has the homotopy type of $\Maps(\NS^1,\NS^2)$ and $\Imm(\NS^1,\NS^3)$ has the homotopy type of $\Maps(\NS^1,\NS^3)\times\Maps(\NS^1,\NS^2)$. Moreover, the map induced by the inclusion $\Emb_{p,q}(\NS^1,\R^3)\hookrightarrow\Imm(\NS^1,\R^3)$ at $\pi_2$--level sends the homotopy class $[\gamma^z]\in\pi_2(\Emb_{p,q}(\NS^1,\R^3))$ to $[(\gamma^z)']\in\pi_2(\Imm(\NS^1,\R^3))\cong\pi_2(\Maps(\NS^1,\NS^2))$; ie coincides with the diagonal map $d:\pi_2 (\FLegO(\R^3))\rightarrow \pi_1 (\F)\cong\pi_2(\Maps(\NS^1,\NS^2))$. 
		
		Consider the induced commutative diagram at $\pi_2$--level 
		\begin{displaymath} 
		\xymatrix@M=10pt{
			\pi_2 (\Emb_{p,q}(\NS^1,\R^3)) \ar[r]\ar[d]& \pi_2 (\Emb_{p,q}(\NS^1,\NS^3)) \ar[d]\\
			\pi_2 (\Imm(\NS^1,\R^3)) \ar[r]& \pi_2 (\Imm(\NS^1,\NS^3)) }
		\end{displaymath}
		since $\pi_2 (\Emb_{p,q}(\NS^1,\NS^3))$ is trivial (see Theorem \ref{thm:HatcherTorusKnots}) it is sufficient to show that the homomorphism $\pi_2 (\Imm(\NS^1,\R^3))\rightarrow \pi_2 (\Imm(\NS^1,\NS^3))$ is injective to conclude the proof; but this is clear, by using the $h$--principle for immersions, since the degree of the induced map $\NS^2 \times \NS^1$ to $\NS^3$ is zero and the induced map for the derivative from $\NS^2\times \NS^1$ to $\NS^2$ is sent to itself by the inclusion.	
	\end{proof}
	
	\section{An application: \kalmans loop.}\label{section:KalmanLoop}
	
	\subsection{\kalmans loop.} K\'{a}lm\'{a}n has constructed a series of examples of loops of Legendrian positive torus knots non--contractible in the space $\Leg(\R^3)$, though contractible in $\Emb(\NS^1,\R^3)$ \cite{kalman}. Let us prove that \kalmans examples are non trivial even as loops of formal Legendrian embeddings; that is, in the space $\FLeg(\R^3)$. We will prove that they are not contractible for any choice of parametrization, thus they are not contractible as loops of unparametrized oriented knots.
	\begin{figure}[h]
		\centering
		\includegraphics[width=0.8\textwidth]{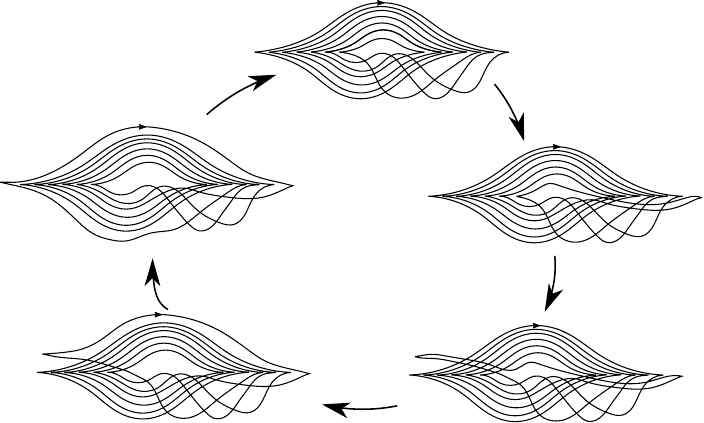}
		\caption{The loop in front projection $(p=3, q=7)$. Note that we have to cycle $2p=6$ times.
		\label{FrontKalman}}
	\end{figure}
	Consider a Legendrian positive $(p,q)$ torus knot, a loop is described in Figure \ref{FrontKalman}. The loop takes the $p$ strands of the knot to the cyclic rotation of them. This geometrically corresponds to a $2\pi/p$ rotation along the core of the defining torus. Let us consider $2p$ concatenations of this loop. Thus, it is generated by two full rotations along the core of the torus.	
	\subsubsection{Simplified position.}
	First, we will deform through formal loops the initial loop into a formal loop in a ``simplified'' position.\newline
	
	\textbf{Step 1.} Consider the contactomorphism $f(x,y,z)=(x/r,y/r,z/r^2)$ in the standard $(\R^3,\ker(dz-ydx))$. By using it, we assume that the defining torus for the loop has arbitrarily small meridional radius. Therefore, the knot is $C^1$--close to the core $\beta$ of the torus. We are not using the standard notion of $C^1$--closeness, but a weaker one. Ie we mean that a sequence of immersions $\hat \gamma_k$ is $C^1$-- close to an immersion $\tilde{\gamma}$ if for any $\varepsilon>0$ and for any point $t\in \NS^1$: for every $k\in \Z$ large enough there exists a point $\tau(t) \in \NS^1$ such that $|\hat \gamma_k(t)-\tilde{\gamma}(\tau(t))|\leq \varepsilon$ and $|\frac{\hat \gamma'_k(t)}{||\hat \gamma'_k(t)||}-\frac{\tilde \gamma'(\tau(t))}{||\tilde \gamma'(\tau(t))||}|\leq \varepsilon$, 
	
	Moreover, by further shrinking, the knot and the core $\beta$, they may be assumed to be arbitrarily close to $\ker dz$; ie $C^1$--close to their Lagrangian projections.
	
	\begin{figure}[h]
		\centering
		\subfloat[Front projection of the knot and the core (in red colour).]{
			\label{fig:KnotAndCore}
			\includegraphics[width=0.4\textwidth]{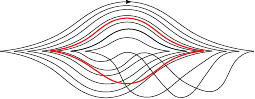}}
		\subfloat[Knot $C^1$--close to the core.]{
			\label{fig:KnotCloseCore}
			\includegraphics[width=0.4\textwidth]{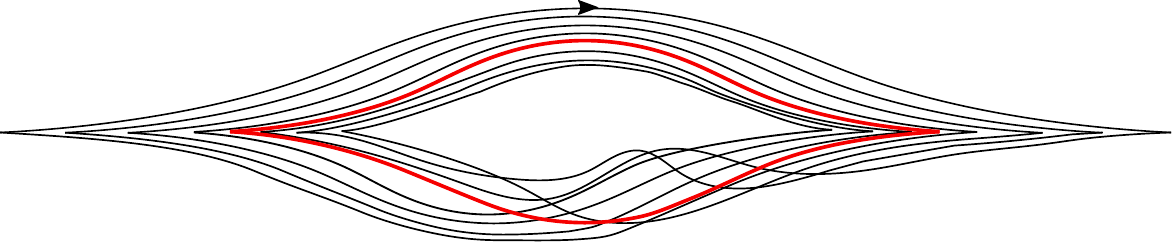}}
		\caption{$C^1$--approximation of the knot to the core, shown in the front projection.}
	\end{figure}
	
	\textbf{Step 2.} Denote $\gamma^\theta$ the initial loop of Legendrian embeddings. Understood as a formal loop, it is written as $(\gamma^\theta, F_{s}^{\theta})$, where $F_{s}^{\theta}=(\gamma^\theta)'$. Let us construct a $1$--parametric family of formal loops $(\gamma^{\theta,u}, F^{\theta,u}_s)$, $u\in[0,1]$, defined as follows
	\begin{itemize}
		\item [(i)] $\gamma^{\theta,u}=\gamma^\theta$,
		\item [(ii)] $F_s^{\theta,u}=(1-s)(\gamma^\theta)'+s((1-u)(\gamma^\theta)'+u\partial_y)$.
	\end{itemize}
	This is a family of formal loops, since $(\gamma^\theta)'$ is never a negative multiple of $\partial_y$.\footnote{We are assuming an orientation of the Legendrians. The argument with the opposite orientation runs in the same way by changing $\partial_y$ by $-\partial_y$.} To check that, note that $\beta$ is $C^1$--close to $\gamma^\theta$.
	
\begin{figure}[h]
		\centering
		\subfloat[Step 2.]{
			\label{fig:Step2}
			\includegraphics[width=0.3\textwidth]{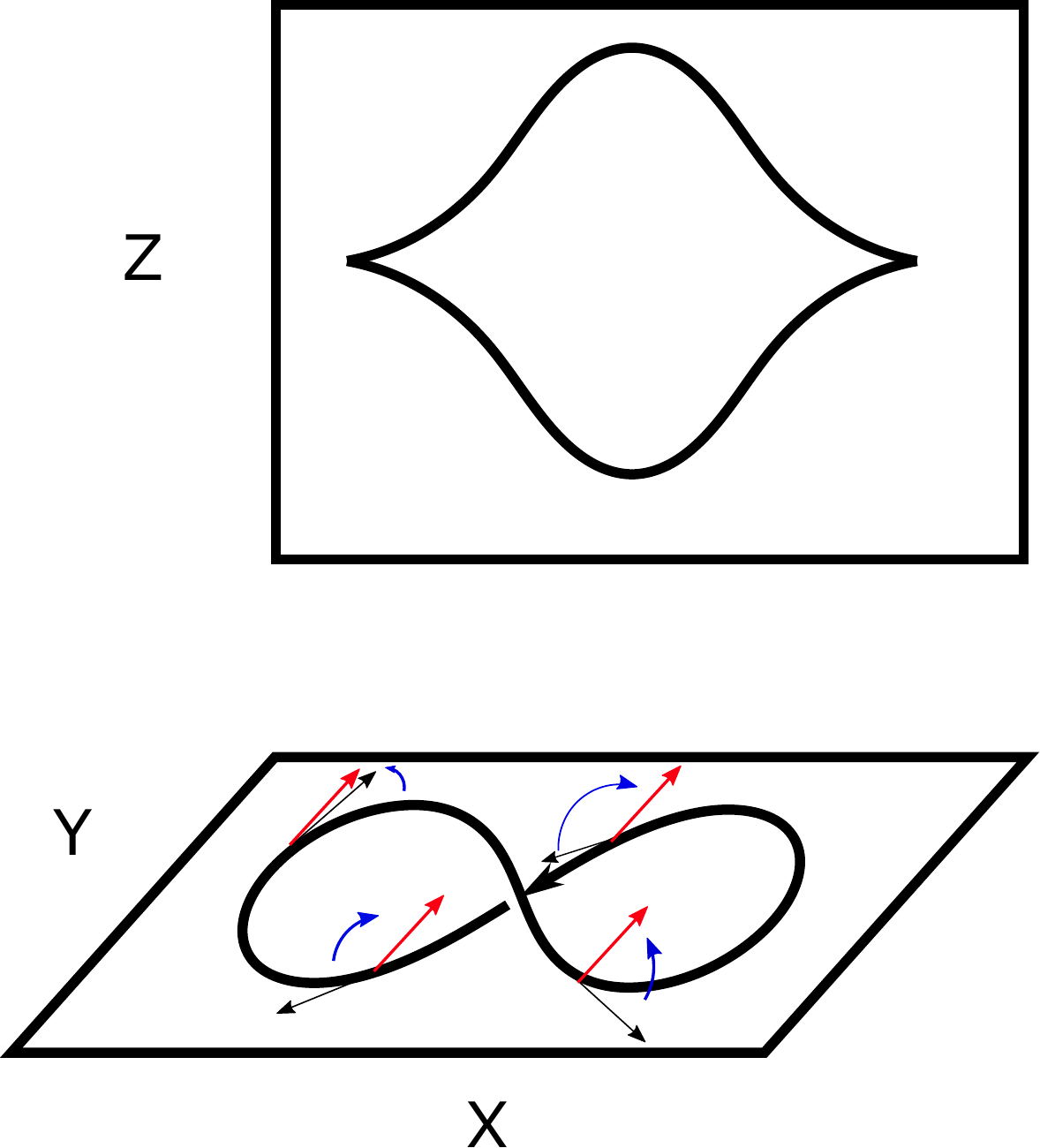}}
		\subfloat[Step 3.]{
			\label{fig:Step3}
			\includegraphics[width=0.3\textwidth]{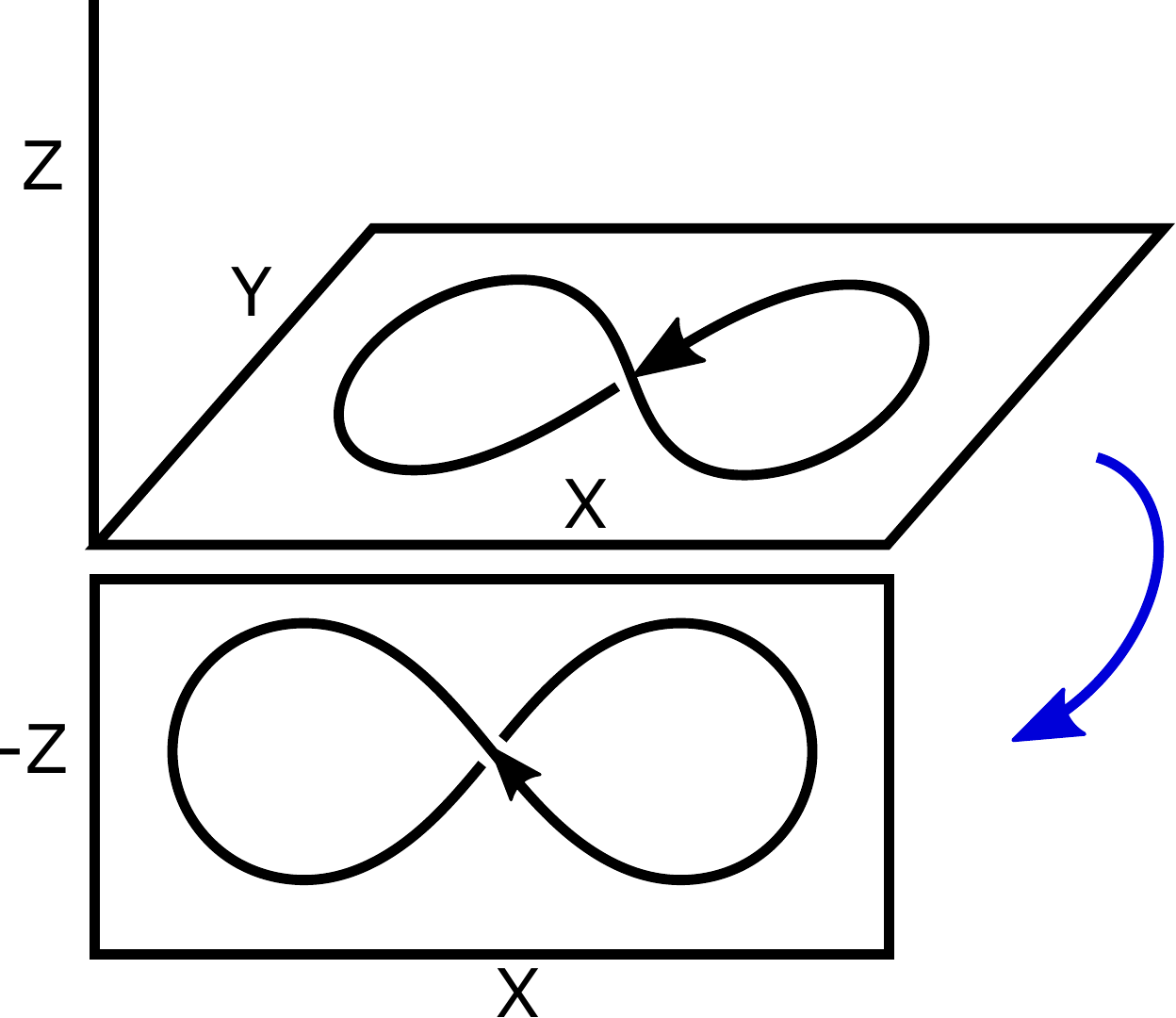}}
		\subfloat[Step 4.]{
			\label{fig:Step4}
			\includegraphics[width=0.3\textwidth]{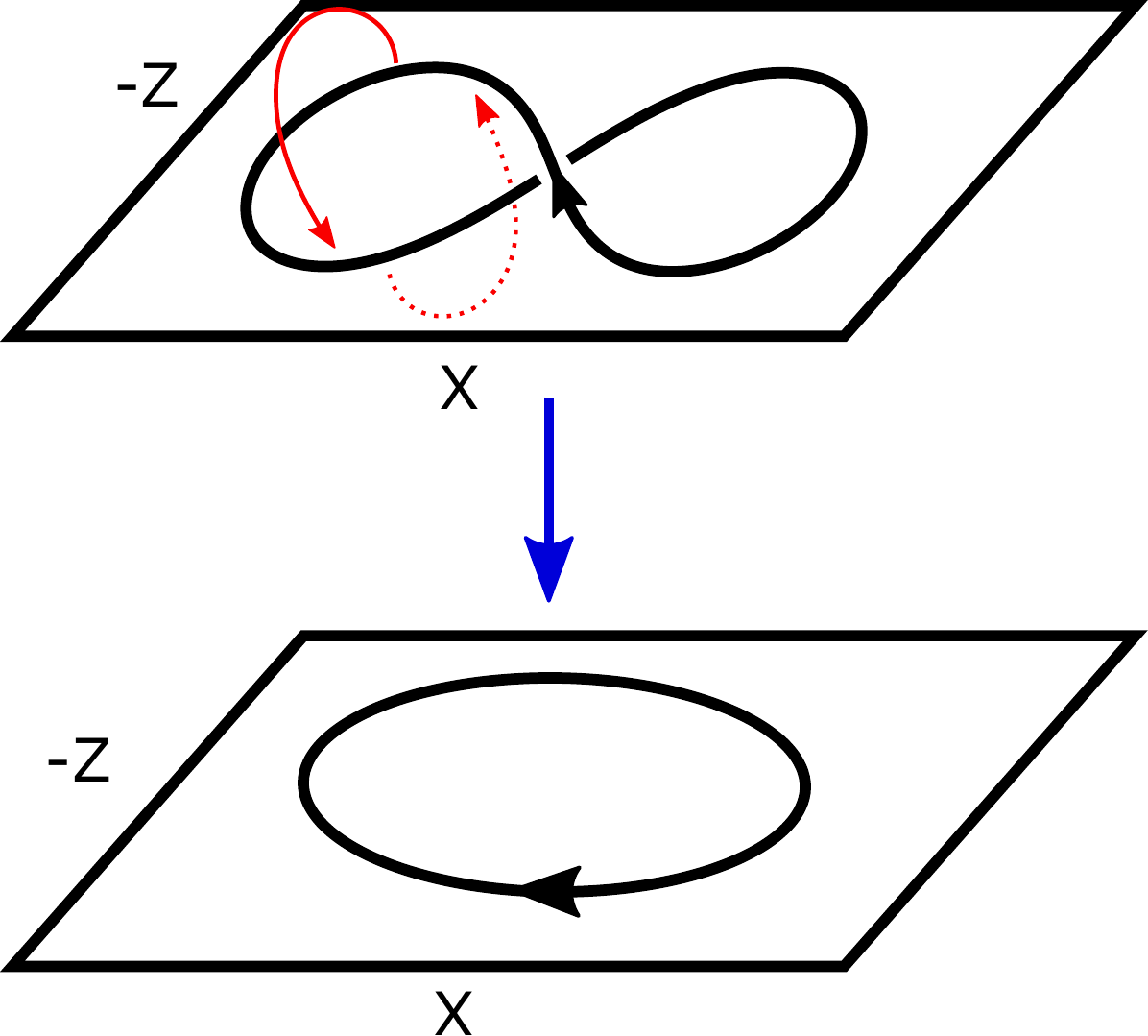}}
		\caption{Construction of the path of loops. We represent the moves of the core $\beta$.\label{SimplifiedKnot}}
	\end{figure}

	\textbf{Step 3.} We consider the family of rotations $\lbrace r_v\in \SO(3) \rbrace _{v\in [0,1]}$ taking the quadrant $XY$ to the quadrant $-ZX$, see Figure \ref{SimplifiedKnot}. Construct a family of formal loops $(\gamma^{\theta,u}, F^{\theta,u}_s), u\in[1,2]$, as follows
	
	\begin{itemize}
		\item[(i)] $\gamma^{\theta,u}=r_{u-1}\cdot\gamma^\theta$,
		\item [(ii)] $F_s^{\theta,u}=(1-s)(\gamma^{\theta,u})'+s\partial_y$.
	\end{itemize}
	
	Again, this is a family of formal loops because $(\gamma^{\theta,u})'$ is never a negative multiple of  $\partial_y$. Note that we are using that $\gamma^\theta$ is $C^1$--close to $\beta$.
	
	\textbf{Step 4.} Finally, we turn over the left lobe\footnote{See definition \ref{lobe}.} of the unknot core $r_1\cdot\beta$ by an isotopy defined as follows. Take polar coordinates $(r,\varphi)$ in the plane $YZ$ and $\varepsilon>0$ small enough. We define the isotopy as:
	
	\begin{equation*}
		f_u(x,r,\varphi)=\begin{cases}
			(x,r,\varphi+u\cdot \chi(x)\cdot\pi)& \text{if $x\leq\varepsilon$,} \\
			(x,r,\varphi)& \text{if $\varepsilon\leq x$} 
		\end{cases}
	\end{equation*}
	where $0\leq u\leq 1$ and $\chi:\R\rightarrow\R$ is a non--decreasing smooth function satisfying
	\begin{itemize}
		\item $\chi(x)=1$ for all $x\leq 0$,
		\item $\chi(x)=0$ for all $x\geq \varepsilon$.
	\end{itemize}
	
	We apply the isotopy to $r_1\cdot\gamma^\theta$, see Figure \ref{SimplifiedKnot}. Again, the derivative is never tangent to $-\partial_y$ and thus we can interpolate to $\partial_y$.
	
	We have proven that our initial loop of Legendrian embeddings is homotopic to the loop of formal Legendrian embeddings $(\tilde{\gamma}^\theta, \tilde{F}^{\theta}_{s})$ defined as follows: 
	\begin{itemize}
		\item[(i)] $\tilde{\gamma}^\theta$ is the loop of parametrized $(p,q)$ torus knots supported in the torus associated to the unknot contained in the plane $XZ$. The loop is obtained by a rotation of $4\pi$ radians of the standard $(p,q)$--embedding in the direction of the parallel of the supporting torus.
		\item[(ii)] $\tilde{F}_{s}^{\theta}=(1-s)(\tilde{\gamma}^\theta)'+s\partial_y$.
	\end{itemize}
	
	\subsection{Set of parametrizations of the family of loops.}
	As an outcome of the previous discussion, we may assume that our formal Legendrian parametrized $(p,q)$ torus knot can be written as $(\gamma, F_s)$, where 
	\begin{itemize}
		\item $\gamma(t)=\left( {\begin{array}{c}
			(\cos(2\pi p t)+2)\cos(2\pi qt)\\
			\sin(2\pi p t)\\
			(\cos(2\pi p t)+2)\sin(2\pi qt)
			\end{array} }\right)$,
		\item $F_s(t)=s\partial_y+(1-s)(\gamma)'(t).$
	\end{itemize}
	
	One particular parametrization of the loop can be written as $( \gamma^{\theta}, F^{\theta}_s)$, where 
	\begin{itemize}
		\item $\gamma^{\theta}(t)=\left( {\begin{array}{ccc}
			\cos(4\pi\theta) & 0 & -\sin(4\pi\theta)\\
			0 & 1 & 0 \\
			\sin(4\pi\theta) & 0 & \cos(4\pi\theta) \end{array} } \right)\left( {\begin{array}{c}
			(\cos(2\pi p t)+2)\cos(2\pi qt)\\
			\sin(2\pi p t)\\
			(\cos(2\pi p t)+2)\sin(2\pi qt)
			\end{array} }\right)$,
		\item $F^{\theta}_s(t)=s\partial_y+(1-s)(\gamma^{\theta})'(t)$.
	\end{itemize}
	
	We will show that any possible parametrization of the loop gives raise to a non--trivial loop of parametrized formal Legendrian embeddings. Up to homotopy, the possible parametrizations of the formal Legendrian loop are given by:
	\begin{itemize}
		\item$\gamma^{\theta,k}(t)=\gamma^{\theta}(t+k\theta),$
		\item $F^{\theta,k}_{s} (t)=F^{\theta}_{s} (t+k\theta)=s\partial_y+(1-s)(\gamma^{\theta,k})'(t),$
	\end{itemize}
	where $k\in\Z$. This is because $\pi_1 (\Diff^+(\NS^1))=\pi_1 (\SO(2))=\Z$.\footnote{Remember that the Legendrian knots are oriented.}
	
	We will prove the following statement.
	
	\begin{proposition}\label{nontrivial}
		The loop of formal Legendrian embeddings $(\gamma^{\theta,k},F_{s}^{\theta,k})$ is non trivial for any $k\in\Z$.
	\end{proposition}
	
	This proves that the loop is non trivial as a loop of non parametrized formal Legendrian knots. 
	
	\subsection{Proof of Proposition \ref{nontrivial}.}
	It follows from the previous discussion that the loop $\gamma^{\theta,k}(t)$ of smooth embeddings lies in $\SO(4)\subset\Emb_{p,q}(\NS^1,\NS^3)$, ie $\gamma^{\theta,k}(t)=A_{\theta,k}\gamma(t)$, where $A_{\theta,k}\in\SO(4)$. More specifically, on $\NS^3(\sqrt2)$, the $\sqrt{2}$ radius sphere in $\C^2$, we have \[ \gamma^{\theta,0}(t)=
	\left( {\begin{array}{cc}
		1  & 0 \\
		0 & e^{4\pi i \theta} \\
		\end{array} } \right)\left( {\begin{array}{c}
		e^{2\pi i p t} \\
		e^{2\pi i q t}\\
		\end{array} }\right). \]
	Thus, in these coordinates the other parametrizations are given by \[ \gamma^{\theta,k}(t)=\gamma^{\theta,0}(t+k\theta)=
	\left( {\begin{array}{cc}
		e^{2\pi ipk\theta} & 0 \\
		0 & e^{2\pi i(2+qk) \theta} \\
		\end{array} } \right)\left( {\begin{array}{c}
		e^{2\pi i p t} \\
		e^{2\pi i q t}\\
		\end{array} }\right). \]
	
	By Theorem \ref{thm:HatcherTorusKnots} the parametrized loop $\gamma^{\theta,k}$ is trivial in $\Emb_{p,q}(\NS^1,\NS^3)$ if and only if $2+k(p+q)$ is even. From now on we will assume that this is the case. Thus, there is a family $\lbrace \tilde{A}_{(r, \theta),k}\rbrace_{(r,\theta)\in\mathbb{D}}$ such that $\tilde{A}_{(1, \theta),k}=A_{\theta,k}$. Since $\pi_2(\SO(4))=0$, the disk $\tilde{A}_{(r, \theta),k}$ is unique up to homotopy fixing the boundary and the same holds for the disk $\tilde{A}_{(r, \theta),k}\cdot\gamma$ in $\Emb_{p,q}(\NS^1,\NS^3)$. 
	
	By Theorem \ref{pi1FLeg}, we have the following exact sequence:
	
	\begin{equation} \label{eq:shortseq}
	\xymatrix@M=10pt{
		& 0 \ar[r] & \pi_1(\mathcal{F})  \ar[r]^<<<<<{m_1} & \pi_1(\FLeg(\R^3))   \ar[r]^{m_2} & \pi_1(\FLegO(\R^3)) \ar[r]& 0}
	\end{equation}
	
	Thus, in order to prove that our loop $\ell_k=[(\gamma^{\theta,k}, F^{\theta,k}_s)]\in \pi_1(\FLeg(\R^3))$ is non trivial we distinguish two cases:
	
	\subsubsection{Case 1. $k\neq0$.}
	
	We claim that $m_2(\ell_k)\neq 0$, ie $[\gamma^{\theta,k}]\neq 0\in \pi_1 (\Emb_{p,q}(\NS^1,\R^3))$\footnote{Observe that $\Rot_L (\gamma^{\theta,k},F^{\theta,k}_{s})=0.$}.
	
	Recall from Corollary \ref{cor:FundamentalGroupKnots} that $\pi_1 (\Emb_{p,q}  (\NS^1,\R^3))\cong \pi_1 (\Stereo_{p,q})$ and that we have an exact sequence 
	\begin{equation*}
	\xymatrix@M=10pt{
		0\ar[r] & G_{p,q}\ar[r] & \pi_1 (\Stereo_{p,q}) \ar[r] & \pi_1 (\Emb_{p,q}(\NS^1,\NS^3)) \ar[r]& 0} 
	\end{equation*} 
	
	Thus, we must show that $[(A_{\theta,k} \gamma,\infty)]\in\pi_1 (\Stereo_{p,q})$ is non trivial. Ie the family of loops that is trivial by hypothesis in $\pi_1(\Emb_{p,q}(\NS^1,\NS^3))$ does not admit a capping disk whose evaluation map avoids $\infty\in\NS^3$. We check it by composing with the $1$--parametric family of loops $\tilde A_{(r, \theta),k}^{-1}\in\SO(4)$, we obtain a $1$--parametric family of loops $(\tilde A_{(r, \theta),k}^{-1}A_{\theta,k} \gamma,\tilde A_{(r, \theta),k}^{-1}\infty)$, $r\in[0,1]$. For $r=0$ we obtain the initial loop and for $r=1$ we obtain the loop  $(\gamma, A_{\theta,k}^{-1}\infty)$. Thus, these two loops represent the same element of $\pi_1(\Stereo_{p,q})$. Moreover, $[(\gamma,A_{\theta,k}^{-1}\infty)]$ can be lifted to $G_{p,q}$, since it lies on the fiber defined by the element $[\gamma]\in\pi_1(\Emb_{p,q}(\NS^1,\NS^3))$. So, we are reduced to check whether $[A_{\theta,k}^{-1}\infty]\in G_{p,q}$ is trivial. The knot group of the $(p,q)$ torus knot is $G_{p,q}= \langle a,b :a^p=b^q\rangle$. Thus, $[A_{\theta,k}^{-1}\infty]=b^{pk}\neq0$, since $b$ is a non torsion element of $G_{p,q}$.

\begin{figure}[h!]
	\centering
	\subfloat{
		\label{fig:02Torus}
		\includegraphics[scale=0.3]{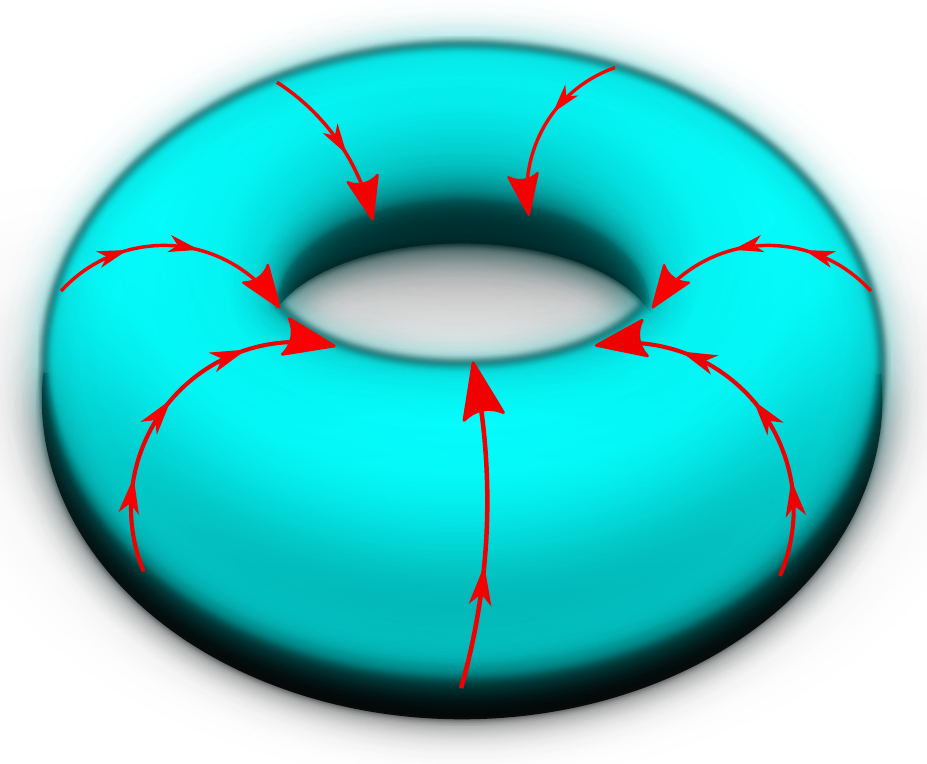}}
	\subfloat{
		\label{fig:02TorusKnot}
		\includegraphics[scale=0.20]{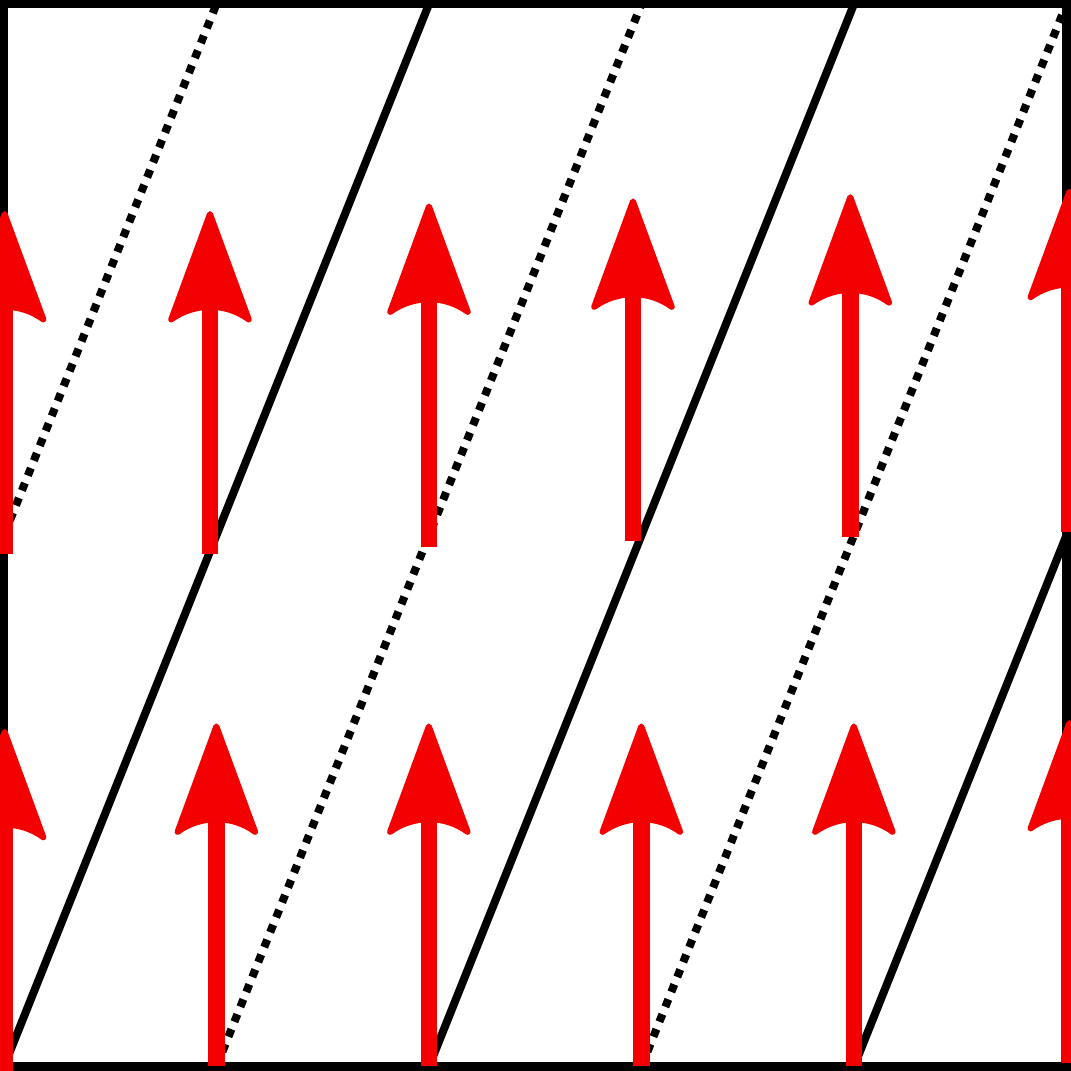}}
	\caption{Visualization of the loop $\gamma^{\theta,k}$ for a $(5,2)$ torus knot $(k\neq 0)$.}
	\label{fig:02Movement}
\end{figure}
	
	\subsubsection{Case 2. $k=0$.}
	Since $m_2(\ell_0)\in\pi_1(\FLegO(\R^3))$ is zero, there exists $A\in\pi_1(\mathcal{F})$ such that $m_1(A)=\ell_0$. We are going to geometrically check that $A\neq 0$ and therefore, by the injectivity of $m_1$ provided by the sequence (\ref{eq:shortseq}), the non triviality of $\ell_0$ follows.
	
	Write $(\gamma^\theta,F^{\theta}_s)=(\gamma^{\theta,0}, F^{\theta,0}_s)$. Note that in $(\R^3(x,y,z),\ker(dz-ydx))$ the parametrized loop is written as:
	\begin{itemize}
		\item $\gamma^{\theta}(t)=B_{\theta}\gamma(t)=
		\left( {\begin{array}{ccc}
			\cos(4\pi\theta) & 0 & -\sin(4\pi\theta)\\
			0 & 1 & 0 \\
			\sin(4\pi\theta) & 0 & \cos(4\pi\theta) \end{array} } \right)\left( {\begin{array}{c}
			(\cos(2\pi p t)+2)\cos(2\pi qt)\\
			\sin(2\pi p t)\\
			(\cos(2\pi p t)+2)\sin(2\pi qt)
			\end{array} }\right)$,
		\item $F^{\theta}_s(t)=s\partial_y+(1-s)(\gamma^{\theta})'(t)$.
	\end{itemize}

	\begin{figure}[h!]
		\centering
		\subfloat{
			\label{fig:20Torus}
			\includegraphics[scale=0.3]{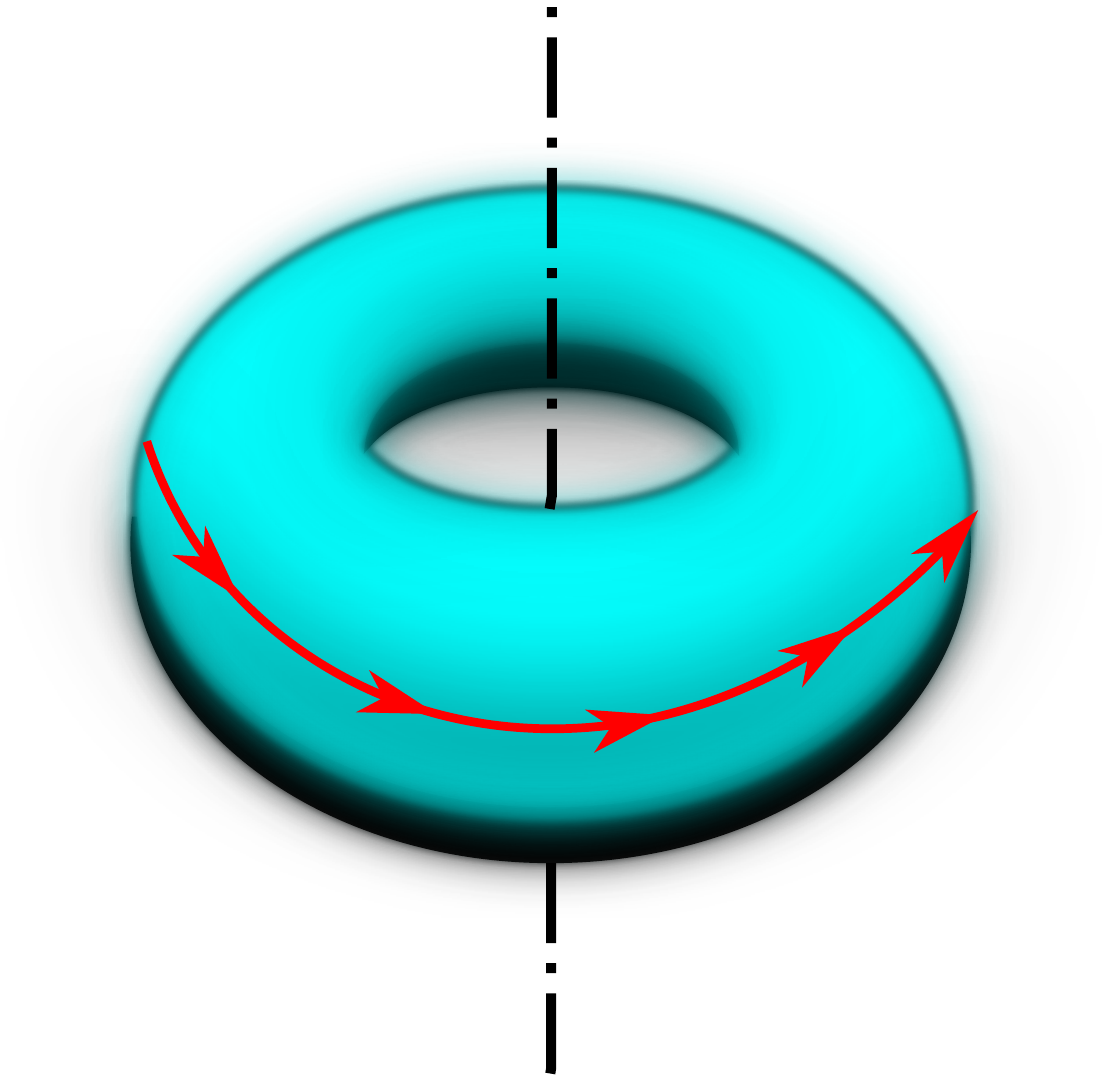}}
		\subfloat{
			\label{fig:20TorusKnot}
			\includegraphics[scale=0.20]{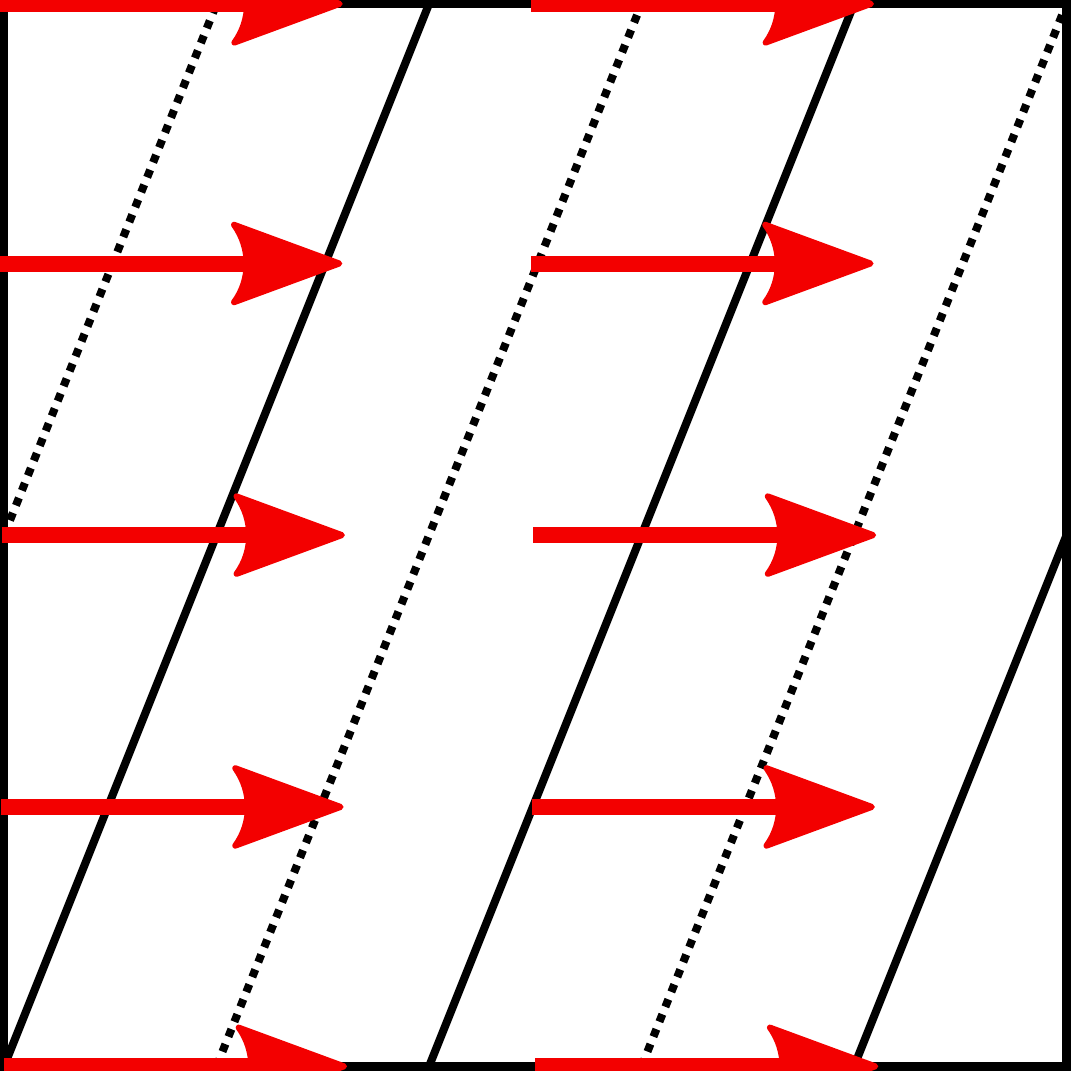}}
		\caption{Visualization of the loop $\gamma^{\theta,0}$ for a $(5,2)$ torus knot $(k=0)$.}
		\label{fig:20Movement}
	\end{figure}
	
	Moreover, $\gamma^{\theta}$ is bounded by $\tilde{\gamma}^{(r,\theta)}(t)=\tilde{B}_{(r,\theta)}\gamma(t)$, where $\tilde{B}_{(r,\theta)}\in\SO(3)$, such that $\tilde{B}_{(1,\theta)}=B_\theta$ and is homotopic to $\tilde{A}_{(r,\theta)}$ inside $\SO(4)$. Thus, we have a disk $\SD(r, \theta)=\{(\tilde{B}_{r,\theta} \gamma,\partial_y)\}$ in $\FLegO(\R^3)$ that bounds $m_2 (\ell_0)$. We try to lift it to a disk in $\FLeg(\R^3)$ that bounds $\ell_0$. There is no homotopical obstruction to lifting it for the punctured disk $\SD(r, \theta)|_{r>0}$. Therefore, the homotopy obstruction is represented by an element of $\pi_1(\mathcal F)$. So we obtain a loop $(\gamma_0,\tilde{F}^{\theta}_s)$ over the fiber  of $(\gamma_0,\partial_y)$ where $\gamma_0=B_{0}\gamma$. It follows by construction that $A=[(\gamma_0,\tilde{F}^{\theta}_s)]$. Let us perform the computation.
	\begin{equation} \label{eq:chunga}
		\tilde{F}^{\theta}_{s}=\begin{cases}
			\tilde{B}_{2s,\theta} (\gamma') & \text{if } 0\leq s\leq \frac{1}{2}, \\
			F^{\theta}_{2s-1}=(2s-1)\partial_y + (2-2s)B_\theta (\gamma') & \text{if } \frac{1}{2}\leq s\leq 1.
		\end{cases}
	\end{equation}  
	Observe that $\tilde{F}^{\theta}_{0}=\gamma'_{0}$ and $\tilde{F}^{\theta}_{1}=\partial_y$. Thus, we can understand $\tilde{F}^{\theta}_{s}$ as a map $\tilde{F}:\NS^2\rightarrow\Maps(\NS^1,\NS^2)$. It follows that $A=[\tilde{F}]\in\pi_1(\mathcal{F})\cong\pi_2(\Maps(\NS^1,\NS^2))\cong\pi_2(\NS^2)\oplus\pi_2(\Omega_p (\NS^2))$, that is the fundamental group of the fiber. We already computed this group in Theorem \ref{pi1FLeg}. Moreover, we also showed that the morphism to $\pi_1(\FLeg(\R^3))$ induced by the inclusion is injective.

	To conclude the proof we must check that $[\tilde{F}]\neq0$. In order to see this, we will verify that $\deg(\tilde{F}^{\theta}_{s}(0):\NS^2\rightarrow\NS^2)$ is nonzero. This degree is the first coordinate  of $[\tilde{F}] \in \Z\oplus\Z \cong \pi_1\left(\mathcal{F}\right)\cong\pi_2(\NS^2)\oplus\pi_2(\Omega_p(\NS^2))$.
	
\begin{figure}[h]
		\centering
		\subfloat{
			\includegraphics[scale=0.15]{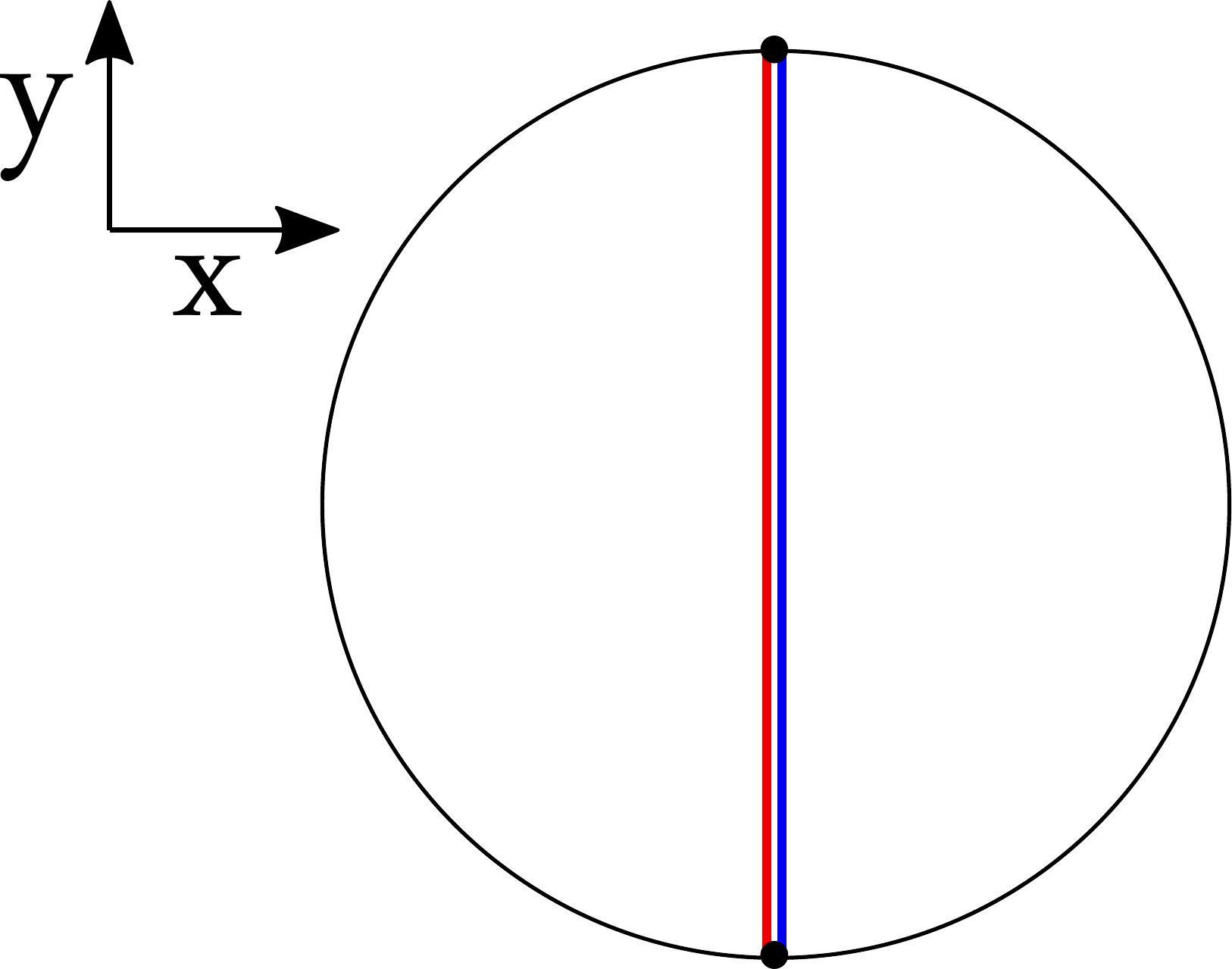}}
		\subfloat{
			\includegraphics[scale=0.15]{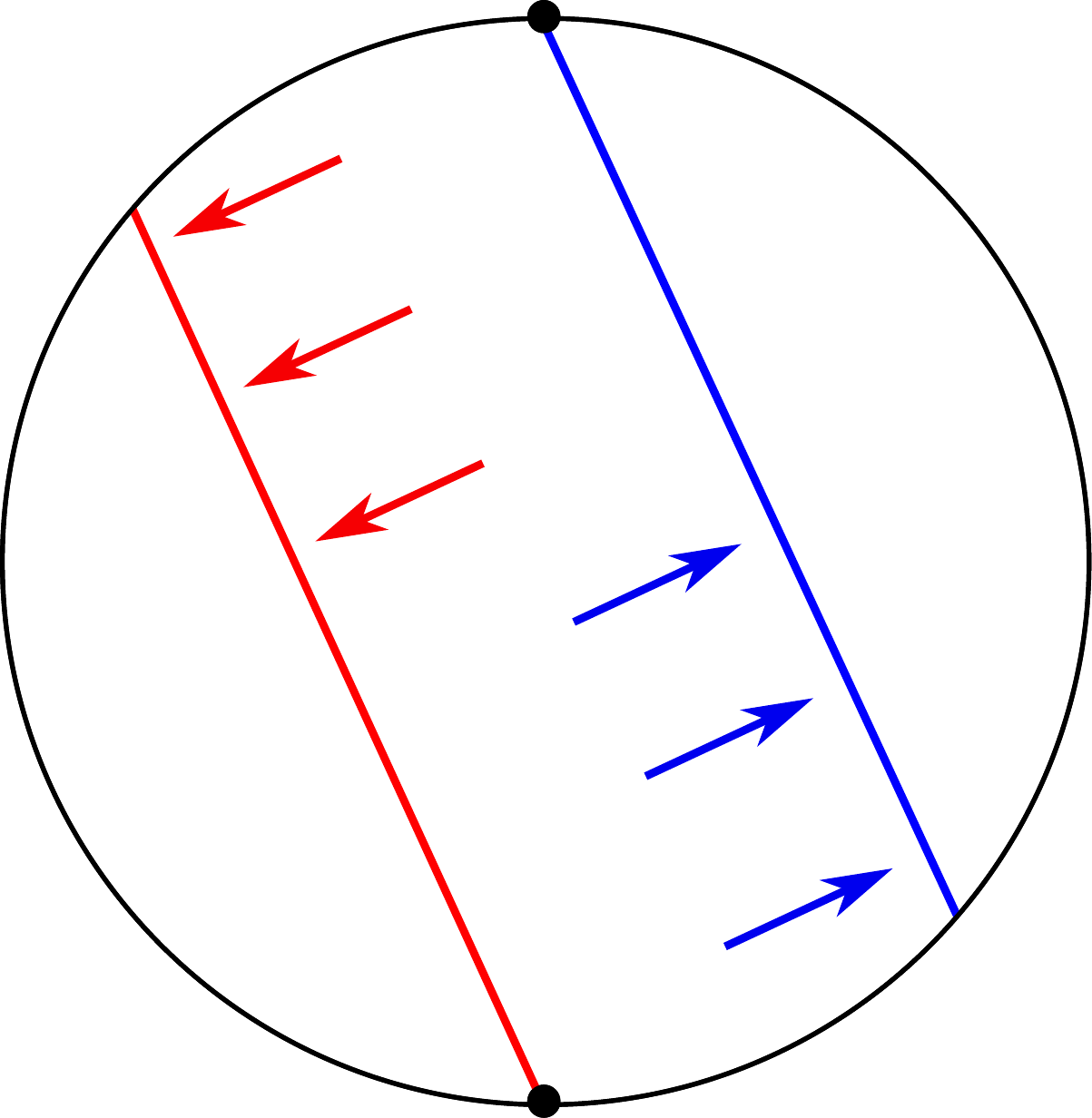}}
		\subfloat{
			\includegraphics[scale=0.15]{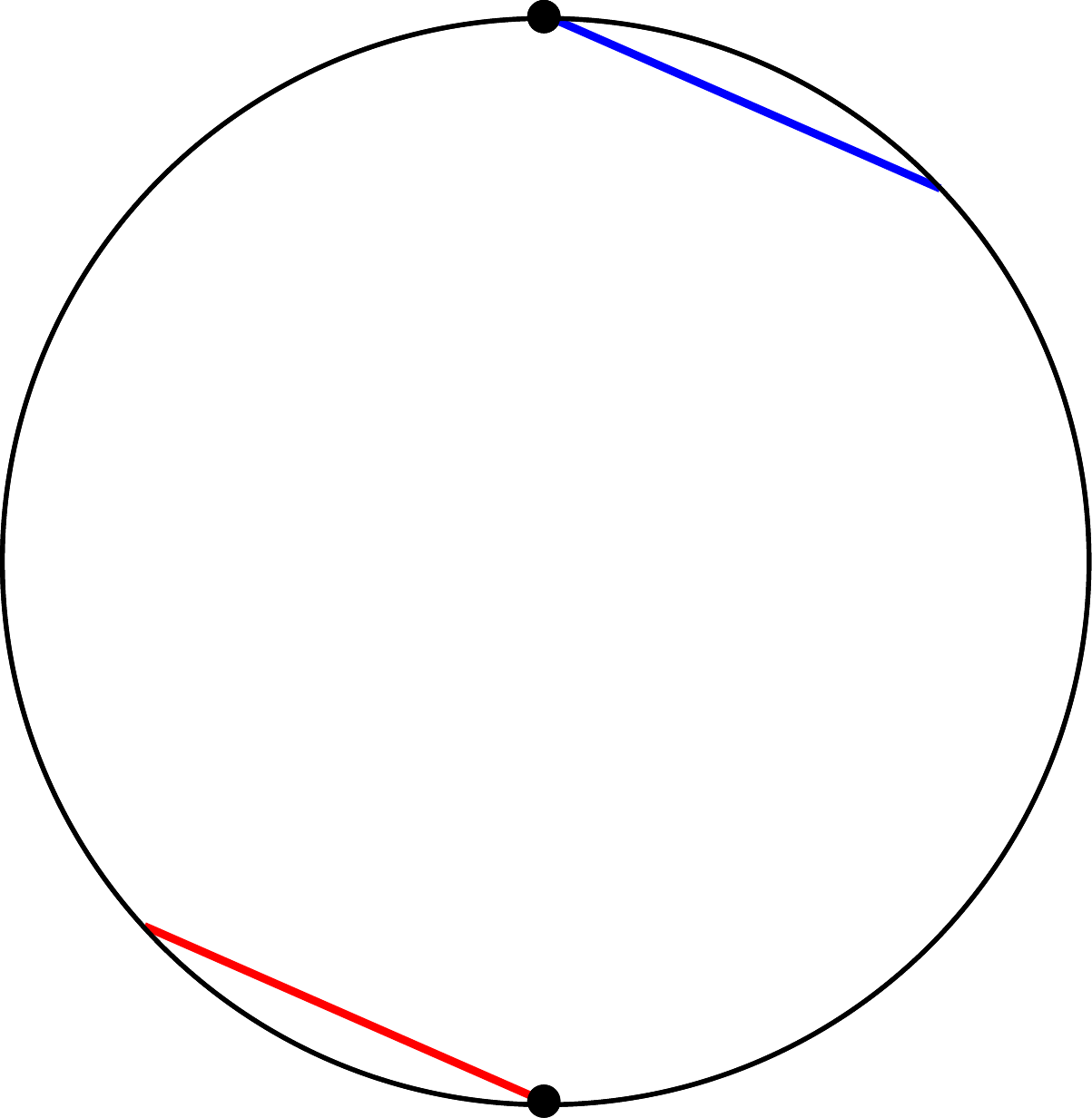}}
		\subfloat{
			\includegraphics[scale=0.15]{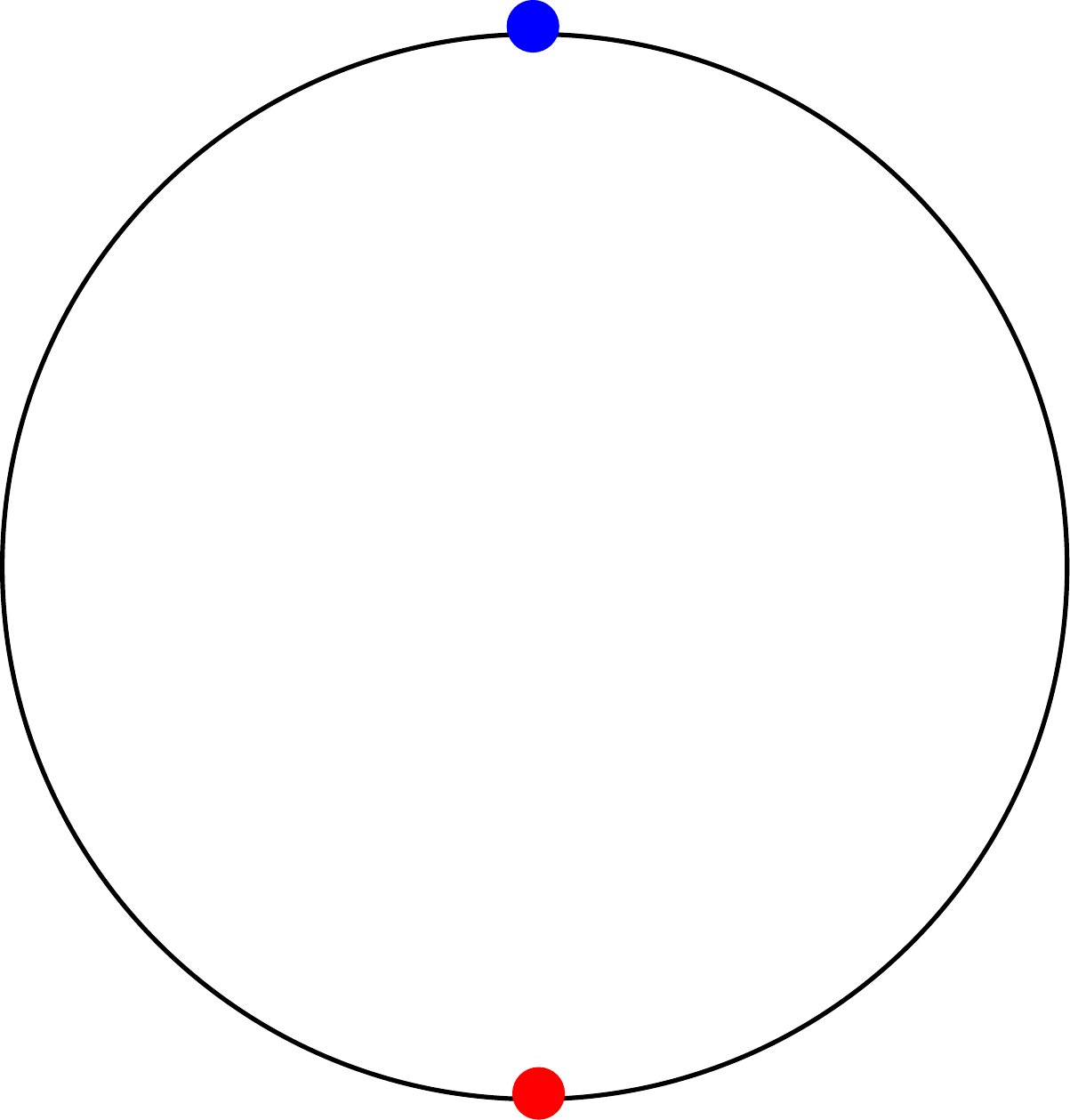}}
		\caption{Explicit construction of the capping disk.}
		\label{fig:capping}
	\end{figure}

	We give an explicit description of $\{\tilde{B}_{r,\theta}\}$ in $\SO(3)$. Identify $\SO(3)$ with $\mathbb{RP}^3$ in the usual way. Ie understand $\mathbb{RP}^3$ as the $3$--ball of radius $\pi$ with its boundary points identified via the antipodal map. Then a point $p=(p_1,p_2,p_3)$ in the described $3$--ball corresponds in $\SO(3)$ to the rotation of angle $\sqrt{p_1^2+p_2^2+p_3^2}$ in $\R^3(x,y,z)$ around the axis described by its position vector $p$. In these coordinates, see the left drawing in Figure \ref{fig:capping}, the loop $B_\theta$ is given by
	\begin{equation*}
		B_\theta=\begin{cases}
			(0, 4\pi\theta,0)& \text{if } 0\leq \theta\leq \frac{1}{4}, \\
			(0,-2\pi+4\pi\theta,0) & \text{if } \frac{1}{4}\leq \theta\leq \frac{3}{4}, \\
			(0,-4\pi+4\pi\theta,0) & \text{if } \frac{3}{4}\leq \theta\leq1.
		\end{cases}
	\end{equation*}
	
	Define the disk $\{\tilde{B}_{r,\theta}\}$ as the intersection of the plane $\{ z =0 \}$ with the $\pi$--radius ball. It produces an $\mathbb{RP}^2 =  \mathbb{RP}^3 \bigcap \{ z =0 \} \subseteq \mathbb{RP}^3$. We have $\hat{B} =\{B_{\theta} : \theta \in \NS^1\} \subseteq \mathbb{RP}^2$. We obtain $\mathbb{RP}^2 \setminus \hat{B}$ is an embedded $2$-disk. See Figure \ref{fig:capping}.
	
	We have that $\gamma'(0) = 2 \pi p \partial_y+ 6\pi q \partial_z$. Substituting in equation (\ref{eq:chunga}) for $t=0$, we get
\begin{equation*} 
		\tilde{F}^{\theta}_{s}(0)=\begin{cases}
			\tilde{B}_{2s,\theta} (2 \pi p \partial_y+ 6\pi q \partial_z) & \text{if } 0\leq s\leq \frac{1}{2}, \\
			F^{\theta}_{2s-1}=(2s-1)\partial_y + (2-2s)B_\theta (2 \pi p \partial_y+ 6\pi q \partial_z) & \text{if } \frac{1}{2}\leq s\leq 1.
		\end{cases}
	\end{equation*} 	 
	
	 Moreover, the maps, $u\in [0,1]$, 
\begin{equation*} 
		G(s,\theta,u)=\begin{cases}
			\tilde{B}_{2s,\theta} ((1-u)(2 \pi p \partial_y+ 6\pi q \partial_z)+ u\partial_z) & \text{if } 0\leq s\leq \frac{1}{2}, \\
			F^{\theta}_{2s-1}=(2s-1)\partial_y + (2-2s)B_\theta ((1-u)(2 \pi p \partial_y+ 6\pi q \partial_z)+ u\partial_z) & \text{if } \frac{1}{2}\leq s\leq 1
		\end{cases}
	\end{equation*} 	 
	 are always non zero. Thus, the map $\tilde{F}_{s}^{\theta}(0)=G(s,\theta,0)$ is homotopic to
	\begin{equation*}
		G(s,\theta)= G(s,\theta,1)=\begin{cases}
			\tilde{B}_{2s,\theta} (\partial_z) & \text{if } 0\leq s\leq \frac{1}{2}, \\
			(2s-1)\partial_y + (2-2s)B_\theta(\partial_z) & \text{if } \frac{1}{2}\leq s\leq 1.
		\end{cases}
	\end{equation*} 
	In order to compute the degree of $G(s,\theta)$, we write $G(s,\theta)=(g_x(s,\theta),g_y(s,\theta),g_z(s,\theta))\in\NS^2(\partial_x,\partial_y,\partial_z)$ and we check that $\#G^{-1}(-\partial_y)=1$:
	\begin{itemize}
	\item   If $\frac12\leq s\leq 1$ then $G(s,\theta)$ is a linear combination with positive coefficients between $\partial_y$ and $B_\theta(\partial_z)\in\NS^1(\partial_x,\partial_z)\subseteq\NS^2(\partial_x,\partial_y,\partial_z)$ thus $g_y(s,\theta)\geq0$.
	
	\item If $0\leq s\leq \frac12$ then $G(s,\theta)$ is just a rotation around an axis in the $XY$--plane acting over $\partial_z$. The rotation of angle $\frac{\pi}{2}$ around the $X$--axis is the unique rotation that sends $\partial_z$ to $-\partial_y$. Thus, $G^{-1}(-\partial_y)=\{(s_0,\theta_0)\}$ where $(s_0,\theta_0)$ is the only point that satisfies that $\tilde{B}_{2 s_0, \theta_0}$ is the mentioned rotation.
	\end{itemize}
	
	The map $G$ is a local diffeomorphism in a neighborhood of the point $(s_0,\theta_0)$. Thus, $-\partial_y$ is a regular value for $G$ and  $|\deg(\tilde{F}^{\theta}_{s}(0))|=|\deg(G(s,\theta))|=1\neq0$.
	
	\section{Formal Horizontal Embeddings in $\R^4$.}\label{FormalHorizontal}
	In this Section, we denote by $\SD$ the standard Engel structure in $\R^4(x,y,z,w)$ given by $\SD=\ker(dy-zdx)\cap\ker(dz-wdx)$. Throughout the Section we fix the framing of $\SD$ given by the \em kernel \em of the Engel structure $\SW=\langle \partial_w \rangle$.
	
	\subsection{Formal Horizontal Embeddings in $\R^4$.}
	
	\begin{definition}
		An immersion $\gamma:\NS^1\rightarrow\R^4$ is said to be \em horizontal \em if $\gamma'(t)\in\SD_{\gamma(t)}$ for all $t\in\NS^1$. When $\gamma$ is an embedding, we say it is a \em horizontal embedding\em. 
	\end{definition}
	
	\begin{remark}
		The projection \[\pi_{G}: \R^4(x,y,z,w)\rightarrow \R^3 (x,z,w)\] is called the \em Geiges projection\em. Notice that the Geiges projection maps horizontal embeddings to Legendrian immersions in $(\R^3(x,z,w), \ker(dz-wdx))$. Moreover, since the front projection for $(\R^3(x,z,w),\ker(dz-wdx))$ is the Lagrangian projection for $(\R^3(x,y,z),\ker(dy-zdx))$ the Geiges projection $\gamma_G$ of a horizontal embedding $\gamma\in\Hor(\R^4)$ is a Legendrian immersion which satisfies the following area condition \[\int_{\gamma} zdx=0.\] The composition of the Geiges projection and the front projection is called the \em front Geiges projection\em. The front Geiges projection of $\gamma$ is denoted by $\gamma_{FG}$. The front projection of a Legendrian immersion $\beta$ is denoted by $\beta_F$ and the Lagrangian projection is denoted by $\beta_L$.
	\end{remark}
	
	\begin{definition}\label{lobe}
		Let $\gamma\in\Imm(\NS^1,\R^2)$. By a \em lobe of $\gamma$ \em we mean any segment of the curve that encloses a topological disk in $\R^2\backslash\gamma(\NS^1)$. 
	\end{definition}
	When we do not specify it explicitly, by a lobe of a horizontal embedding $\gamma$ we mean a lobe of $\gamma_{FG}$.
	
	\begin{definition}
		\item [(a)]	A \em formal horizontal immersion \em in $\R^4$ is a pair $(\gamma,F)$ such that: 
		\begin{itemize}
			\item [(i)] $\gamma:\NS^1\rightarrow\R^4$ is a smooth map.
			\item [(ii)] $F:\NS^1\rightarrow \gamma^*(T\R^4\backslash\{0\})$ satisfies $F(t)\in\SD_{\gamma(t)}$.
		\end{itemize}
		\item[(b)] A \em formal horizontal embedding \em in $\R^4$ is a pair $(\gamma,F_s)$ satisfying: 
		\begin{itemize}
			\item [(i)] $\gamma:\NS^1\rightarrow\R^4$ is an embedding.
			\item [(ii)] $F_s:\NS^1\rightarrow \gamma^*(T\R^4\backslash\{0\})$, is a $1$--parametric family, $s\in[0,1]$, such that $F_0=\gamma'$ and $F_1(t)\in\SD_{\gamma(t)}$.
		\end{itemize}
	\end{definition}
	
	Use $\langle \partial_w \rangle $ to identify $\SD\equiv\R^2$. From now on, we will understand the family $F_s:\NS^1\rightarrow\NS^3$ with $F_1:\NS^1\rightarrow\NS^1\equiv\NS^3\cap\R^2$. 
	
	Denote by $\HorImm(\R^4)$ the space of horizontal immersions and by $\Hor(\R^4)$ the space of horizontal embeddings in $\R^4$. All these spaces are endowed with the $C^r$--topology. Denote by $\FHorImm(\R^4)$ the space of formal horizontal immersions in $\R^4$ and by $\FHor(\R^4)$ the space of formal horizontal embeddings in $\R^4$. All these spaces are endowed with the $(C^r, C^{r-1})$--topology. These definitions make sense for immersions and embeddings of the interval. We define $\HorImm([0,1],\R^4)$, $\Hor([0,1],\R^4)$, $\FHorImm([0,1],\R^4)$ and $\FHor([0,1],\R^4)$ analogously.
	\begin{remark}
		In the same vein as in the Legendrian case (see Remark \ref{rem:FrechetImm}) we can equip the spaces $\HorImm(\R^4),\Hor(\R^4),\FHorImm(\R^4)$ and $\FHor(\R^4)$ with a structure of Banach manifolds.
		
		Like in the Legendrian case it is important to note that any $\gamma\in\HorImm^r(\R^4)$ that is transverse to the kernel $\mathcal{W}$ of the Engel structure admits a particular type of immersed chart $\varphi_{\gamma}$ which identifies via an immersed Engelmorphism a tubular neighborhood of the zero section of the jet space $J^2(\NS^1)$ (see \cite[Lemma 1]{PinoPresas}) with a tubular neighborhood of $\gamma(\NS^1)$. This construction provides local charts in the $C^r$--topology just by defining the map
		\begin{equation}\label{eq:ExponentialHorImm}
		\begin{array}{rccl}
		\mathcal{C}\colon & \mathcal{U}\subseteq \Maps^r(\NS^1, \R) \times \Maps^r(\NS^1, \R) & \longrightarrow &  \mathcal{V}\subseteq\HorImm^r(\R^4) \\
		& (f,g)& \longmapsto & \varphi_{\gamma}^{-1}\circ j^2(g\circ\tilde{f}).
		\end{array}
		\end{equation}
		This proves that the open subset of horizontal immersed curves that are transverse to the kernel of the Engel structure is a Banach manifold. 
	\end{remark}

	\begin{remark}\label{rem:RotationsHorizontal}
		In \cite{PinoPresas} it is proved that the inclusion $\HorImm(\R^4)\hookrightarrow\FHorImm(\R^4)$ is a weak homotopy equivalence. Thus, $\pi_0 (\HorImm(\R^4))\cong\Z$, $\pi_1 (\HorImm(\R^4))\cong\Z$ and $\pi_k (\HorImm(\R^4))=0$ for all $k\geq 2$. As in the Legendrian case the connected components of $\HorImm(\R^4)$ are classified by a \em rotation number\em. The rotation number of $\gamma\in\HorImm(\R^4)$ is $\Rot(\gamma)=\deg(\gamma':\NS^1\rightarrow\NS^1)$. In the same way the homotopy type of a loop $\gamma^\theta$ in $\HorImm(\R^4)$ is determined by the integer $\Rot_L (\gamma^\theta)=\deg(\theta\mapsto(\gamma^\theta)'(0))$ which is called \em rotation number of the loop\em. Moreover, it is known that horizontal embeddings in $\R^4$ are also classified by their rotation number. This result was proven by Adachi \cite{Ada} and Geiges \cite{GeigesLoops}. Note that these invariants are defined in the formal setting in the obvious way. 
	\end{remark}
	
	\subsection{The space $\FHor(\R^4)$.}
	Consider the natural fibration $\FHor(\R^4)\rightarrow\FHorO(\R^4)=\{(\gamma,F):\gamma\in\Emb(\NS^1,\R^4),F\in\Maps(\NS^1,\NS^1)\}$. Take $\gamma\in\Hor(\R^4)$ and fix $(\gamma,\gamma')\in\FHorO(\R^4)$ as base point to compute the homotopy of $\FHor(\R^4)$. The fiber over this point is $\F=\F_{(\gamma,\gamma')}=\Omega_{\gamma'} (\Maps(\NS^1,\NS^3))$. To compute some homotopy groups of $\F$ use the fibration $\Maps(\NS^1,\NS^3)\rightarrow\NS^3$ defined by the \em evaluation map. \em Observe that, since every element $[f]\in\pi_n (\NS^3)$ can be lifted to an element $[f_n]\in\pi_n (\Maps(\NS^1,\NS^3))$, defined as
	\begin{equation*}
		f_n(p)(t)=f(p), t\in\NS^1, p\in\NS^n,
	\end{equation*}
	all the diagonal maps in the associated exact sequence are zero and the associated short exact sequences are right split and, thus, split for $n>1$. Since $\NS^3$ is $2$--connected we conclude that $\pi_0 (\F)\cong\pi_1(\Maps(\NS^1,\NS^3))=0$, $\pi_1 (\F)\cong\pi_2(\Maps(\NS^1,\NS^3))\cong\pi_2 (\Omega_p (\NS^3))\cong\pi_3 (\NS^3)\cong\Z$ and $\pi_2 (\F)\cong\pi_3(\Maps(\NS^1,\NS^3))\cong\pi_3 (\NS^3)\oplus\pi_3 (\Omega_p (\NS^3))\cong\pi_3 (\NS^3)\oplus\pi_4 (\NS^3)\cong\Z\oplus\Z_2$.
	
	Furthermore, the space $\FHorO(\R^4)$ has the homotopy type of $\Emb(\NS^1,\R^4)\times\NS^1\times\Z$. Thus,  $\pi_0 (\FHorO(\R^4))\cong\Z$, $\pi_1 (\FHorO(\R^4))\cong\pi_1(\Emb(\NS^1,\R^4))\oplus\Z$ and $\pi_2 (\FHorO(\R^4))\cong\pi_2(\Emb(\NS^1,\R^4))\cong\Z[S_U]\oplus\pi_2(\LEmb(\R,\R^4))$ (see Lemma \ref{lem:HomotopyKnotsR4}).
	
	The exact sequence associated to the fibration $\FHor(\R^4)\rightarrow\FHorO(\R^4)$ takes the following form:
	\begin{displaymath}
	\xymatrix@M=10pt{
		& \cdots\ar[r] & \Z\oplus\pi_2  (\LEmb(\R,\R^4)) \ar[dll]\\
		\Z \ar[r] & \pi_1 (\FHor(\R^4)) \ar[r] & \pi_1 (\Emb(\NS^1,\R^4))\oplus \Z \ar[dll] \\
		0 \ar[r] & \pi_0 (\FHor(\R^4)) \ar[r] & \Z \ar[r] & 0 }
	\end{displaymath}
	
	In particular, this proves that \em formal horizontal embeddings are classified by their rotation number\em.
	
	We can state the next result concerning the fundamental group of each connected component of $\FHor(\R^4)$ :
	
	\begin{lemma}\label{Pi1FHor}
		The sequence
		\begin{displaymath}
		\xymatrix@M=10pt{
			0\ar[r] & \Z_2\ar[r] & \pi_1 (\FHor(\R^4)) \ar[r] & \pi_1 (\Emb(\NS^1,\R^4))\oplus\Z\ar[r] & 0 }
		\end{displaymath}
		is exact.
	\end{lemma}
	\begin{proof}
		It is sufficent to show that the image of the diagonal map $d:\pi_2 (\FHorO(\R^4))\rightarrow\pi_1 (\F)$ is $2\Z$. Observe that $d$ measures the obstruction to lifting an element $[(\gamma^z,F^{z}_{1})]\in\pi_2 (\FHorO(\R^4))$ to $\pi_2 (\FHor(\R^4))$. In other words, let $F_0=(\gamma^z)':\NS^2\times\NS^1\rightarrow\NS^3$, $(z,t)\mapsto F_0(z,t)=(\gamma^z)'(t)$, be the derivative map. The homomorphism $d$ measures the obstruction to find a homotopy between the derivative map $F_0=(\gamma^z)'$ and the map $F_1:\NS^2\times\NS^1\rightarrow\NS^3$, $(z,t)\mapsto F^{z}_{1}(t)$. Since the map $F_1$ is null--homotopic, we conclude that $d[(\gamma^z,F^{z}_{1})]=\deg(F_0)$.\footnote{We are arguing as in Lemma \ref{componentsFLeg}, but now the antipodal map in $\NS(\R^4)=\NS^3$ has degree $1$.}
		
		Recall that $\pi_2 (\FHorO(\R^4))\cong\pi_2(\Emb(\NS^1,\R^4))\cong\Z[S_U]\oplus\pi_2 (\LEmb(\R,\R^4))$, see Lemma \ref{lem:HomotopyKnotsR4}. We have that $d[S_U]=2$. Indeed, $d[S_U]$ is the degree of the map 
		\begin{center}
			$\begin{array}{rccl}
			F_{\gamma}=\frac{\gamma'}{||\gamma'||} \colon & \NS^2\times\NS^1 & \longrightarrow &  \NS^3 \\
			& ((\lambda_1,\lambda_2,\lambda_3),t)& \longmapsto & \cos(2\pi t)e_1+\sin(2\pi t)(\lambda_1 e_2 +\lambda_2 e_3 +\lambda_3 e_4),
			\end{array}$
		\end{center}
		where $\{e_1,e_2,e_3,e_4\}$ is the canonical basis of $\R^4\supseteq\NS^3$. It follows that \[F_{\gamma}^{-1}(e_4)=\lbrace q_1=((0,0,1),\frac{1}{4}), q_2= ((0,0,-1),\frac{3}{4})\rbrace.\] We want to compute the orientation class of the image of the tangent space into the tangent space. We understand $\NS^2 \times \NS^1$ as a submanifold of $\R^3 \times \NS^1$ and $\NS^3$ as a submanifold of $\R^4$. We extend the tangent spaces of those two submanifolds by adding the outward normal vectors to both spheres $T \NS^2  \times T\NS^1 \subset T (\NS^2  \oplus r \partial r) \times T\NS^1$ and $T \NS^3\subset T (\NS^3  \oplus r \partial r)$. We declare that the orientations in the submanifolds are induced from the orientations in $\R^3 \times \NS^1$ and $\R^4$. In order to compute the degree of the differential at the two points, we realize that $F_{\gamma}$ linearly extends to $\bar{F}_\gamma: \R^3 \times \NS^1 \to \R^4$. The degree of $F_{\gamma}$ is computed by that of $\bar{F}_{\gamma}$. We obtain	
		\begin{equation*}	
		\det(dF_{\gamma})=\det\begin{pmatrix}
		0   &        0     &        0     & -2\pi \sin(2\pi t) \\
		\sin(2\pi t)&        0     &        0     &  2\pi\cos(2\pi t)\lambda_1\\
		0   & \sin(2\pi t) &        0     & 2\pi\cos(2\pi t)\lambda_2\\
		0   &        0     & \sin(2\pi t) &  2\pi\cos(2\pi t)\lambda_3
		\end{pmatrix}=2\pi (\sin(2\pi t))^4
		\end{equation*}
		is non zero and positive in $(F_\gamma)^{-1}(e_4)$. Thus, $e_4$ is a regular value for $F_\gamma$ and $d[S_U]=\deg F_\gamma=2$.
				
		To conclude the proof we must check that any $2$--sphere $\gamma^{z}\in\Emb(\NS^1,\R^4)$ which comes from the space of long embeddings satisfies that $d[\gamma^z]\in 2\Z$. In order to see this observe that if $\gamma(t)=(x(t),y(t),z(t),w(t))^\intercal \in\Emb(\NS^1,\R^4)$ comes from a long embedding then we may assume the following: 
		\begin{itemize}
			\item [(i)] $\gamma((\frac{1}{4},\frac{3}{4}))\subseteq\Int((-1,1)^4)$,
			\item [(ii)] $\gamma(\frac{1}{4})=- e_1$, $\gamma(\frac{3}{4})=e_1$ and $\gamma'(\frac{1}{4})=\gamma'(\frac{3}{4})=e_1$,
			\item [(iii)] $\gamma(\NS^1\backslash(\frac{1}{4},\frac{3}{4}))$ is contained in the plane spanned by $e_1$ and $e_2$.
			\item[(iv)] $\gamma_{|\NS^1\backslash(\frac{1}{8},\frac{7}{8})}$ parametrizes the semicircle of radius $R$ and center $R\cdot e_2$ joining $\gamma(\frac{1}{8})=-R \cdot e_1+R\cdot e_2$ and $\gamma(\frac{7}{8})=R\cdot e_1+R\cdot e_2$  which passes trough  $\gamma(0)=2R\cdot e_2$,
		\end{itemize}
		where $R>>1$ is a constant independent of $\gamma$.
		
		Let $\gamma^{z}$ be any $2$--sphere which comes from the space of long embeddings. The associated map $\gamma'=(\gamma^z)'$ is homotopic to $\gamma^{h}(z,t)=\frac{ \gamma^z (t+h)-\gamma^z(t)}{|| \gamma^z (t+h)-\gamma^z(t)||}$, for any $h\in(0,1)$. Thus, $\gamma'$ is homotopic to $\gamma^{\tilde{h}}(z,t)=\frac{ \gamma^z (t+\tilde{h}(t))-\gamma^z(t)}{|| \gamma^z (t+\tilde{h}(t))-\gamma^z(t)||}$, where $\tilde{h}:[0,1]\rightarrow(0,1)$ is any continuous map which satisfies $\tilde{h}(0)=\tilde{h}(1)$. Consider the family of functions $\tilde{h}_\varepsilon$, $0<\varepsilon<\frac{1}{4}$, defined as
		
		\[\tilde{h}_\varepsilon (t)= \left\{ \begin{array}{lcc}
		\varepsilon &   \text{if}  & t\in[0,\frac{1}{4}-\varepsilon]\cup[\frac{7}{8}-\varepsilon,1]\\
		
		\frac{7}{8}-t & \text{if} &t\in[\frac{1}{4},\frac{7}{8}-\varepsilon]\\
		
		(1-\frac{t-\frac{1}{4}+\varepsilon}{\varepsilon})\varepsilon + \frac{t-\frac{1}{4}+\varepsilon}{\varepsilon} \frac{5}{8} &  \text{if} & t\in[\frac{1}{4}-\varepsilon,\frac{1}{4}]\\
		\end{array}
		\right.
		\]
		
		It follows that $e_3\notin\im(\gamma^{\tilde{h}_\varepsilon})$. Hence, $d[\gamma^z]=\deg(\gamma')=\deg(\gamma^{\tilde{h}_\varepsilon})=0\in2\Z$.
	\end{proof}
	
	\subsection{The $\Z_2$ factor in $\pi_1 (\FHor(\R^4))$.}\label{subseccionAreatwist}
	Let us explain the $\Z_2$ factor in the exact sequence described in Lemma \ref{Pi1FHor}. Observe that $\Z_2$ is just the \em kernel \em of the homomorphism 
	\begin{equation*}
		f:\pi_1 (\FHor(\R^4))\rightarrow\pi_1 (\FHorO(\R^4)),
	\end{equation*}
	induced by the fibration $\FHor(\R^4)\rightarrow\FHorO(\R^4)$ and that $\ker(f)\cong\pi_1(\F)/d(\pi_2(\FHorO))$ (see Lemma \ref{Pi1FHor}). Consider the unique isomorphism 
	\[\mu:\ker(f)\rightarrow\Z_2.\]
	\begin{definition}
		The $\mu$--invariant of a loop $(\gamma^\theta,F^{\theta}_{s})\in\ker(f)$ is $\mu(\gamma^\theta,F^{\theta}_{s})\in\Z_2$.
	\end{definition}
	
	Take a loop $(\gamma^\theta,F_{s}^{\theta})$ in $\FHor(\R^4)$ and assume that $(\gamma^\theta,F_{1}^{\theta})$ is trivial in $\FHorO(\R^4)$, ie $\gamma^\theta$ is trivial as a loop of smooth embeddings and $\Rot_L (\gamma^\theta,F_{s}^{\theta})=0$. Then, there is a disk $(\gamma^z,F_{1}^{z})$ in $\FHorO(\R^4)$, $z\in\D$, such that $(\gamma^{e^{2\pi i\theta}},F_{1}^{e^{2\pi i\theta}})=(\gamma^\theta,F_{1}^{\theta})$. Now we try to lift it to a disk $(\gamma^z,F_{s}^{z})$ in $\FHor(\R^4)$. This can be done for all the values of the parameter except for $z=0$, where we have $\NS^1$ different limits. So, we have a loop $(\gamma^0,\tilde{F}_{s}^{\theta})$ in the fiber of $(\gamma^0,F_{1}^{0})$. We define $\tilde{F}$ as
	\begin{center}
		$\begin{array}{rccl}
		\tilde{F}\colon & [0,1]\times\NS^1 & \longrightarrow &  \Maps(\NS^1,\NS^3) \\
		& (s,\theta)& \longmapsto & \tilde{F}(s,\theta)=\tilde{F}_{s}^{\theta}.
		\end{array}$
	\end{center}
	
	Since $\tilde{F}_{0}^{\theta}=(\gamma^0)'$ and $\tilde{F}_{1}^{\theta}=F^{0}_{1}$, we can understand $\tilde{F}$ as a map from $\NS^2$ to $\Maps(\NS^1,\NS^3)$. Ie $[\tilde{F}]\in\pi_1 (\F)\cong\pi_2(\Maps(\NS^1,\NS^3))\cong\pi_3(\NS^3)$. The $\mu$--invariant measures the \em degree mod $2$ \em of the associated map from $\NS^2\times\NS^1$ to $\NS^3$, ie $\mu(\gamma^\theta,F_{s}^{\theta})=\deg_2 (\tilde{F})$. This is because the elements in $\pi_1 (\F)\cong\pi_3(\NS^3)$ are classified by the degree but the construction of $\tilde{F}$ depends on the choice of disk, inside $\FHorO(\R^4)$, bounding $(\gamma^{\theta},F_{1}^{\theta})$. The difference between two disks is measured by the diagonal map $d:\pi_{2}(\FHorO(\R^4))\rightarrow\pi_2 (\F)$; ie $\tilde{F}$ is unique, up to homotopy, inside $\pi_1 (\F)/d(\pi_{2}(\FHorO(\R^4)))$. It follows from Lemma \ref{Pi1FHor} that the elements $\pi_1 (\F)/d(\pi_{2}(\FHorO(\R^4)))\cong\Z_2$ are classified by the degree mod $2$ of the associated map.
	
	\begin{figure}[h]
		\includegraphics[scale=1.5]{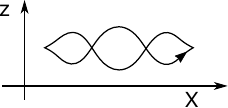}
		\centering
		\caption{The base point (a horizontal knot) of the are twist loop in the front Geiges projection.}\label{estabilizado}
	\end{figure}
	
\begin{figure}[h]
		\includegraphics[scale=0.5]{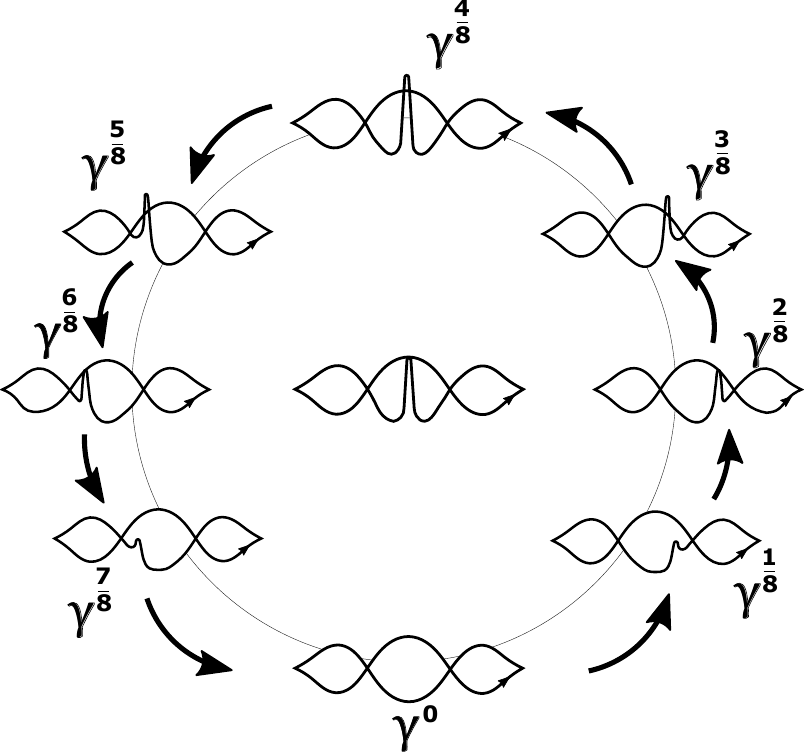}
		\centering
		\caption{The area twist loop in the front Geiges projection with the standard capping disk through horizontal immersions. Note that the unique strict horizontal immersion in the disk is the center of it. \label{AT}}
	\end{figure}

	Let us describe explicitly $\ker(f)$. The first element of $\ker(f)$ is just the trivial loop and the second one is given by the homotopy class of a loop $\gamma^\theta$ of horizontal embeddings that we name \em area twist loop \em and which is described in Figure \ref{AT} (we give a precise definition below). Observe that this loop is trivial as a loop of horizontal immersions. 
	
Let $\gamma$ be a horizontal immersion with exactly one generic self--intersection. We are going to construct a loop of horizontal embeddings that is supported on a very small neighborhood of the self--intersection. Informally speaking (it could be completely formalized by adapting to the horizontal case the Subsection \ref{localmodelsoflegendrians}, in particular we would need the equivalent of Lemma \ref{lem:SigmaBanach}), on an neighborhood of $\gamma$ inside the space of horizontal immersions, the subset of strict horizontal immersions has codimension $2$. Thus, there exists a disk $\mathbb D^2 \to \HorImm(\R^4)$ that intersects the strict immersion just at the center $\gamma$. There is a natural orientation on it and we take the loop $\partial \mathbb D^2$ with the induced orientation. Let us give a more hands--on approach.

Assume that the self--intersection times are $t_0, t_1$ and fix $t_1=0$. Now, select $\varepsilon,\delta>0$ such that $\delta <<\varepsilon<1$. By \cite[Lemma 1]{PinoPresas} we may assume that 
	\begin{itemize}
		\item $\gamma_{FG} (t)=(t-t_0,0)$ for $t\in\Op(\{t_0\})$ and 
		\item $\gamma_{FG} (t)=(t,-t^2)$ for $t\in(-2\varepsilon,2\varepsilon)$.
	\end{itemize}
	For each $(u,v)\in[-\delta,\delta]^2$ define $\gamma^{u,v}\in\HorImm(\R^4)$ in terms of the front Geiges projection $\gamma^{u,v}_{FG} (t)=(x^{u,v}(t),z^{u,v}(t))$ as follows:
	\begin{itemize}
		\item $\gamma^{0,0}=\gamma$,
		\item $\gamma^{u,v}_{FG}(t)=\gamma_{FG}(t)$ for $t\in\NS^1\backslash(-2\varepsilon,2\varepsilon)$,
		\item $x^{u,v}=x^{0,0}$,
		\item $z^{u,v}(t)$, $t\in(-2\varepsilon,2\varepsilon)$, is defined by the following conditions:
		\begin{itemize}
			\item [(i)] $z^{u,v}(t)=v-(t-u)^2$, $t\in[-\frac{\varepsilon}{2},\frac{\varepsilon}{2}]$,
			\item [(ii)] $z^{u,v}(t)=-t^2$, $t>|\varepsilon|$,
			\item [(iii)] $\int_{-\varepsilon}^{\varepsilon} z^{u,v} (t)dt=\int_{-\varepsilon}^{\varepsilon} z^{0,0} (t)dt$,
			\item [(iv)] $\int_{-\varepsilon}^{u} z^{u,v} (t)dt<\int_{-\varepsilon}^{0} z^{0,0} (t)dt$, for $0<u\leq\delta$,
			\item [(iv)] $\int_{-\varepsilon}^{u} z^{u,v} (t)dt>\int_{-\varepsilon}^{0} z^{0,0} (t)dt$, for $-\delta\leq u<0$.
		\end{itemize}
	\end{itemize}
	Observe that conditions (iii), (iv) and (v) imply $\gamma^{u,v}\in\Hor(\R^4)$ for $(u,v)\neq(0,0)$.
	
	\begin{definition}\label{def:AreaTwist}
		An \em area twist loop \em is any loop of horizontal embeddings which is homotopic to $\gamma^\theta:=\gamma^{\delta\cos(2\pi\theta),\delta\sin(2\pi\theta)}$.
	\end{definition}

	Observe that the notion of \em area twist loop \em is well--defined since the space of choices in the construction of $\gamma^{u,v}$ is contractible.
		
	The following lemma states that the area twist loop is actually the second element of $\ker(f)$.
	\begin{lemma}\label{AT.Z2.Invariant}
		Let $\gamma^\theta$ be an area twist loop. Then, $(\gamma^\theta,(\gamma^\theta)')\in\ker(f)$ and $\mu(\gamma^\theta,(\gamma^\theta)')=1$.  
	\end{lemma}
	\begin{proof}
		We prove the result for the area twist $\gamma^\theta$ described in Figure \ref{AT}. The general case is analogous.
		
		Consider the $1$--parametric family of formal horizontal embeddings $(\gamma^{\theta,u},F_{s}^{\theta,u})$, $u\in[0,1]$, defined by
		\item [(i)]$\gamma^{\theta,u}=\gamma^\theta$,
		\item [(ii)] $F_{s}^{\theta,u}=(1-s)(\gamma^\theta)'+s((1-u)(\gamma^\theta)'+u\partial_w)$.
		
		Note that, $F_{0}^{\theta,u}=(\gamma^\theta)'$, $u\in [0,1]$, and $(\gamma^\theta)'$ is never a positive multiple of $-\partial_w$. Hence, we can assume that the area twist loop in $\FHor(\R^4)$ is $( \gamma^\theta,F_{s}^{\theta,1})$. To compute the $\Z_2$--invariant of the area twist loop we are going to construct a disk $\{(\tilde{\gamma}^{r,\theta},\partial_w):r\in[0,1],\theta\in\NS^1\}$ in $\FHorO(\R^4)$ such that $\tilde{\gamma}^{1,\theta}=\gamma^\theta$. 
		
		For each $\lambda\in\R$ consider the $3$--plane $\Pi_\lambda\subseteq\R^4$ with equation $x=\lambda$. The plane $\Pi_\lambda$ intersects each embedding $\gamma^\theta$ in zero, one or two points. Consider the segment $G=\{\lambda: \# \Pi_\lambda\cap\gamma^\theta =1 \text{ or } 2 \}$ and write $\Pi_\lambda\cap\gamma^\theta =\{A_{\lambda}^{\theta},B_{\lambda}^{\theta}\}$ for all $\lambda\in G$. Denote by $\overline{A_{\lambda}^{\theta}B_{\lambda}^{\theta}}$ the unique segment (in $\R^4$) defined by the pair of points  $A_{\lambda}^{\theta}$, $B_{\lambda}^{\theta}$ . The union of all these segments allows us to construct an embedded $2$--disk $\D^\theta= \cup_{\lambda\in G} \overline{A_{\lambda}^{\theta}B_{\lambda}^{\theta}}$, whose boundary is $\gamma^\theta$. We call the foliation, in the so constructed $2$--disk, provided by the segments the ``ruling".
		
		Assume that the left cusp of $\gamma^{\theta}$ happens at time $t_0\in\NS^1$. Observe that $\gamma^{\theta}(t)=\gamma^{0}(t)$ for $t\in\Op(\{t_0\})$. Thus, by construction, the intersection  $\D^{\theta}\cap\Op(\{\gamma^{0}(t_0)\})=\D$ is independent of $\theta$. Fix $q_0\in\Int(\D)$ and carefully take polar coordinates $\D^{\theta}=\{p^{\rho,t}_{\theta}:\rho\in[0,1],t\in\NS^1\}$ such that
		\begin{itemize}
			\item $p^{1,t}_{\theta}=\gamma^{\theta}(t)$,
			\item $p^{0,t}_{\theta}=q_0$.
			\item for any $\rho$ and $\theta$ fixed, the loop $\{ p^{\rho,t}_{\theta} \}_{t \in \NS^1}$ intersects the ``ruling" provided by $\overline{A_{\lambda}^{\theta} B_{\lambda}^{\theta}}$ in zero, one or two points. In other words, the curve $\{ p^{\rho,t}_{\theta} \}_{t \in \NS^1}$ becomes tangent to the ruling just two times.
		\end{itemize}
		
		Moreover, fix $\varepsilon>0$ small enough and let $Q_0: \D(\varepsilon)\rightarrow\Int(\D)$ be an embedding such that $Q_0 (0)=q_0$. Assume that $Q_0 (\rho e^{2\pi i t})=p^{\rho,t}_{\theta}$, $0\leq\rho\leq\varepsilon$, $\theta\in\NS^1$. Finally, consider a smooth increasing function $\tilde{\rho}:[0,1]\rightarrow[\varepsilon,1]$ such that $\tilde{\rho}(0)=\varepsilon$ and $\tilde{\rho}(1)=1$. Then, the disk $\{\tilde{\gamma}^{r,\theta}\}$ is defined as $\tilde\gamma^{r,\theta}(t)=p^{\tilde{\rho}(r),t}_{\theta}$.
		
		Observe that the tangent space at each point of $\D^\theta$ is given by \begin{equation}\label{eq:TangentDiskArea}
		\langle \partial_x+\epsilon, \delta\partial_z+\beta\partial_y+\alpha\partial_w \rangle,
		\end{equation} where $\epsilon\in\Span\{ \partial_y,\partial_z,\partial_w\}$ and $\delta\partial_z+\beta\partial_y+\alpha\partial_w=B_{\lambda}^{\theta}-A_{\lambda}^{\theta}$ as vectors in $\R^3$. Obviously the embedding $\tilde{\gamma}^{r,\theta}$ is tangent to $\D^\theta$.
		\begin{figure}[h]
			\includegraphics[scale=2.5]{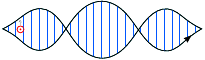}
			\centering
			\caption{Disk $\D^\theta$ with $\tilde{\gamma}=\tilde{\gamma}^{0,\theta}$ in its first lobe in the front Geiges projection.\label{reglada}}
		\end{figure}
	
	Consider the map $\tilde{F}:[0,2]\times\NS^1\times\NS^1(s,\theta,t)\rightarrow\NS^3$ defined by
	\[\tilde{F}(s,\theta,t)= \left\{ \begin{array}{lcc}
	(\tilde{\gamma}^{s,\theta})'(t) &   \text{if}  & s\in[0,1]\\
	
	F_{s-1}^{\theta,1}(t)& \text{if} &s\in[1,2].\\

	\end{array}
	\right.
	\]
	Since $\tilde{F}(0,\theta,-)=(\tilde{\gamma})'(-)$ and $\tilde{F}(1,\theta,-)=\partial_w$, for $\theta\in\NS^1$, the map $\tilde{F}$ can be regarded as a map from $\NS^2\times\NS^1\rightarrow\NS^3$. The $\mu$--invariant of the area twist is given by the degree mod $2$ of this map. In order to compute the degree of this map we compute the preimages of $\partial_z$. Since $\partial_z\notin\SD$ and $F(s,\theta,t)\in\SD$, for $s\geq1$, the equality $F(s,\theta,t)=\partial_z$ is only possible if $s\in[0,1)$. Ie we need to study the equation $(\tilde{\gamma}^{s,\theta})'(t)=\partial_z$, $s\in[0,1)$, this equality implies that 
		\item [(a)] $(\tilde{\gamma}^{s,\theta})'(t)\in\ker(dx)$, ie the front Geiges projection of $(\tilde{\gamma}^{r,\theta})$ has either a cusp or a vertical tangency at time $t$, in other words, $\tilde\gamma^{r,\theta}$ becomes tangent to the ``ruling'' (see Equation (\ref{eq:TangentDiskArea}));
		\item[(b)] $(\tilde{\gamma}^{s,\theta})'(t)\in\ker(dy)$, ie $y(B_{\lambda}^{\theta})-y(A_{\lambda}^{\theta})=0$, in other words, the area enclosed by the front Geiges projection of $\gamma^\theta$ between $A_{\lambda}^{\theta}$ and $B_{\lambda}^{\theta}$ is zero; and
		\item[(c)] $(\tilde{\gamma}^{s,\theta})'(t)\in\ker(dw)$, ie $w(B_{\lambda}^{\theta})-w(A_{\lambda}^{\theta})=0$, in other words, the slopes of the tangent lines of the front Geiges projection of $\gamma^\theta$ at $A_{\lambda}^{\theta}$ and $B_{\lambda}^{\theta}$ are the same.
		
		\begin{figure}[h!]
			\includegraphics[scale=2]{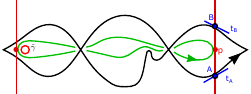}
			\centering
			\caption{Example where the first condition is fulfilled for a point $p$ but not the other two. \label{kerdx}}
		\end{figure}	
		
		After these observations, it is clear that the $\mu$--invariant of the area twist loop is $1$, see Figure \ref{partialz}.	
		\begin{figure}[h]
			\centering
			\subfloat[The tangent is $-\partial_z$.]{
				\label{fig:-PartialZ}
				\includegraphics[width=0.45\textwidth]{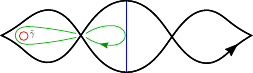}\label{partialz1}}
			\subfloat[The tangent is $\partial_z$.]{
				\label{fig:PartialZ}
				\includegraphics[width=0.45\textwidth]{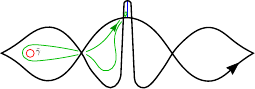}\label{partialz2}}
			\caption{Cases where the three conditions are fulfilled. \label{partialz}}
		\end{figure}
	\end{proof}

\section{The Area Invariant.}\label{AreaInvariantSection}
	In this Section we define a $\Z_2$-invariant, called \em Area Invariant\em, which coincides with the $\mu$--invariant defined in the previous Section in the horizontal case. Ie for any loop of horizontal embeddings $\gamma^\theta$ the \em Area Invariant \em satisfies $\Ar(\gamma^\theta)=\mu(\gamma^\theta,(\gamma^\theta)')\in\Z_2$. This invariant is defined studying the space $\Hor(\R^4)$ as a subspace of the space of horizontal immersions $\HorImm(\R^4)$, which is a flexible space \cite{PinoPresas}. The key point is that the strict horizontal immersions are encoded by means of the zeroes of a ``suitable" function, called the \em area function. \em We keep the same notation as in Section \ref{FormalHorizontal}.
	
	\subsection{Local models of Legendrians.}\label{localmodelsoflegendrians}
	Our goal is to understand generic local models of families of Legendrian immersions. Before studying the different models of Legendrian immersions depending on the codimension, we first define a proper subset of Legendrian immersions which has infinite codimension in the space $\LegImm(\R^3)$ and, thus, does not intersect a generic finite dimensional family. A point in a Legendrian is \em not injective \em if the preimage of the image has cardinality greater than one. The set of \em non injective images \em is the image of the set of non injective points. A Legendrian is \em special \em if the cardinality of the non injective images is non--finite. We denote the space of special Legendrian immersions by $\SpLegImm(\R^3)$.
	
If a Legendrian is special, there exists a Weinstein chart $(x,y,z)$ such that the Legendrian has at least two branches on the chart, one of them being the curve $(t,0,0)$ and the other one is $(t,f(t),f'(t))$, where $f(t)$ is a function with infinite zeroes such that $t=0$ is a zero and is limit of a sequence of zeroes. This implies that $f$ has an infinite order zero at $t=0$. This is an infinite codimension condition in the space of $\LegImm(\R^3)$. Therefore, the set $\SpLegImm(\R^3)$ is a subset of an infinite codimension stratified submanifold of the space of Legendrian immersions. We can dismiss this set when considering finite dimensional families, since any finite dimensional family can be generically perturbed to be disjoint from an infinite codimension subset. 
	
	First, consider the following space
	\begin{align*}
		\SSLegImm(\R^3)=&\lbrace\left(\gamma,v,t_0,t_1\right)\in\left(\SLegImm(\R^3)\setminus\SpLegImm(\R^3)\right)\times\R^3\times\NS^1\times\NS^1:\\
	&\gamma(t_0)=\gamma(t_1)=v\rbrace,
	\end{align*}
	equipped with the natural inclusion $\SSLegImm(\R^3)\subseteq\LegImm(\R^3)\times\R^3\times\NS^1\times\NS^1$. Denote by $\pi:\LegImm(\R^3)\times\R^3\times\NS^1\times\NS^1\to\LegImm(\R^3)$ the projection onto the first factor. Observe that $\pi(\SSLegImm(\R^3))=\SLegImm(\R^3)\setminus\SpLegImm(\R^3)$.
	
	The projection $\pi_{|\SSLegImm(\R^3)}:\SSLegImm(\R^3)\to\SLegImm(\R^3)$ is a local diffeomorphism. We just need to check that the fiber of $\pi$ is transverse to $\SSLegImm(\R^3)$; ie given $\alpha=(\gamma,v,t_0,t_1)\in\SSLegImm(\R^3)$, then $T_{\alpha}\mathcal{F}_{\gamma}\cap T_{\alpha}\SSLegImm(\R^3)=\lbrace 0\rbrace$, where $\mathcal{F}_\gamma=\pi^{-1}(\gamma)$. This is obvious since we are considering Legendrians with isolated non injective points. Thus, we have the following diagram:
	\[
	\xymatrix{
		\SSLegImm(\R^3)\ar@{->}[d]^{\text{loc. diff.}} \ar@{^{(}->}[r]
		  &*+[r]{\LegImm(\R^3)\times\R^3\times\NS^1\times\NS^1}\ar[d]\\
		\SLegImm(\R^3) \ar@{^{(}->}[r] &\LegImm(\R^3)}
	\]
	The space $\SLegImm(\R^3)\setminus\SpLegImm(\R^3)$ is not a stratified submanifold while 
	 $\SSLegImm(\R^3)$ it is. Thus, the first set can be regarded as the image by a local diffeomorphism of a stratified submanifold.  We can now compute local models of $R$--versal deformations (see Arnold, Gusein-Zade and Varchenko \cite[page  147]{Arnold}) upstairs and translate them downstairs just by using the local diffeomorphism.	
	
\begin{definition}\cite[page 49]{Arnold}
Consider a germ of a smooth real function $f:\R\rightarrow \R$ satisfying $f(0)=0$. We say that $f\in\Sigma^{\stackrel{k}{1,\ldots,1}}$ if $f^{(j)}(0)=0$, for $j\leq k$. If $f^{(k+1)}(0)\neq 0$, then we say that $f\in\Sigma^{\stackrel{k}{1,\ldots,1},0}\subseteq \Sigma^{\stackrel{k}{1,\ldots,1}}$.
\end{definition}

Fix $\gamma(t_0)=\gamma(t_1)$ a non injective image for a curve $\gamma\in\LegImm(\R^3)$. By \em Weinstein's tubular neighborhood theorem, \em we can assume that there is a chart  in which the branch of $\gamma$ through $t_0$ becomes the zero section in the Weinstein model. So, we can locally regard the other branch as a smooth real map with an order 2 zero at $t=0$. 

Let us now introduce some concepts that will be useful along the upcoming discussion.
\begin{remark} \cite[pages 147--148]{Arnold}
Given $f:\left(\R^n,0\right)\to\R$ a smooth map germ, we say that a deformation with parameter space $K=\R^k$ is a germ at the origin of a map-germ $F:\left(\R^n\times\R^k,0\right)\to\R$ such that $F(x,0)\equiv f(x)$.
We say that a deformation $G$ is $R-$equivalent to $F$ if
\[G(x,k)\equiv F\left(h(x,k),k\right) \text{ with }h:\left(\R^n\times\R^k,0\right)\to(\R^m,0)\text{ a smooth germ with }h(x,0)\equiv x.\]
We say that a deformation $F$ of $f$ is $R-$versal if any other deformation $G$ can be expressed as
\[G(x,\kappa)\equiv F(h(x,\kappa),\varphi(\kappa)) \text{ with }h(x,0)\equiv x\text{ and }\varphi(0)=0.\]

In order to understand the local models that we will produce, we consider the space of Legendrians up to reparametrizations. Note that this is equivalent to say that the local model representing the non injective point of the Legendrian that is given by a smooth function has to be considered up to $R$--versal deformations.
\end{remark}

We can then say that a point $(\gamma,v,t_0,t_1)\in\SSLegImm(\R^3)$ belongs to $\Sigma^I$ if and only if when we express $\gamma$ as the germ of a smooth real function around the origin near the tangency $(t_0,t_1)$, this germ of function belongs to $\Sigma^I$. Note that this is well defined.

Let us introduce some notation. Define the configuration space of two points in $\NS^1$ as $\Conf^2(\NS^1)=\{(\theta_0,\theta_1)\in\NS^1\times\NS^1:\theta_0\neq \theta_1\}$. If $e:\D^k\rightarrow\LegImm(\R^3)$ is a $k$--disk of Legendrian immersions  we want to check that \em Thom's Transversality Theorem \em can be applied to $e$ in order to obtain a disk which is transverse to $\Sigma^I$. By a map transversal to a stratified set we mean a map that is transversal to each stratum. We say that a property P for maps is generic if the set of maps that satisfy the property is a Baire set in the total space of maps. We are going to check that a property of the map $e:\D^k\rightarrow\LegImm(\R^3)$ is generic. Thus, by a standard compactness argument, we can further assume that the disk is very small. Ie take any finite covering for which we check that the property is satisfied in a Baire subset, then we are able to check the property in the initial disk since the finite intersection of Baire subsets is Baire. 

By the previous discussion, there is an associated map $\bar{e}:\D^k\times \NS^1\rightarrow\R$ that is constructed by identifying a tubular neighborhood of $e(0)$ with $J^1(\NS^1,\R)$, where the image of $e(0)$ corresponds to the zero section. Thus, $\bar{e}(z,-)$ is a function whose associated $1$--jet represents the Legendrian immersion $e(z)$ in a Weinstein chart. Define  $\hat{e}:\D^k\times\Conf^2(\NS^2)\rightarrow\R,(z,t_0,t_1)\mapsto \bar{e}(z,t_0)-\bar{e}(z,t_1)$, that is considered as a $k$--parametric family of real functions. Now, fix $r\in\Z^+$ the amount of regularity of the Banach manifold of maps in which we are working, ie we are working in the space of $C^r$ Legendrian immersions, and consider the holonomic lift of $\hat{e}$ into $J^r(\D^k\times\Conf^2(\NS^1),\R)$ and denote it by $j^r\hat{e}:\D^k\times \Conf^2(\NS^1)\rightarrow J^r(\D^k\times\Conf^2(\NS^1),\R)$. We have a stratified submanifold inside $J^r(\D^k\times\Conf^2(\NS^1),\R)$ defined by $F\in\Sigma^{\stackrel{k}{1,\ldots,1}}$, $k<r$, defined by the intersection of hyperplanes $F^{(j)}_{\theta_0}=F^{(j)}_{\theta_1}$ for $j\in\{0,1,\ldots,k\}$. In particular, the holonomic lift of $\hat{e}$ intersects $\Sigma^{\stackrel{k}{1,\ldots,1}}$ if there exist $z\in\D^k$ and $(t_0,t_1)\in\Conf^2(\NS^1)$ such that $e(z)^{(j)}(t_0)=e(z)^{(j)}(t_1)$ for $j\in\{0,1\ldots,k\}$.

Now we further reduce our domain. For this denote by $\mathcal{U}_\Delta$ an open neighborhood of the diagonal of $\NS^1\times\NS^1$. Thus, $j^r \hat{e}_{|\D^k\times(\mathcal{U}_\Delta \backslash{\Delta})}$ does not intersect $\Sigma^I$. So, by compactness, we may assume that there is a finite covering by squares of $\Conf^2(\NS^1)\backslash(\mathcal{U}_\Delta \backslash{\Delta})$ \[ \{ S_j=[t^{j}_{0}-\varepsilon,t^{j}_{0}+\varepsilon]\times[t^{j}_{1}-\varepsilon,t^{j}_{1}+\varepsilon]:[t^{j}_{0}-\varepsilon,t^{j}_{0}+\varepsilon]\cap[t^{j}_{0}-\varepsilon,t^{j}_{1}+\varepsilon]=\emptyset, j\in\{1,\ldots,N\} \}.\] Now, we try to perturb $\hat{e}$ along the domains $\D^k\times S_j$ in such a way that $\hat{e}_{|\D^k\times S_j}$ becomes transverse to $\Sigma^I$. But in that domain perturbations of $\hat{e}$ are in one to one correspondence with perturbations of $\bar{e}$. Therefore, we can apply the standard \em Thom's Transversality Theorem \em (see \cite[Theorem 2.3.2]{EliashMisch}) to conclude that the space of perturbations transverse to the stratified submanifold is a Baire set.

Thus, we have proven the following statement 

\begin{lemma}\label{lem:ThomTransversalityDiskOfImmersions}
Assume $0\leq k \leq r-2$. Let $e:\D^k\rightarrow\LegImm(\R^3)\backslash\SpLegImm(\R^3)$ be a disk. Then, there is a $C^r$--perturbation $\tilde{e}:\D^k\rightarrow\LegImm(\R^3)\backslash\SpLegImm(\R^3)$ of $e$ such that
\begin{itemize}
	\item [(i)] $\tilde{e}$ is transverse to $\SLegImm(\R^3)$ and
	\item [(ii)] the lift of $\tilde{e}$ to $\SSLegImm(\R^3)$ is transverse to each $\Sigma^I$, $|I|\leq r$.
\end{itemize}
\end{lemma}

From the previous discussion we conclude that the classification of generic local models of self--intersections of Legendrians tantamounts to classifying local germs of smooth real maps. In fact, we obtain a systematic method for classifying local models of generic $k$--dimensional disks of Legendrians by taking into account Arnold's classification of local germs of smooth functions (see \cite{Arnold}). In order to be precise we need to work on $\SSLegImm(\R^3)$.

$\bullet$\textit{\textbf{Codimension 0}}. A generic germ of a smooth function vanishing at the origin has generically non zero first derivative at the origin. This just reflects the fact that Legendrian immersions are generically embeddings.

\begin{figure}[h!]
	\includegraphics[scale=2.5]{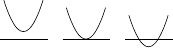}
	\centering
	\caption{Local model in codimension 1 representing a simple crossing. \label{Cod1}}
	\end{figure}

$\bullet$\textit{\textbf{Codimension 1}}. A generic $1$--parametric disk of smooth maps at the origin intersects $\Sigma^{1,0}$ at isolated points and $\Sigma^{1,1}$ in the empty set. By an application of \em Hadamard's lemma \em (see Milnor \cite[Lemma 2.1]{Milnor}) we can assume that the intersection with this stratum takes place at the germ at the origin of the map $f:\R\to\R, x\mapsto x^2$. We have that an $R$--versal deformation of $x^2$ at $0$ is given by $x^2+a$, where $a\in\Op(\{0\})$ (see \cite[page 151]{Arnold} for further details). This corresponds to the local model represented in Figure \ref{Cod1}. 
	
\begin{figure}[h]
		\centering
		\hspace*{\fill}
		\subfloat[Two simultaneous simple tangencies.]{
			\includegraphics[width=0.35\textwidth]{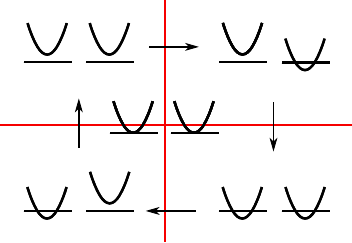}\label{Cod2Tang}}
			\hspace*{\fill}
		\subfloat[Cuspidal curve.]{
			\includegraphics[width=0.35\textwidth]{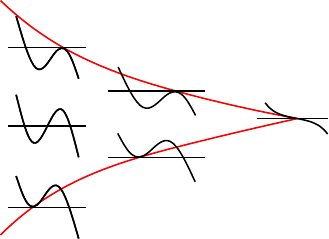}\label{Cod2Cusp}}\hfill
		\caption{Local models in codimension 2.}
		\hspace*{\fill}
	\end{figure}

$\bullet$\textit{\textbf{Codimension 2}}. A generic $2$--parametric family $\D$ of smooth maps at the origin intersects $\Sigma^1$ at curves and $\Sigma^{1,1,0}$ at isolated points. Fix a point in $\Sigma^{1,1,0}\cap\D$, we have that the crossing can be identified with the germ at the origin of the function $f:\R\to\R, t\mapsto t^3$. We know that an $R$--versal deformation of $t^3$ at $0$ is given by $\beta(t,a,b)=t^3+at+b$, where $(a,b)\in\D^2(\varepsilon)$ for  $\varepsilon>0$ small enough (see \cite[page 151]{Arnold} for details). Realize that the intersection of the disk with the stratum $\Sigma^{1}$ is given by the curve (seen as a curve of germs of functions at the origin) determined by the equations
	\[ \left. \begin{array}{lcc}
	         \beta(t,a,b)=0 \Longleftrightarrow t^3+at+b=0\\
	          \\ \frac{d}{dt}\beta(t,a,b)=0 \Longleftrightarrow 3t^2+a=0\\
	           \end{array}
	  \right\} \]
	with solutions $t=\pm\sqrt{\frac{-a}{3}}$ and the parameters satisfying the relationship
	\[ a^3 (3^{-1/2}-3^{-3/2})^2+b^2=0 \]
	which describes a cuspidal curve.  This local model is represented in Figure \ref{Cod2Cusp}.

\begin{figure}[h!]
	\includegraphics[scale=0.4]{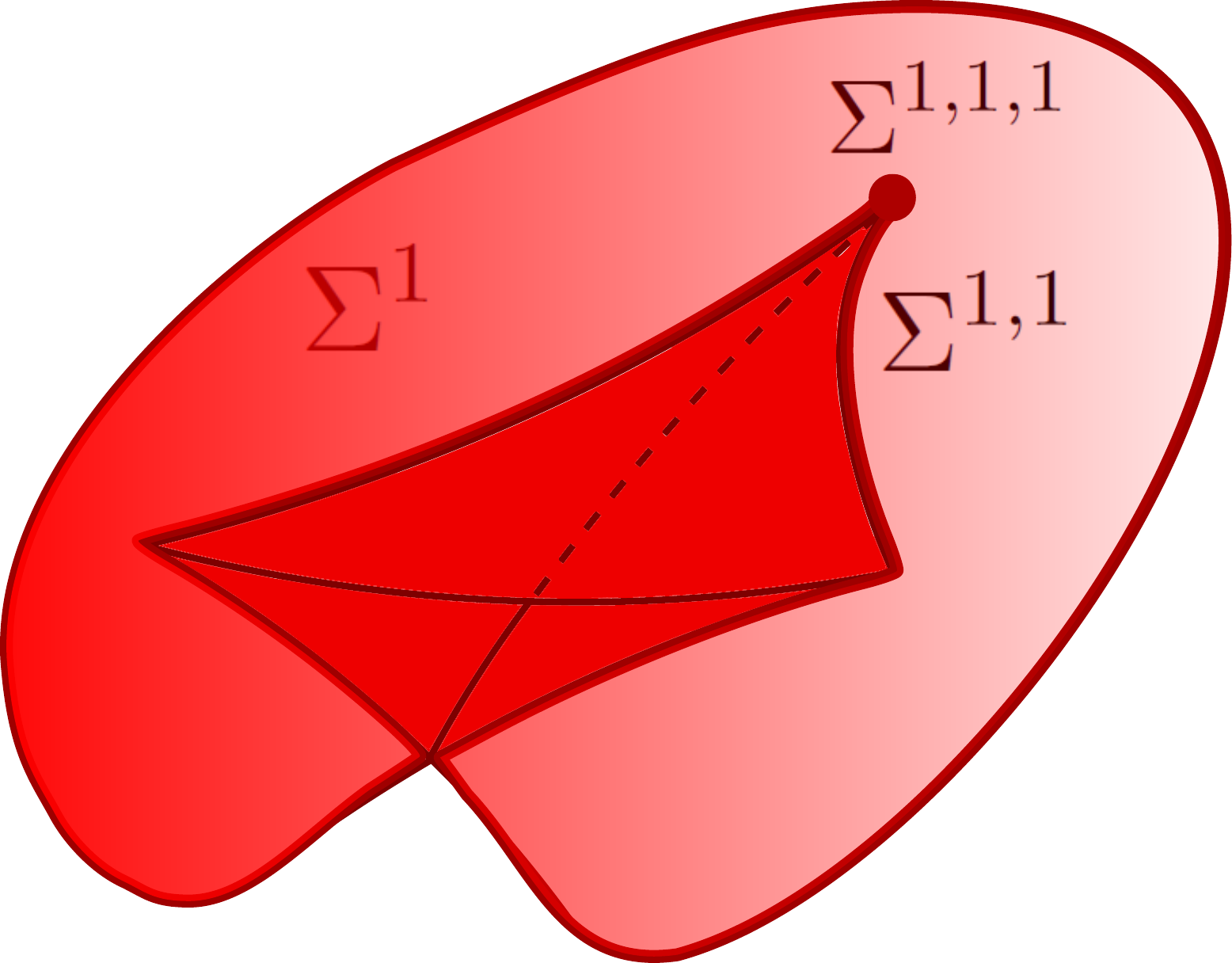}
	\centering
	\caption{Local model in codimension $3$. \label{Swallow}}
\end{figure}

$\bullet$ \textit{\textbf{Codimension 3}}. A generic $3$--parametric family $\D^3$ intersects $\Sigma ^1$ at surfaces, $\Sigma^{1,1}$ at curves and $\Sigma^{1,1,1,0}$ at isolated points. A point in $\Sigma^{1,1,1,0}\cap\D^3$ can be identified with the germ at the origin of the map $f:\R\to\R, x\mapsto x^4$. An $R$--versal deformations of $x^4$ at $0$ is given by $x^4+ax^2+bx+c$, where $(a,b,c)\in\D^3(\varepsilon)$ (see \cite[page 151]{Arnold} for further details). This gives raise to the local model represented in Figure \ref{Swallow} (see \cite[pages 47--48]{Arnold}). 

In order to understand deformations around points of $\SLegImm(\R^3)\setminus\SpLegImm(\R^3)$ we project the models from $\LegImm(\R^3)\times\R^3\times\NS^1\times\NS^1$ to $\LegImm(\R^3)$. Since the projection is a local diffeomorphism, we get the same local models, but different branches upstairs may come together downstairs. So, we have the same models but with the addition of cartesian product models that show up at points $\gamma\in\SLegImm(\R^3)$ that lie at the intersection of $2$ or more branches. For instance, in codimension $2$ we get a new model represented in Figure \ref{Cod2Tang}.

We have shown that any continuous disk in the space of Legendrian immersions can be perturbed to be in generic position. However, in this article we usually have a disk of horizontal immersions and we want to check that a perturbation of this disk renders a generic perturbation of its Geiges projection inside the space of Legendrian immersions. 

\begin{lemma}\label{lem:ThomTransversalityDiskOfHorizontalImmersions}
	Assume that $0\leq k\leq r-2$. Let $e:\D^k\rightarrow\HorImm(\R^4)$ be a disk such that $e_G=\pi_G \circ e:\D^k\rightarrow\LegImm(\R^3)\backslash\SpLegImm(\R^4)$. Then, there is a $C^r$--perturbation $\tilde{e}:\D^k\rightarrow\HorImm(\R^4)$ of $e$ such that
	\begin{itemize}
		\item [(i)] $\tilde{e}_G:\D^k\rightarrow\LegImm(\R^3)\backslash\SpLegImm(\R^3)$,
		\item [(ii)] $\tilde{e}_G$ is transverse to $\SLegImm(\R^3)$ and
		\item [(iii)] the lift of $\tilde{e}_G$ to $\SSLegImm(\R^3)$ is transverse to each $\Sigma^I$, $|I|\leq r$.
	\end{itemize}
\end{lemma}
\begin{proof}
It follows directly from Lemma \ref{lem:ThomTransversalityDiskOfImmersions} and Lemma \ref{DeltaBarquillo}. This last lemma will be proven later on and it allows to lift multiparametric families of Legendrian immersions.
\end{proof}

For later use we want to check that 

\begin{lemma}\label{lem:SigmaBanach}
$\Sigma^{1,0}\subseteq\LegImm(\R^3)$ is a locally closed codimension $1$ immersed Banach submanifold.
\end{lemma}
Note that the lift of $\Sigma^{1,0}$ to $\SSLegImm(\R^3)$ is, in fact, an embedded submanifold. However, it does not help our purposes.
\begin{proof}
	We already noted that $\LegImm(\R^3)$ is an immersed Banach submanifold of $\Imm(\NS^1,\R^3)$ (see Remark \ref{rem:FrechetImm}). In order to check that $\Sigma^{1,0}$ is a codimension $1$ locally closed immersed submanifold we just need to define an atlas of $\LegImm(\R^3)$ with charts satisfying that the pullback by the chart map of the intersection of $\Sigma^{1,0}$ with the chart is a closed codimension $1$ Banach subspace. 
	
	For a point $\gamma\in\Sigma^{1,0}$, the exponential map is a local diffeomorphism because its differential at $\gamma$ is the identity (see \cite[Theorem 5.2, page 15]{Lang}). By using a Weinstein chart around $\gamma$ we just define the exponential by the formula \ref{eq:ExponentialLegImm}. We can simply check that given a vector $(f,g) \in \Maps^r(\NS^1, \R) \times \Maps^r(\NS^1, \R)$, the condition for it  to be tangent to $\Sigma^{1,0}$ reads as $f(\theta_0)=0$, where $\theta_0$ is the self--intersection time. This is obviously a closed codimension $1$ subspace, this concludes the proof. 
\end{proof}

\subsection{The Area Invariant.} 
	Fix a loop $\gamma^\theta$ in $\Hor(\R^4)$ with $\Rot_L (\gamma^\theta)=0$. This means that $\gamma^\theta$ is trivial \em as a loop of horizontal immersions\em. Thus, there exists a disk $\D=\{\tilde{\gamma}^{z}\}$ in $\HorImm(\R^4)$ such that $\tilde{\gamma}^{e^{2\pi i \theta}}=\gamma^{\theta}$. We want to study the obstruction for that disk to belong to $\Hor(\R^4)$.

\begin{figure}[h]
		\includegraphics[scale=0.25]{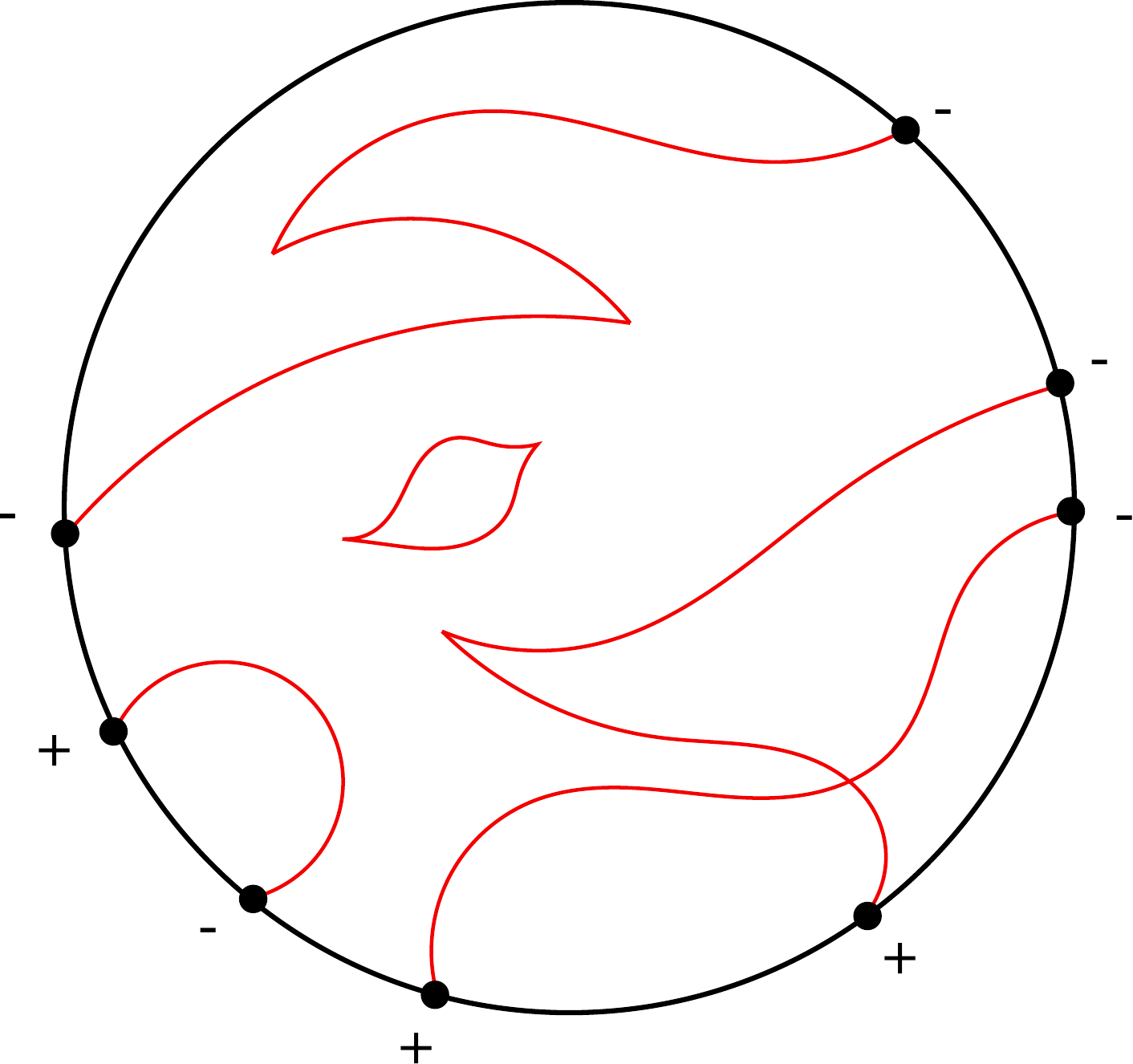}
		\centering
		\caption{Generic disk $\D_{G}$. The red curves represent strict Legendrian immersions. \label{ModelDisk}}
	\end{figure}

	By Lemma \ref{lem:ThomTransversalityDiskOfHorizontalImmersions} we can assume that $\D_G$ is in generic position, ie $\D_G$ contains curves ($1$--dimensional connected manifolds), smooth away from the cusps, of strict Legendrian immersions which may intersect in $\Int(\D_G)$. Moreover, we may assume that in $\partial\D_G\cap\SLegImm(\R^3)$ we have only elements in $\Sigma^{1,0}$, hence the number of strict Legendrian immersions in $\partial\D_G$  is even. See Figure \ref{ModelDisk}.
	
	Let $\tilde{\gamma}^{z}_G\in\D_G$ be a strict Legendrian immersion. Write $\tilde{\gamma}^{z}_G(t)=(x(t),z(t),w(t))^\intercal$ and assume that $p=\tilde{\gamma}^{z}_G (t_0)=\tilde{\gamma}^{z}_G (t_1)$ is a transverse self--intersection point, ie of type $\Sigma^{1,0}$; $t_0,t_1\in\NS^1$. We define 
	\begin{equation*}
		\varepsilon_A (\tilde{\gamma}^{z}_G,p)=\delta\int_{t_0}^{t_1} z(t)x'(t)dt,
	\end{equation*}
	where the factor $\delta\in\{\pm1\}$ corresponds to the following sign convention: \em we always integrate along a segment that starts at the lower branch of the curve in the front projection and follow the orientation of the curve.\em \footnote{The tangent space of the  two branches in $\R^3$ at the intersection point $p\in\R^3$ defines a framing of $\xi_p$. The choice just ensures that the orientation of the framing is always positive.} Notice that this convention is well--defined since we are working with self--intersections of type $\Sigma^{1,0}$. If $\tilde{\gamma}^{z}_G$ has only one tangency point, then we write $\varepsilon_A (\tilde{\gamma}^{z}_G)=\varepsilon_A (\tilde{\gamma}^{z}_G, p)$. We call $\varepsilon_A$ the \em area function\em. The sign of $\varepsilon_A$ over a curve of strict Legendrian immersions is reversed when it passes through a cusp (a point of type $\Sigma^{1,1,0}$). Indeed, note that when crossing these points the lower branch becomes the upper branch and vice versa, thus producing a change of sign in the area function. In particular, the absolute value of the area function is well--defined over a cusp. Note that, by genericity, we may assume that $\varepsilon_A$ is non zero over the cusp points. Hence, the area function is continuous over the curves of strict Legendrian immersions except at the cusp points. Moreover, if $\tilde{\gamma}^{z}_G$ is a strict Legendrian immersion, then $\tilde{\gamma}^{z}$ is a horizontal embedding if and only if $\varepsilon_A (\tilde{\gamma}^{z}_G,\cdot)$ is non zero over all tangency points of $\tilde{\gamma}^{z}_G$. Thus, to study the obstruction for $\D$ to belong to $\Hor(\R^4)$ we need to study the zero set of the area function. Observe that, since $\partial\D_G$ is the Geiges projection of $\gamma^{\theta}\in\Hor(\R^4)$, we obtain that $\varepsilon_A (\tilde{\gamma}^{e^{2\pi i \theta}}_G)$ is non zero. We define the sign of $\tilde{\gamma}^{z}_G$ as the sign of $\varepsilon_A (\tilde{\gamma}^{z}_G)$. Denote by $\PMas(\D_G)$ the set of positive points in $\partial\D_G$. The following lemma is just a simple exercise: 
    \begin{figure}[h]
		\centering
		\subfloat[Capping disk in $\LegImm(\R^3)$.]{
			\label{fig:CappingAT}
			\includegraphics[width=0.35\textwidth]{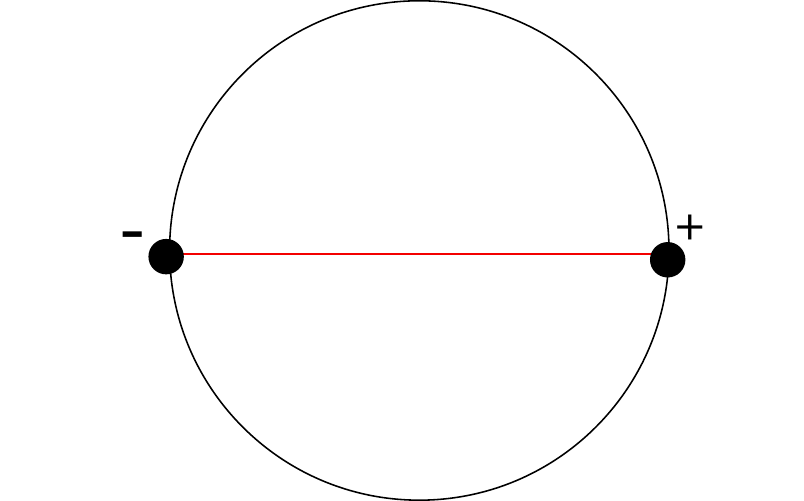}}
		\subfloat[Curve of strict Legendrian immersions in the front Geiges projection.]{
			\label{fig:CappingAT2}
			\includegraphics[width=0.35\textwidth]{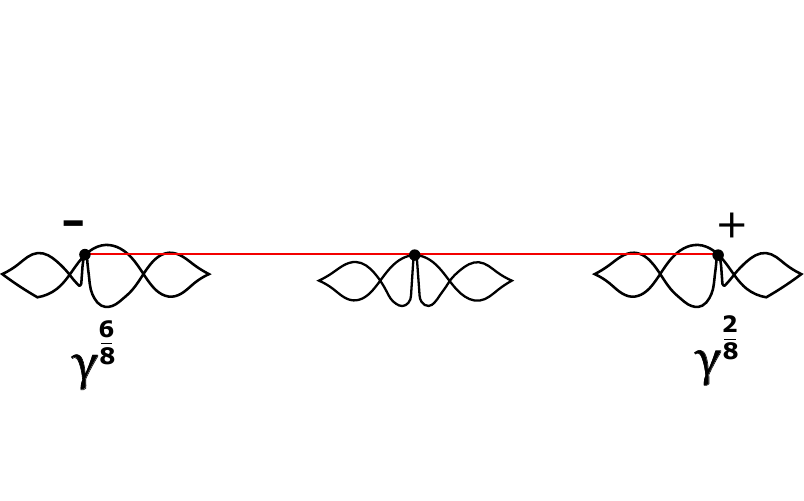}}
		
		\caption{Geiges projection of the standard capping disk through horizontal immersions of the area twist loop.}\label{DiagramAT}
	\end{figure}
	
	\begin{lemma}\label{ObstructionCurves}
		Let $C\subseteq\D_G$ be a curve of strict Legendrian immersions. Assume that one of the three following conditions holds:
		\begin{itemize}
			\item [(i)] The boundary points of $C$ have the same sign and the number of cusps in $C$ is odd.
			\item[(ii)] The boundary points of $C$ have different signs and the number of cusps in $C$ is even.
			\item[(iii)] The boundary of $C$ is empty and the number of cusps in $C$ is odd.
		\end{itemize} 
		Then, $\D$ contains at least one strict horizontal immersion.
	\end{lemma}
	
	\begin{figure}[h]
		\centering
		\subfloat[Curves with the same sign in the boundary and an odd number of cusps.]{
			\label{fig:SameSingBoundary}
			\includegraphics[width=0.47\textwidth]{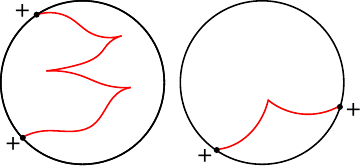}}
		\subfloat[Curves with different signs in the boundary and an even number of cusps.]{
			\label{fig:DifferentSignBoundary}
			\includegraphics[width=0.45\textwidth]{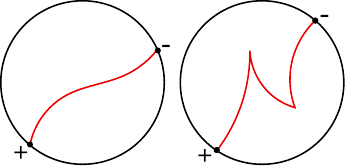}}
		
		\caption{Some examples of curves of strict Legendrian immersions that contains the Geiges projection of a strict horizontal immersion.}
	\end{figure}
	
	\begin{definition}
		Let $\Cusp(\D_G)$ be the set of cusps of $\D_G$. We define the \em area \em of $\D$ as
		\begin{equation*}
			\Ar(\D)=(\#\PMas(\D_G)+\#\Cusp(\D_G))\mod2.
		\end{equation*}
	\end{definition}
	
	\begin{lemma}
		Assume that $\Ar(\D)=1$, then $\D$ contains at least one strict horizontal immersion.
	\end{lemma}
	\begin{proof}
		We proceed by induction over the number of strict Legendrian immersions in the boundary of $\D_G$. Observe that, since $\Ar(\D)=1$ there has to be a curve of strict Legendrian immersions in $\D_G$.
		
		\textit{$0$ points in the boundary.} 
		Since $\Ar(\D)=1$ then $\#\Cusp(\D_G)$ must be odd. Hence, there exists a closed curve with an odd number of cusps. By Lemma \ref{ObstructionCurves}(iii), we conclude this case.
		
		\textit{$2$ points in the boundary.}
		We have two cases:
		\begin{itemize}
			\item [(a)] Assume that the two points in the boundary have the same sign. Further assume there exists a curve of strict Legendrian immersions with the same sign in the boundary and an odd number of cusps. Use Lemma \ref{ObstructionCurves}(i) to conclude. Otherwise, there is a closed curve of strict Legendrian immersions with an odd number of cusps, by Lemma \ref{ObstructionCurves}(iii) we conclude this case.
			\item[(b)] Assume that there exist two points in the boundary which have different signs. If the curve of strict Legendrian immersions that connects them has an even number of cusps, we conclude by Lemma \ref{ObstructionCurves}(ii). In other case, there must be a closed curve with an odd number of cusps. By Lemma \ref{ObstructionCurves}(iii) we conclude this case.   
		\end{itemize}
		\textit{$2n+2$ points in the boundary.}
		We have two cases:
		\begin{itemize}
			\item [(a)] Assume that the points in the boundary connected by a curve of strict Legendrian immersions have the same sign. Further assume there exists a curve with the same sign in the boundary and an odd number of cusps. We conclude by Lemma \ref{ObstructionCurves}(i). Otherwise, there exists a closed curve with an odd number of cusps by Lemma \ref{ObstructionCurves}(iii) and we conclude this case.
			\item[(b)] Assume that we have a curve $C$ with opposite signs in the boundary. Then, if the number of cusps is even, we apply Lemma \ref{ObstructionCurves}(ii) to conclude. On the other hand, if the number of cusps is odd, consider as a new finite set of curves the old one minus $C$ and apply the induction hypothesis.
		\end{itemize}
	\end{proof} 
	
	We want to show that the assignment $\gamma^\theta\mapsto\Ar(\D)$ defines a homotopy invariant. 
	Thus, let $\gamma_{0}^{\theta}$ and $\gamma_{1}^{\theta}$ two homotopic loops of horizontal embeddings with rotation number zero. For that, choose $\D_0=\{\tilde{\gamma}_{0}^{z}\}$ and $\D_1=\{\tilde{\gamma}_{1}^{z}\}$ two disks in $\HorImm(\R^4)$ such that $\tilde{\gamma}_{0}^{e^{2\pi i\theta}} =\gamma_{0}^{\theta}$ and $\tilde{\gamma}_{1}^{e^{2\pi i\theta}}=\gamma_{1}^{\theta}$. We need to study the relationship between the configuration of curves of strict Legendrian immersions in $(\D_0)_G$ and $(\D_1)_G$. Denote by $\gamma_{t}^{\theta}$ the homotopy between $\gamma_{0}^{\theta}$ and $\gamma_{1}^{\theta}$. Since $\pi_2 (\HorImm(\R^4))=0$ we can assume that we have a solid cylinder $\mathcal{C}yl=[0,1]\times\D$ in $\HorImm(\R^4)$ such that $\{t\}\times\partial\D$ coincides with $\gamma_{t}^{\theta}$, $\D_0 =\{0\}\times\D$ and $\D_1 =\{1\}\times\D$.
		\begin{figure}[h]
			\includegraphics[scale=0.04]{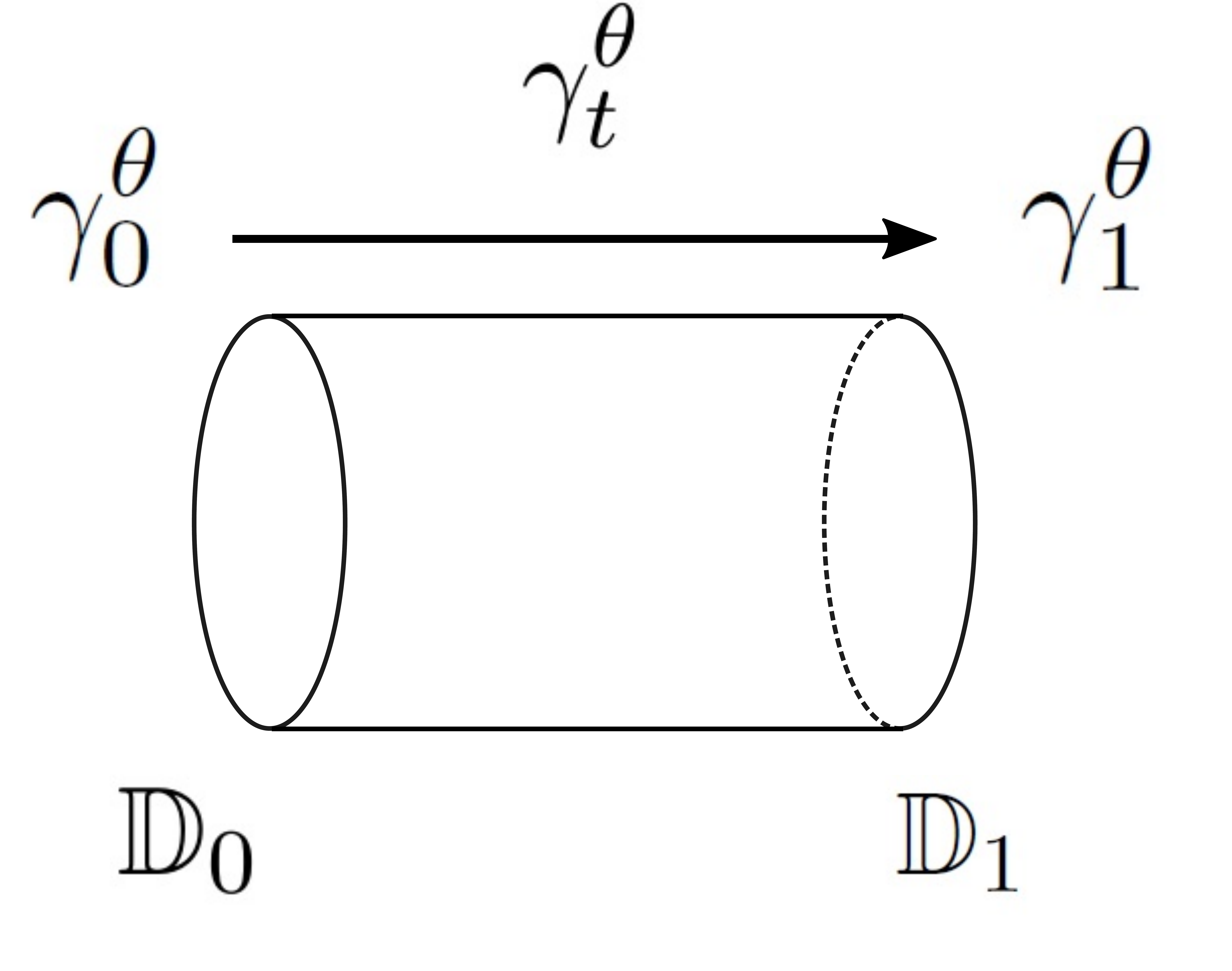}
			\centering
			\caption{Solid cylinder in $\HorImm(\R^4)$ connecting $\D_0$ and $\D_1$}
		\end{figure}	
 Let $\mathcal{C}yl_{G}=[0,1]\times\D_G$ be the Geiges projection of the cylinder. Denote by $\tilde{L}=[0,1]\times\partial\D$ the lateral boundary of the cylinder and $L=\pi_{G}(\tilde{L})$ its Geiges projection. In particular, every point in $L=[0,1]\times\partial\D_G$ is a Legendrian embedding. By genericity, we can assume that every Legendrian in $\mathcal{C}yl_G\subseteq\LegImm(\R^3)$ is the projection of an element in $\SSLegImm(\R^3)$ via the canonical projection \[\pi:\LegImm(\R^3)\times\R^3\times\NS^1\times\NS^1\to\LegImm(\R^3).\]
In addition, we can assume that $S_0:=\pi(\Sigma^{1,1,1,0})\cap\mathcal{C}yl_G$ is a set of isolated points, $S_1:=\pi(\Sigma^{1,1,0})\cap\mathcal{C}yl_G$ a set of immersed curves and $S_2:=\pi(\Sigma^{1,0})\cap\mathcal{C}yl_G$ a set of immersed surfaces. In addition, by genericity, we have that different branches of $S_1$ do not intersect, different branches of $S_1$ and $S_2$ intersect in isolated points $S_{12}=S_1\cap S_2$, and different branches of $S_2$ intersect in a set of embedded curves $S_{22}$. We can also assume that 
\begin{itemize}
\item[(i)] $S_0\cap L=\emptyset$,
\item[(ii)] $S_{1L}:=S_1\cap L$ is a finite set of points,
\item[(iii)] $S_{2L}:=S_2\cap L$ is a set of embedded curves.
\end{itemize}

\begin{figure}[h]
			\includegraphics[scale=0.15]{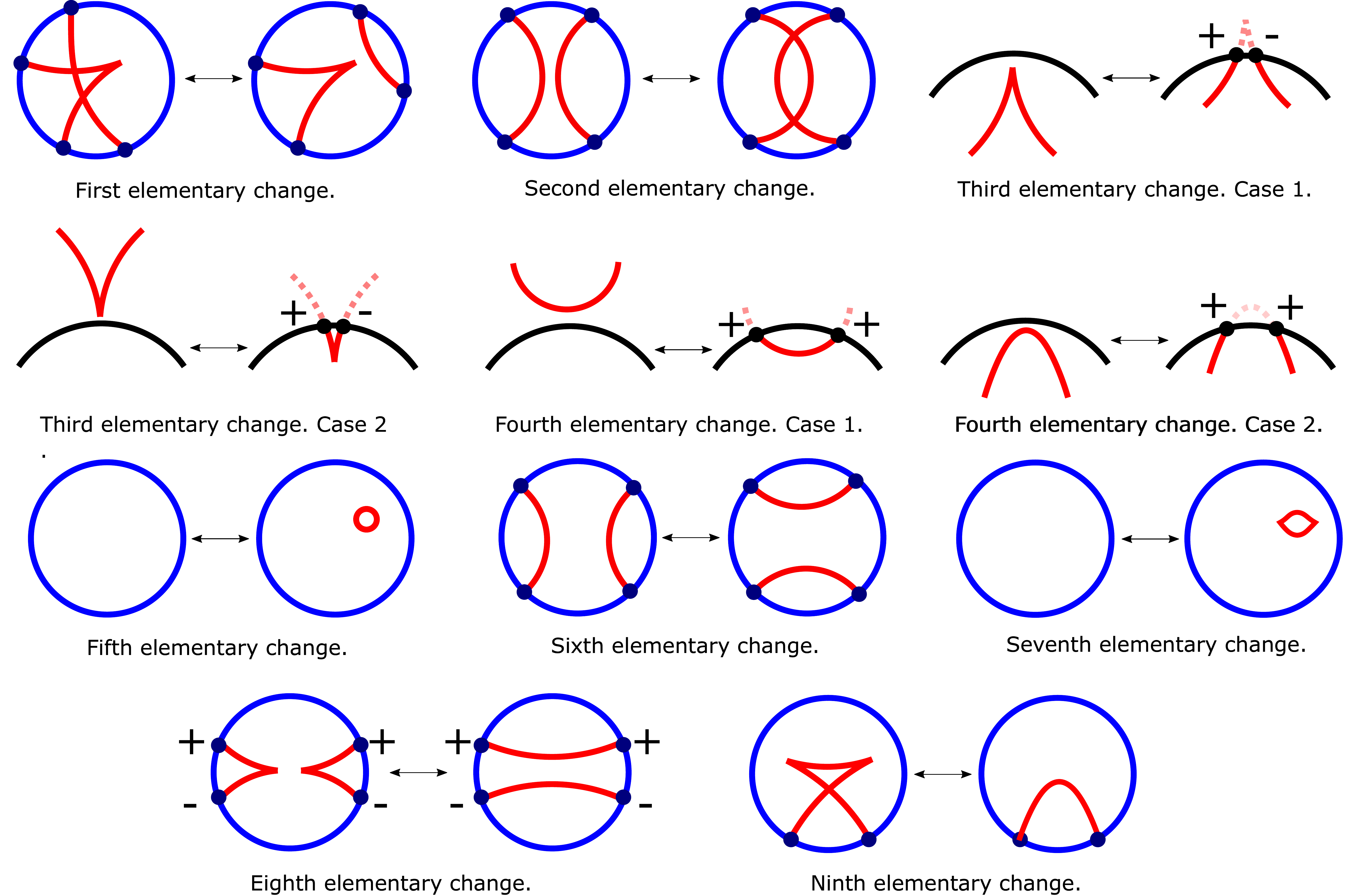}
			\centering
			\caption{List of elementary changes. The boundary of the cylinder is represented in black. Blue circles represent local neighborhoods of the slices. The signs in the boundary points are defined by the area function. Of course the same pictures with signs reversed could happen. \label{ElementaryChanges}}
		\end{figure}

Assume that the ``height function'' for the cylinder $h:\mathcal{C}yl_G=[0,1]\times\D_G\to [0,1]$ is a Morse function for $S_{1}, S_{2L}$ and $S_{22}$. Thus, with these assumptions, the $1$--parametric family of slices $\lbrace t\rbrace\times\D_G, t\in[0,1]$, induces a movie of curves defined by the intersections of $\pi(\Sigma ^{1,0})$ with each slice. Let us discuss the topology of the configuration of $\SLegImm(\R^3)$ at each slice (see Figure \ref{ElementaryChanges}). The changes of the topology of the curves happen at specific times due to the following phenomena:
\begin{itemize}
\item Points in $S_{12}$. If such a point lies on a slice $\{t\}\times\D_G$ the movie around that time corresponds to the crossing of a curve with a cusp with another immersed curve (see Figure \ref{ElementaryChanges}, First elementary change).
\item Points in $S_{22}$ This corresponds to the crossing of two immersed curves (see Figure \ref{ElementaryChanges}, Second elementary change).
\item Points in $S_{1L}$. This corresponds to the crossing of a cusp with the boundary (see Figure \ref{ElementaryChanges}, Third elementary change, cases 1 and 2).
\item Appearance (disappearance) of points in $S_{2L}\cap\left(\lbrace t\rbrace\times\D_G\right)$ for some $t\in[0,1]$. This generically corresponds to the appearance (or disappearance) of an immersed curve during the $1$--parametric family. The time $t\in[0,1]$ represents the time where this curve appears in the movie of slices (Fourth elementary change, cases 1 and 2).
\item Critical points of $h$ for $S_2$. A critical point in $S_2$ for $h$ represents either the appearance/disappearance of a closed embedded curve: the case of a minimum/maximum;  or an index $1$ surgery in the set of curves: the case of an index $1$ critical point (see Figure \ref{ElementaryChanges},  Fifth and Sixth elementary changes).
\item Critical points of $h$ for $S_1$. This corresponds to the appearance/disapearance of an embedded curve (respectively minimun/maximum). What is special is that the curve has a double cusp (See Figure \ref{ElementaryChanges}, Seventh elementary change and Eighth elementary change).
\item Points in $S_0$. This corresponds to the codimension $3$ swallow--tail singularity (see Figure \ref{ElementaryChanges}, Ninth elementary change). 
\end{itemize}		
		
	\begin{corollary}\label{AreaEsInvariante}
		Let $\gamma_{0}^{\theta}$ and $\gamma_{1}^{\theta}$ be two homotopic loops of horizontal embeddings with rotation number zero. Let $\D_0=\{\tilde{\gamma}_{0}^{z}\}$ and $\D_1=\{\tilde{\gamma}_{1}^{z}\}$ be two disks in $\HorImm(\R^4)$ such that $\tilde{\gamma}_{0}^{e^{2\pi i\theta}} =\gamma_{0}^{\theta}$ and $\tilde{\gamma}_{1}^{e^{2\pi i\theta}}=\gamma_{1}^{\theta}$. Then,
		\begin{equation*}
			\Ar(\D_0)=\Ar(\D_1).
		\end{equation*}
	\end{corollary}
	\begin{proof}
		By the previous discussion, any generic homotopy between two disks is a composition of the nine elementary changes: they come with suitable signs in order to be geometrically realizable, see Figure \ref{ElementaryChanges}. Since these nine elementary changes preserve the area the result follows.
	\end{proof}
	We claim:
	
	\begin{theorem}\label{AreaInvariant}
		Let $\gamma^\theta$ be a loop in $\Hor(\R^4)$ such that $\Rot_L (\gamma^\theta)=0$. Let $\D=\{\tilde{\gamma}^{z}\}$ be some disk in $\HorImm(\R^4)$ such that $\tilde{\gamma}^{e^{2\pi i\theta}}=\gamma^{\theta}$. The assignment
		\begin{equation*}
			\gamma^\theta\rightarrow\Ar(\gamma^\theta):=\Ar(\D),
		\end{equation*}
		defines a homotopy $\Z_2$--invariant, called the \em Area Invariant. \em Alternatively, we define the \em Area Invariant \em as the number of strict horizontal immersions in $\D\mod2$. Both definitions coincide.
	\end{theorem}
\begin{proof}
We are left with checking that the Area invariant is computed by the number of strict horizontal immersions mod $2$. We proceed by induction in the number of zeroes of the Area function over the strict immersed curves. If there are no zeroes, then Lemma \ref{ObstructionCurves} implies that the Area invariant is zero. Assume that it is true for $k=0, \ldots, n$. Then, we want to study a configuration with $n+1$ zeroes. If all of them lie in the same curve, we have that Lemma \ref{ObstructionCurves} implies the result. If not, the curves can be divided into two subsets such that they have $n \geq k_1>0$ and $n \geq k_2>0$ zeroes each with $n+1=k_1+k_2$. By induction, we conclude the proof.
\end{proof}	
	
\begin{remark}
Let us check that the Area Invariant coincides with the $\mu$--invariant defined in the previous Section in the horizontal case. Ie for any loop of horizontal embeddings $\gamma^\theta$ the \em Area Invariant \em satisfies \[\Ar(\gamma^\theta)=\mu(\gamma^\theta,(\gamma^\theta)')\in\Z_2.\] The $\mu$--invariant describes the $\Z_2$ factor given by the \em kernel \em of the homomorphism
		\begin{equation*}
			f:\pi_1 (\FHor(\R^4))\rightarrow\pi_1 (\FHorO(\R^4)).
		\end{equation*}
Consider the natural inclusion $i:\Hor(\R^4)\hookrightarrow\FHor(\R^4)$ and let $\pi_1(\Hor(\R^4))_f:=\pi_1(i)^{-1}(\ker(f))$. In order to prove that the Area Invariant computes the $\mu$--invariant in the formal case, we must check that the next diagram commutes:
		\begin{displaymath}
		\xymatrix@M=10pt{
			\pi_1(\Hor(\R^4))_f\ar[r]^{\pi_1(i)}\ar[dr]_{\Ar} & \ker(f)\ar[d]^{\mu}\\
			& \Z_2}
		\end{displaymath}
		
		Let $\gamma^\theta$ be a loop of horizontal embeddings which lies in $\pi_1(\Hor(\R^4))_f$. Let $\D=\{\tilde{\gamma}^z\}$ be a disk in $\HorImm(\R^4)$ such that  $\tilde{\gamma}^{e^{2\pi i\theta}}=\gamma^{\theta}$. Assume that $\D$ has exactly $n$ strict horizontal immersions $\beta^1,\ldots,\beta^n$. Let $\D_G$ be the Geiges projection of $\D$, note that by genericity, we can assume that $\beta^{1}_{G},\ldots,\beta^{n}_{G}$ belong to $\Sigma^{1,0}$. Because of Definition \ref{def:AreaTwist} any small disk whose diagrams of strict Legendrians is as in Figure \ref{fig:CappingAT}  has as boundary circle an area twist. Then collapsing the disk into small neighborhoods of the curves $\beta^{j}_{G}$, as we show in Figure \ref{HomotopyAT}, $\gamma^\theta$ is homotopic to a concatenation of $n$ area twist loops. Observe that the Area Invariant and the $\mu$--invariant of an area twist loop is one. Thus,  $\Ar(\gamma^\theta)=\mu(\gamma^\theta,(\gamma^\theta)')=n \mod 2$ (see Lemma \ref{AT.Z2.Invariant}). This concludes the argument.
		\begin{figure}[h]
			\includegraphics[scale=0.14]{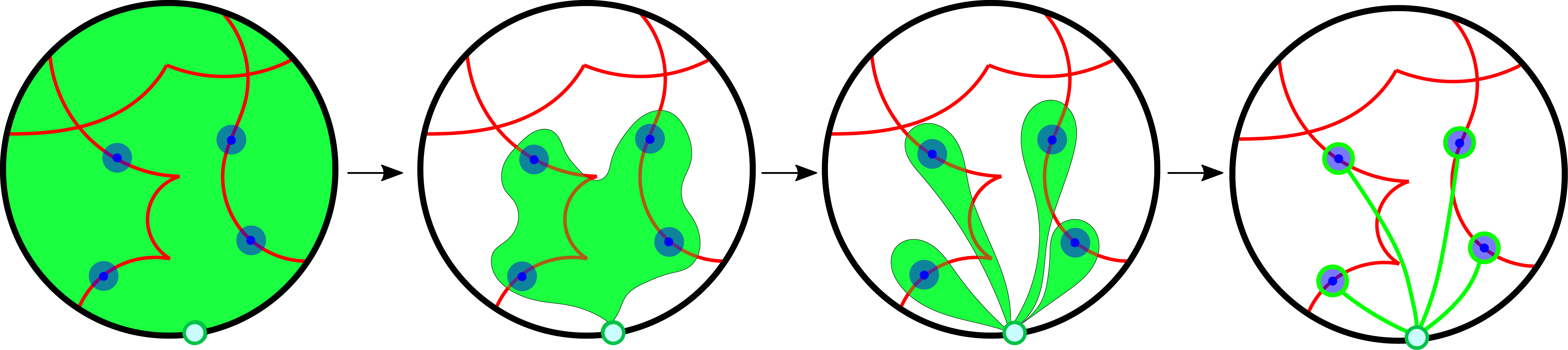}
			\centering
			\caption{Homotopy between $\gamma^\theta$ and the concatenation of $n$ area twist loops for the case $n=4$ represented in $\D_G$. We have drawn the zeroes of the area function in blue. The concatenation of $4$ area twist loops is represented by the green curve in the last step. \label{HomotopyAT}}
		\end{figure}
	\end{remark}

\section{$h$--principle at $\pi_1$--level for horizontal embeddings.}\label{h--Principle}
We keep the same notation as in Sections \ref{FormalHorizontal} and \ref{AreaInvariantSection}. 

\subsection{The main theorem.}
\begin{theorem}\label{hPrincipleHorizontalEmbeddings}
The inclusion $\Hor(\R^4)\hookrightarrow\FHor(\R^4)$ induces an isomorphism of fundamental groups.
\end{theorem}

We are assuming a base point in the previous inclusion. We claim that the isomorphism works for any choice of base point, ie for all the connected components: they are given by the rotation number of the horizontal embedding (see Remark \ref{rem:RotationsHorizontal}). To prove this result we proceed as follows. Consider the homomorphism 
\begin{center}
	$\begin{array}{rccl}
	\Rot_L\colon & \pi_1 (\Hor(\R^4)) & \longrightarrow &  \Z \\
	& [\gamma^\theta]& \longmapsto & \Rot_L(\gamma^\theta)
	\end{array}$
\end{center}
and define $\pi_1(\Hor(\R^4))_0$ as the \em kernel \em of $\Rot_L$. The homomorphism $\Rot_L$ is surjective. Indeed, it is enough to check that for any horizontal embedding $\beta\in\Hor(\R^4)$ and any integer $k\in\Z$ there exists a loop of horizontal embeddings $\beta^\theta$, based at $\beta$, with $\Rot_L(\beta^{\theta})=k$. Define $\beta^{\theta}_{G}$ as the Legendrian immersion obtained by a rotation of angle $2\pi k\theta$ of $\beta_G$ in the Lagrangian projection. The lift $\beta^{\theta}$ of $\beta^{\theta}_{G}$ to $\Hor(\R^4)$ satisfies the desired property.

It follows from Theorem \ref{AreaInvariant} that the Area Invariant defines an homomorphism $\Ar:\pi_1(\Hor(\R^4))_0\rightarrow \Z_2$. Moreover, we will prove the following

\begin{theorem}\label{AreaIsomorphism}
	The Area Invariant defines an isomorphism
	\begin{equation} \label{eq:super}
	\Ar:\pi_1(\Hor(\R^4))_0\rightarrow \Z_2.
	\end{equation}
 In particular, $\pi_1(\Hor(\R^4))\cong\Z\oplus\Z_2$.
\end{theorem}

Moreover, the techniques developed in this Section will allow us to prove (it will be proven at the end of Subsection \ref{subsec:flexible})
\begin{lemma}\label{HorizontalAproximationOfLoops}
	The homomorphism $\pi_1(\Hor(\R^4))_0\rightarrow\pi_1 (\Emb(\NS^1,\R^4))$, induced by the inclusion $\Hor(\R^4)\hookrightarrow\Emb(\NS^1,\R^4)$, is surjective.
\end{lemma}

Thus, since the two elements in $\pi_1 (\Hor(\R^4))_0$ (the trivial loop and the area twist loop) are trivial as loops of smooth embeddings and the space $\Emb(\NS^1,\R^4)$ is path--connected, we conclude the following
\begin{corollary}\label{EmbeddingsSimplyConnected}
	The space $\Emb(\NS^1,\R^4)$ is simply connected.
\end{corollary}

This is a ``Legendrian'' proof of a well known result (see \cite{Budney} Proposition $3.9 (4)$). We expect to be able to compute higher rank homotopy groups of this space generalizing the discussion below since we know the space of Legendrian immersions behaves nicely. Ie the crossings types correspond to tangencies of smooth functions with the horizontal line. It is known that these tangencies can be assumed to be either Morse or birth--death, see \cite{Igusa}. For higher dimensional families this is a rather non trivial result. In a forthcoming project, we will use this fact to collect non trivial information about the high rank homotopy groups of the space $\Emb(\NS^1,\R^4)$.

It is clear that Theorem \ref{AreaIsomorphism} and Corollary \ref{EmbeddingsSimplyConnected}, together with the exact sequence stated in Lemma \ref{Pi1FHor}, implies Theorem \ref{hPrincipleHorizontalEmbeddings}. 

The purpose of this Section is to reduce the proof of Theorem \ref{AreaIsomorphism} to a technical cancellation Theorem, see Proposition \ref{TwoComponents}, that is of independent interest and is proven in Section \ref{sec:Core}. This Section presents several results needed in orden to reduce and to develop the proof of the cancellation Theorem. 

An overview of what follows is
\begin{itemize}
\item In Subsection \ref{subsec:flexible}, we introduce several results: Lemma \ref{DeltaBarquillo} states that the area function over a connected component (see below for a precise definition of {\em connected component}) can be easily deformed without affecting the other connected components; then we focus on deforming the immersed Legendrians produced by Geiges projection. The main result states that the only obstruction to having an $h$--principle for embeddings of Legendrians lies on the fact that the double stabilization of a Legendrian (main tool to approximate a smooth curve) changes the embedding class of the Legendrian. The key remark (Lemma \ref{DeltaSeta}) shows that this is not the case in the horizontal case. We give a multiparametric proof of the fact that (double) stabilizations allow to approximate any smooth curve: Lemma \ref{LegendrianAproximation}. What is left to complete a full $h$--principle is to deal with the positions of the cusps points. This is done just in $1$--parametric families that is enough for the purposes of this article. The proof contained in \cite{CasalsdelPino} controls those cups in higher dimensional families.
\item In Subsection \ref{sub:same}, we start studying the injectivity of the morphism (\ref{eq:super}). This will be the task of the remaining part of the article. For that we need to show that the parity of the zeroes (Area invariant) of a capping disk on $\HorImm(\R^4)$ for a loop on $\Hor(\R^4)$ fully recovers the homotopy class of the loop. Thus, we need to cancel the zeroes by pairs. In this Subsection, we assume that the two zeroes to be cancelled do lie in the same connected component. In that case, a careful use of Lemma \ref{DeltaBarquillo} is all we need to conclude.
\item In Subsection \ref{sub:changedia}, we explain how to change the diagram of the capping disk associated to an element in $\pi_1(\Hor(\R^4))_0$ in order to place two zeroes in the same connected component. We prove in these pages a previous key result (Proposition \ref{TwoComponents}) explaining how to deform a path connecting the two zeroes through curves that, except for a finite number of points, are Legendrian embeddings in Geiges projection. The idea is to build a map of the square into $\Hor(\R^4)$. The bottom side of the square is the initial path. The other three edges correspond to a ``much better" path that connects the two zeroes  by a curve of strict immersions; ie we have found a curve making them live in the same connected component. Finally, we explain how to exploit the success in order to use the previous square to change the diagram of curves in a suitable way (creating an index one surgery in the diagram). 
\item In Subsection \ref{sub:global}, we take a suitable path joining the two zeroes, then the index $1$ surgery previously built allows to {\em connect} the connected components of the two different zeroes.
\end{itemize}

\subsection{Creating flexibility.}\label{subsec:flexible}

\begin{figure} 
	\includegraphics[scale=0.06]{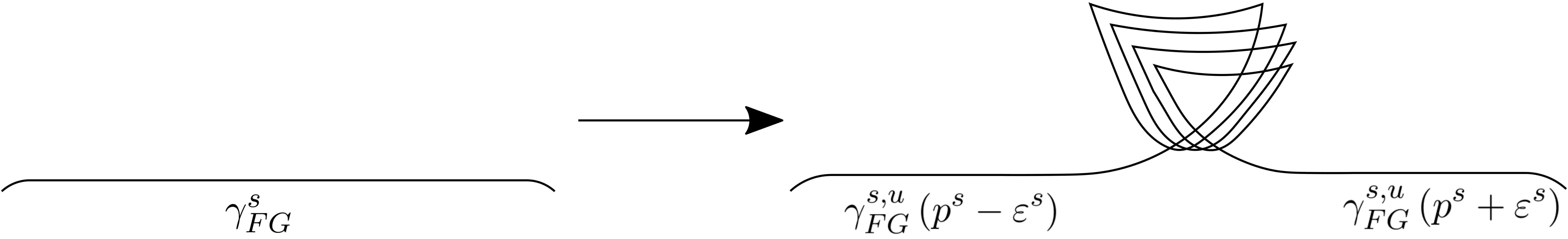}
	\centering
	\caption{Construction of $\gamma^{s,u}$ via its front Geiges projection.\label{DeltaBarcos}}
\end{figure}

 \begin{lemma}\label{DeltaBarquillo}
	Let $K$ be a compact parameter space and $\gamma^s\in\HorImm([0,1],\R^4)$, $s\in K$. Let $A:K\rightarrow\R$ be a continuous map. Let $p^s\in[0,1]$, $s\in K$, be any continuous family of regular points for the front Geiges projection of $\gamma^s$; and $\varepsilon>0$ sufficiently small. Then, there exists a $1$--parametric family $\gamma^{s,u}\in\HorImm([0,1],\R^4)$, $u\in[0,1]$, satisfying:
	\begin{itemize}
		\item [(i)] $\gamma^{s,0}=\gamma^{s}$.
		\item[(ii)] The Geiges projection of $\gamma^{s}$ and $\gamma^{s,u}$ are $C^0$--close and \footnote{We are working in the standard Engel structure in $\R^4(x,y,z,w)$ given by $\SD=\ker(dy-zdx)\cap\ker(dz-wdx)$.}
		\begin{equation*}
			\gamma^{s,u}-\gamma^s=\begin{cases}
				(0,0,0,0) & \text{in $[0,p^s-\varepsilon]$,} \\
				(0,uA(s),0,0) & \text{in $[p^s+\varepsilon,1]$.}
			\end{cases}
		\end{equation*} 
	\end{itemize} 
	Moreover, if $\gamma_{G}^{s}$ is an embedding then $\gamma_{G}^{s,u}$ is also an embedding.
\end{lemma}
\begin{proof}
	Define $\gamma^{s,u}=\gamma^s$ in $[0,p^s-\varepsilon]$ and $\gamma_{G}^{s,u}=\gamma_{G}^{s}$ in $[p^s+\varepsilon,1]$. It is enough to define the front Geiges projection of $\gamma^{s,u}$ in $[p^s-\varepsilon,p^s+\varepsilon]$. Let $N\in\Z^+$, consider the equispaced partition $\{a_{0}^{s} = p^s-\varepsilon<a_{1}^{s}<\ldots<a_{N+1}^{s} =p^s+\varepsilon\}$ of $[p^s-\varepsilon,p^s+\varepsilon]$. Define $\gamma_{FG}^{s,u}$ in $[p^s-\varepsilon,p^s+\varepsilon]$ just adding a Reidemeister I move of area $\frac{uA(s)}{N}$ to $\gamma_{FG}^{s}$ at each point $\gamma_{FG}^{s}(a_{i}^{s})$, $i=1,\ldots,N$; see  Figure \ref{DeltaBarcos}. The $C^0$--closeness follows from choosing $N$ large enough.	
\end{proof}

\begin{remark}
	Lemma \ref{DeltaBarquillo} also applies to horizontal loops, provided that the zero total area condition is preserved. We do it as follows: let $\gamma\in\HorImm(\R^4)$ be a horizontal immersion, $p,n\in\NS^1$ different regular points of the front Geiges projection of $\gamma$ and $A\in\R$. We write 
	\[ \gamma\#\delta_p (A)\#\delta_n (-A) \in\HorImm(\R^4) \]
	to denote the horizontal immersions obtained from $\gamma$ when, using Lemma \ref{DeltaBarquillo}, we \em add area \em $A$ near the point $\gamma(p)$ and $-A$ near the point $\gamma(n)$. Moreover, if $\beta\in\LegImm(\R^3)$  satisfies the zero total area condition, ie $\beta=\gamma_G$ for $\gamma\in\HorImm(\R^4)$, and $p,n\in\NS^1$ are regular points of the front projection of $\beta$ we write 
	\[ \beta\#\delta_p (A)\#\delta_n (-A)\in\LegImm(\R^3) \]
	to denote the Geiges projection of $\gamma\#\delta_p (A)\#\delta_n (-A)\in\HorImm(\R^4)$.
\end{remark}

\begin{figure}[h]
	\includegraphics[scale=0.85]{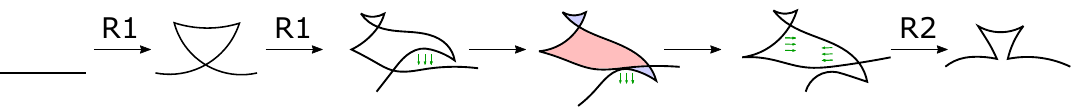}
	\centering
	\caption{Sequence of moves in the front projection of $\gamma$ creating a double stabilization. }\label{CreacionSeta}
\end{figure}

\begin{definition}\label{dobleEstab}
	Let $\gamma\in\LegImm(\R^3)$. A double stabilization (DS) of $\gamma$ is a $1$--parametric family $\gamma^u\in\LegImm(\R^3)$, $u\in[0,1]$, such that 

\begin{itemize}
\item[(i)] $\gamma^0=\gamma$,
\item[(ii)] $\gamma^1$ is $\gamma$ stabilized once negatively and once positively.
\end{itemize}	
The homotopy $\gamma^u, u\in[0,1]$, is explicitly described in Figure \ref{CreacionSeta}.
\end{definition}
We can define the DS of a horizontal immersion $\gamma$ via the Geiges projection, just being careful enough not to change the total area condition. To do this pick $t_0,t_1\in\NS^1$ two regular points of the front Geiges projection of $\gamma$. Use $t_0$ to add a DS and use Lemma \ref{DeltaBarquillo} to add the necessary area over $\gamma(t_1)$ to ensure that the area condition is fulfilled. We can do it even better as the next key result shows.

\begin{lemma}\label{DeltaSeta}
	The double stabilization is a well--defined operation for horizontal embeddings.	
\end{lemma}
\begin{proof}
The total area enclosed by the only tangency point created in a DS is nonzero.
\end{proof}

We do not claim that this construction is unique. In fact, it is possible to produce two paths not formally homotopic, connecting the end of the DS with the initial loop: just draw a different picture movie with a different area balance and then it is clear that the concatenation of the first path with the inverse of the second one is a loop with non-trivial area invariant.
However, this is enough for our purposes.

The next result shows that DSs can be used in order to approximate families of smooth curves by Legendrian curves.

\begin{theorem}\label{LegendrianAproximation}
	Given $K=[0,1]\times[0,1]$, $C=[0,1]\times\{0\}\subseteq K$ and a family $\gamma^k:[0,1]\rightarrow\R^3$, $k\in K$, such that 
	\begin{itemize}
		\item[(i)] $\gamma^k$ are embeddings,
		\item [(ii)] $\gamma^k(t)$ are Legendrian for $t\in\Op\{0,1\}$,
		\item[(iii)] $\gamma^k$ are Legendrian for every $k\in C$.
	\end{itemize}
Then, there exists a family of curves $\gamma^{k,u}:[0,1]\rightarrow\R^3$, $k\in K$, $u\in[0,1]$, such that
\begin{itemize}
	\item [(i)] $\gamma^{k,0}=\gamma^k$,
	\item[(ii)] $\gamma^{k,u}$ are immersions and $C^0$--close to $\gamma^k$,
	\item[(iii)] $\gamma^{k,1}$ are Legendrian embeddings,
	\item[(iv)] $\gamma^{k,u}(t)=\gamma^k (t)$ for $t\in\Op\{0,1\}$,
	\item[(v)] there exist a finite sequence $0<u_1<u_2<\ldots<u_N<1$ and $\delta>0$ small enough such that 
	\begin{itemize}
		\item [(a)] $\gamma^{k,u}$ are Legendrian embeddings for every $k\in C$, $u\neq u_1,\ldots,u_N$,
		\item[(b)] $\gamma^{k,u}$, $u\in[u_i-\delta,u_i+\delta]$, is a DS of $\gamma^{k,u_i -\delta}$ for every $k\in C$ and $i=1,\ldots,N$.
	\end{itemize}
\end{itemize}
\end{theorem}

This Theorem together with Lemma \ref{DeltaBarquillo} easily implies the following

\begin{corollary}\label{cor:RelojHorizontal}
	Given $K=[0,1]\times[0,1]$, $C=[0,1]\times\{0\}\subseteq K$ and a family $\gamma^k:[0,1]\rightarrow\R^4$, $k\in K$, such that 
	\begin{itemize}
		\item[(i)] $\gamma^{k}_{G}$ are embeddings,
		\item [(ii)] $\gamma^k(t)$ are horizontal for $t\in\Op\{0,1\}$,
		\item[(iii)] $\gamma^k$ are horizontal for every $k\in C$.
	\end{itemize}
	Then, there exists a family of curves $\gamma^{k,u}:[0,1]\rightarrow\R^4$, $k\in K$, $u\in[0,1]$, such that
	\begin{itemize}
		\item [(i)] $\gamma^{k,0}=\gamma^k$,
		\item[(ii)] $\gamma^{k,u}$ are embeddings and $C^0$--close to $\gamma^k$,
		\item[(iii)] $\gamma^{k,1}$ are horizontal embeddings,
		\item[(iv)] $\gamma^{k,u}(t)=\gamma^k (t)$ for $t\in\Op\{0,1\}$.
	\end{itemize}
\end{corollary}
\begin{proof}
	Apply the previous Theorem to the family $\gamma^{k}_{G}$ to obtain a family $\gamma^{k,u}_{G}$, $u\in[0,1],$ of immersed curves in $\R^3$. To conclude the proof we just need to approximate the area coordinate (the $w$--coordinate). This can be easily done by using Lemma \ref{DeltaBarquillo}.
\end{proof}

We can relax the hypothesis (i) for later use. We will prove a particular result just with the interval as the parameter space.  

\begin{lemma}\label{lem:TechnicalAprox}
	Let $\gamma^{s}\in\Imm(\NS^1,\R^3)$, $s\in[-1,1]$, be a $1$--parametric family of immersions such that 
	\begin{itemize}
		\item[(i)] $\gamma^{s}\in\Emb(\NS^1,\R^3)$ for $s\neq0$,
		\item[(ii)] $\gamma^{0}$ has a generic self--intersection.
	\end{itemize}
	Then, there exists a deformation $\gamma^{s,u}\in\Imm(\NS^1,\R^3)$, $u\in[0,1]$, satisfying
	\begin{itemize}
		\item [(a)] $\gamma^{s,0}=\gamma^{s}$,
		\item[(b)] $\gamma^{s,u}\in\Emb(\NS^1,\R^3)$ for $s\neq0$, 
		\item[(c)] $\gamma^{0,u}$, $u\in[0,1]$, has a generic self--intersection,
		\item[(d)] $\gamma^{s,1}\in\LegImm(\R^3)$ for $s\in\Op(\{0\})$,
		\item[(e)] $\gamma^{s,u}=\gamma^{s}$ for $s\notin\Op(\{0\})$.
	\end{itemize}
\end{lemma}
\begin{proof}
Let $t_0,t_1\in\NS^1$ be the self--intersection times of $\gamma^0$. After applying an isotopy we may assume that there exists $0<\varepsilon$ such that 
\begin{itemize}
	\item There exists a time--dependent vector field $v_{0}$, supported near $\gamma^{0}(t_0)$, such that $\gamma^{0,\pm u}$, $u\in(0,\varepsilon)$, is obtained from $\gamma^{0}$ by pushing the branch through $\gamma(t_0)$ in the direction of $\pm v^{0}$.
	\item $\gamma^{0}_{|\Op(\{t_0,t_1\})}$ is Legendrian. 
\end{itemize}
Apply Theorem \ref{LegendrianAproximation} to the two smooth embedded arcs $\gamma^{0}_{\NS^1\backslash\Op(\{t_0,t_1\})}$ to construct a path $\tilde\gamma^{0,u}\in\Imm(\NS^1,\R^3)$, $u\in[0,1]$, where $\tilde\gamma^{0,1}\in\SLegImm(\R^3)$ has a generic self--intersection. Define a family of vector fields $\tilde{v}^{u}$, $u\in[0,1]$, supported near $\gamma^{0,u}(t_0)$ such that $\tilde{v}^0=v^0$ and $\tilde{v}^{1}$ is a contact vector field. Define $\tilde{\gamma}^{\pm s,u}\in\Emb(\NS^1,\R^3)$, $(s,u)\in(0,\varepsilon)\times[0,1]$, by pushing the branch through $\tilde\gamma^{0,u}(t_0)$ in the direction of $\pm v^{u}$. It follows that $\tilde{\gamma}^{\pm s,1}\in\Leg(\R^3)$ for $s\neq0$. Finally, consider a smooth bump function $\chi:[-\varepsilon,\varepsilon]\rightarrow[0,1]$ such that $\chi_{|\Op(\{0\})}\equiv 1$ and $\chi_{|\Op(\{-\varepsilon,\varepsilon\})}\equiv 0$. The desired family is
\begin{equation*}
	\gamma^{s,u}=\begin{cases}
		\gamma^{s}& \text{ if $s\in[-1,1]\backslash[-\varepsilon,\varepsilon]$} \\
		\tilde{\gamma}^{s,\chi(s)} & \text{ if $s\in[-\varepsilon,\varepsilon]$.} 
	\end{cases}
\end{equation*} 
\end{proof}

\begin{corollary}\label{cor:AproxHorizontalRot0}
	Let $\gamma^s:[0,1]\rightarrow\R^4$ be a $1$--parametric family of embeddings such that the curves $\gamma^i$, $i\in\{0,1\}$; are horizontal and $\gamma^s(t)$, $s\in[0,1]$; are horizontal for $t\in\Op\{0,1\}$. Then, there exists a homotopy $\gamma^{s,u}:[0,1]\rightarrow\R^4$, $(s,u)\in[0,1]^2$, satisfying 
	\begin{itemize}
		\item [(i)] $\gamma^{s,0}=\gamma^s$,
		\item[(ii)] $\gamma^{s,u}$ are embeddings and $C^0$--close to $\gamma^s$,
		\item[(iii)] $\gamma^{s,1}$ are horizontal embeddings,
		\item[(iv)] $\gamma^{s,u}(t)=\gamma^s (t)$ for $t\in\Op\{0,1\}$,
		\item [(v)] $\gamma^{i,u}$, $i\in\{0,1\}$, are horizontal embeddings.
	\end{itemize}
\end{corollary}
\begin{proof}
	We may assume that there exist $0<s_1<\cdots<s_N<1$ such that
	\begin{itemize}
		\item $\gamma^{s}_{G}$ are embeddings for $s\in[0,1]\backslash\{s_1,\ldots,s_N\}$,
		\item $\gamma^{s_i}_{G}$, $i\in\{1,\ldots,N\}$, has a generic self--intersection.
	\end{itemize} Moreover, after applying  Lemma \ref{lem:TechnicalAprox} we may also assume that $\gamma^{s}_G$ is Legendrian for $s\in\Op(\{s_1,\ldots,s_N\})$. Now the result follows by an application of Corollary \ref{cor:RelojHorizontal} to the paths  $\gamma^s$, $s\in[0,1]\backslash\Op\{s_1,\ldots,s_N\}$.
\end{proof}

\begin{remark}
	The result can be extended to any compact parameter space $K$ and any closed subset $C\subseteq K$. A direct consequence is a parametric version of the Fuchs--Tabachnikov's result (see \cite{FuchsTabachnikov}) that asserts that any two Legendrian embeddings which are smoothly isotopic are, after a  finite number of stabilizations, Legendrian isotopic. If we assume that we have a $k$--disk $\D^k\subseteq\FHor(\R^4)$ such that $\partial\D^k\subseteq\Hor(\R^4)$, it follows that we can homotope this ball inside $\FHor(\R^4)$ into a ball in $\Hor(\R^4)$. For more details, see \cite{CasalsdelPino}. Nevertheless, our argument does not assume that we have a disk of formal horizontal embeddings but just a disk of horizontal immersions; ie $\D^k\subseteq\HorImm(\R^4)$. Thus, we need to homotope  this disk, through the space of horizontal immersions, into a disk of horizontal embeddings. In other words, we are not using as initial datum the homotopy type of the space of smooth embeddings of the circle into $\R^4$: part of the data provided by having a disk of formal horizontal embeddings. In a sense, we are studying the relative homotopy type of the space of horizontal embedding inside the space of horizontal immersions.
\end{remark}

The following result works only for $1$-parametric families over the segment. This is the main obstacle in our proof to generalize to a complete $h$-principle. In \cite{CasalsdelPino} a generalization  to deal with the cusp points in higher dimensional families is provided and this proves a complete $h$-principle. Our result just claims that for a $1$--parametric family of Legendrian immersions the cusps can be assumed to be constant (in time and position) for the whole family.

\begin{lemma}\label{EliminationRI}
	Let $\gamma^s$, $s\in[0,1]$, be a path in $\Leg([0,1],\R^3)$ such that $\Op\{0,1\}$ are regular points for the front projection. Then, there exists a $2$--parametric family $\gamma^{s,u}\in\Leg([0,1],\R^3)$, $s,u\in[0,1]$; satisfying
	\begin{itemize}
		\item [(i)] $\gamma^{s,0}=\gamma^s$ and $\gamma^{s,u}(t)=\gamma^s(t)$ for $t\in\Op\{0,1\}$.
		\item[(ii)] $\gamma^{s,u}$ and $\gamma^{s}$ are $C^0$--close.
		\item[(iii)] $\gamma^{s,1}$ is generated by a sequence of Legendrian Reidemeister moves of (only) type $II $and $III$.
	\end{itemize} 
\end{lemma}
\begin{proof}
	The family $\gamma^s$ is understood as $\gamma:[0,1]\times[0,1]\rightarrow\R^3$, $(t,s)\mapsto \gamma^s(t)$. By \em Thom's Transversality Theorem \em (see \cite[Theorem 2.3.2]{EliashMisch}) we assume that the sets of cusp points in $[0,1]\times[0,1]$ are embedded curves. Moreover, we assume that the height function $h(t,s)=s$, restricted to these curves, has a finite number of maximum and minimum points: all of them non degenerate. These points correspond to Reidemeister I moves. Take the point with the lowest height among all the minima. Since there is not any other minimum with a lower height, we can find a curve $C$ joining this point with $\{s=0\}$ that does not intersect any other curve of cusp points and which does not have any critical point. Now, we can remove this minimum adding a ($C^0$--small) family of Reidemeister I moves over $C$. Keep going to remove all the minima. To remove the maxima we do likewise.
\end{proof}

Thus we obtain
\begin{corollary}\label{CuspSameTime}
	Up to reparametrization, there exist $c_1,\ldots,c_m\in[0,1]$ such that the cusp points of $\gamma^{s,1}$ are at fixed times $c_1,\ldots,c_m$; for all $s\in[0,1]$.
\end{corollary}

\begin{proof}[Proof of Theorem \ref{LegendrianAproximation}.]
	Assume first that $\gamma^k$, $k\in C$, does not have any cusp point. Since the $h$--principle for smooth immersions works relative to the domain and the parameters and is $C^0$--dense (see \cite{EliashMisch}), we can assume that the front projection of each $\gamma^k$ is immersed. We declare the slope of the formal derivative of the front projection $\gamma_F$ of $\gamma$ to be the value of the coordinate $z$ ($\partial_z$ is the projection direction). If $\gamma_F$ is the projection of a Legendrian this definition of formal derivative coincides with the derivative of the front projection. Therefore, this equips the whole family $\gamma^k$ with a formal Legendrian structure.

Take an equidistant partition $u_0=0<u_1 =\frac{1}{N+1}<\ldots<u_N =\frac{N}{N+1}<u_{N+1}=1$ of $[0,1]$. Consider the collection of intervals $I_i =[u_i, u_{i+1}]$, $i=0,\ldots,N$ and let $K(N)\in\Z^+$. Define $\gamma^{k,u}$ inductively over $\{I_i\}$ as follows
\begin{itemize}
	\item [(i)] In $I_0$ and $I_N$ define  $\gamma^{k,u}=\gamma^k$.
	\item[(ii)] In $I_i$, $i=1,\ldots,N-1$, define $\gamma^{k,u}$ to be $K(N)$ DSs of $\gamma^{k,u_i}$ approximating the formal derivate of $\gamma^{k,u_i}$ at time $t=\frac{2i+1}{2(N+1)}$ as described in Figures \ref{reloj} and \ref{InterpolationFina} (for $K(N)=2$). Let us provide the details of the construction. The depicted blue segment represents $\gamma_{F}^{k}(I_i)$ that can be assumed to be ($C^{\infty}$--close to) a straight segment for $N$ large enough. The green segments represent the formal derivative, again for $N$ large enough is $C^0$--close to a constant. We further assume that the formal derivative has angle zero with respect to the horizontal axis, but we are reduced to that case composing with a linear transformation in $\R^2(x,y)$ to make the formal derivative zero, ie  the formal derivative is defined as $(1,z(t))$, we just change coordinates by using the unique linear transformation
	$$ A_{z(t)} : \R^2 \to \R^2 $$
	satisfying $A_{z(t)}\left( \begin{array}{c} 0 \\ 1 \end{array} \right)= \left( \begin{array}{c} 0 \\ 1 \end{array} \right)$ and  $A_{z(t)}\left( \begin{array}{c} 1 \\ z(t)  \end{array} \right)= \left( \begin{array}{c} 1 \\ 0 \end{array} \right)$.
This canonically lifts to a contact transformation. This allows us to assume that the formal derivative is zero. This is the reason why we are depicting only the case of formal derivative horizontal.
	
	Coming back to the figure, the number below each image is the angle between the curve (straight segment) and the formal derivative. The black curves are the $K(N)$ DSs that $C^0$--approximate $\gamma^k$ with prescribed derivative given by the green curve (in our pictures slope zero). The geometric idea is to perform first the negative stabilizations\footnote{We follow the standard convention that a negative stabilization is an upward zig--zag and viceversa for the positive ones.} and later on the positive ones. The rule is for the case in which the curve to be approximated points downwards we make the negative stabilizations small and the positive ones big. In particular, for $K(N)$ sufficiently large, we can make the derivative of the approximating curve arbitrarily close to the formal one (horizontal).
	\item[(iii)] The rule for the case in which the angle between the curve to be approximated and the formal derivative (horizontal in our pictures) tends to $\pi$ radians is the same. First negative stabilizations and afterwards positive ones. The key point is that the front projected curve cannot remain an embedding when reaching $\pi$: increasing the angle. The reason is that the $y$ coordinate of the ending point of the curve for angle slightly over $\pi$ is below the $y$ coordinate of the starting point. Since the end of the curve is decreasing (positive stabilizations), there is a moment in which the positive stabilization has to cross (through  $K(N)-1$ Reidemeister type II moves) over the curve that connects the negative with the positive stabilizations. Note that no front tangencies are created. This is done in order to be able to balance the DS and place them in the middle at angle $\pi$ and then, push them to the beginning when going over $\pi$ radians. This makes possible to proceed simetrically when coming from angles bigger than $\pi$. In other words, the approximation of a curve with angle $\alpha$ is the symmetry with respect to the $y$ axis of the one with angle $2\pi-\alpha$.
\end{itemize}

\begin{figure}[h!]
	\includegraphics[scale=0.17]{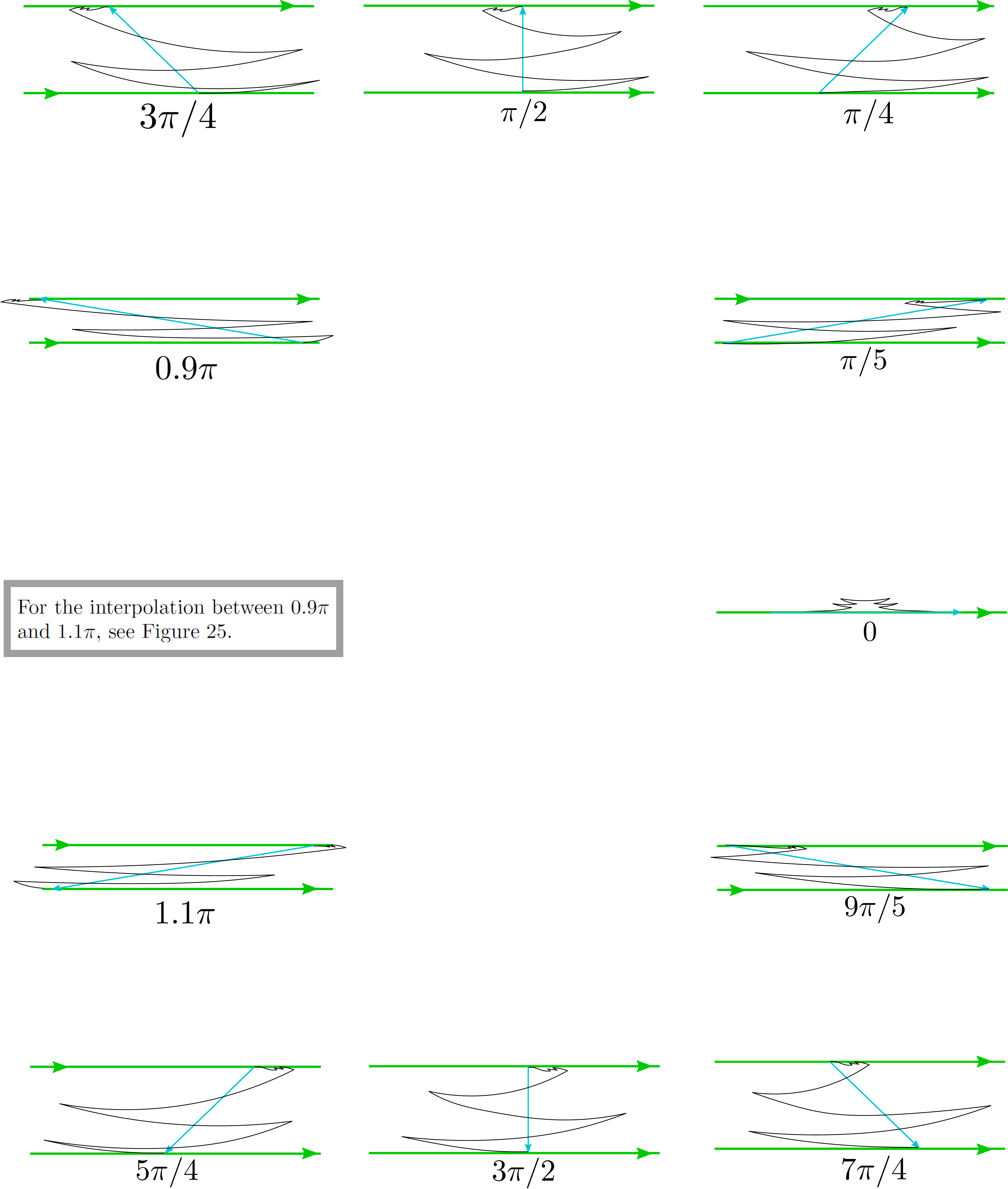}
	\centering
	\caption{Interpolation. Note that the orientation of the formal derivative (green) is determined by the Legendrian segments in $C$.}\label{reloj}
\end{figure}	

\begin{figure}[h!]
	\includegraphics[scale=0.3]{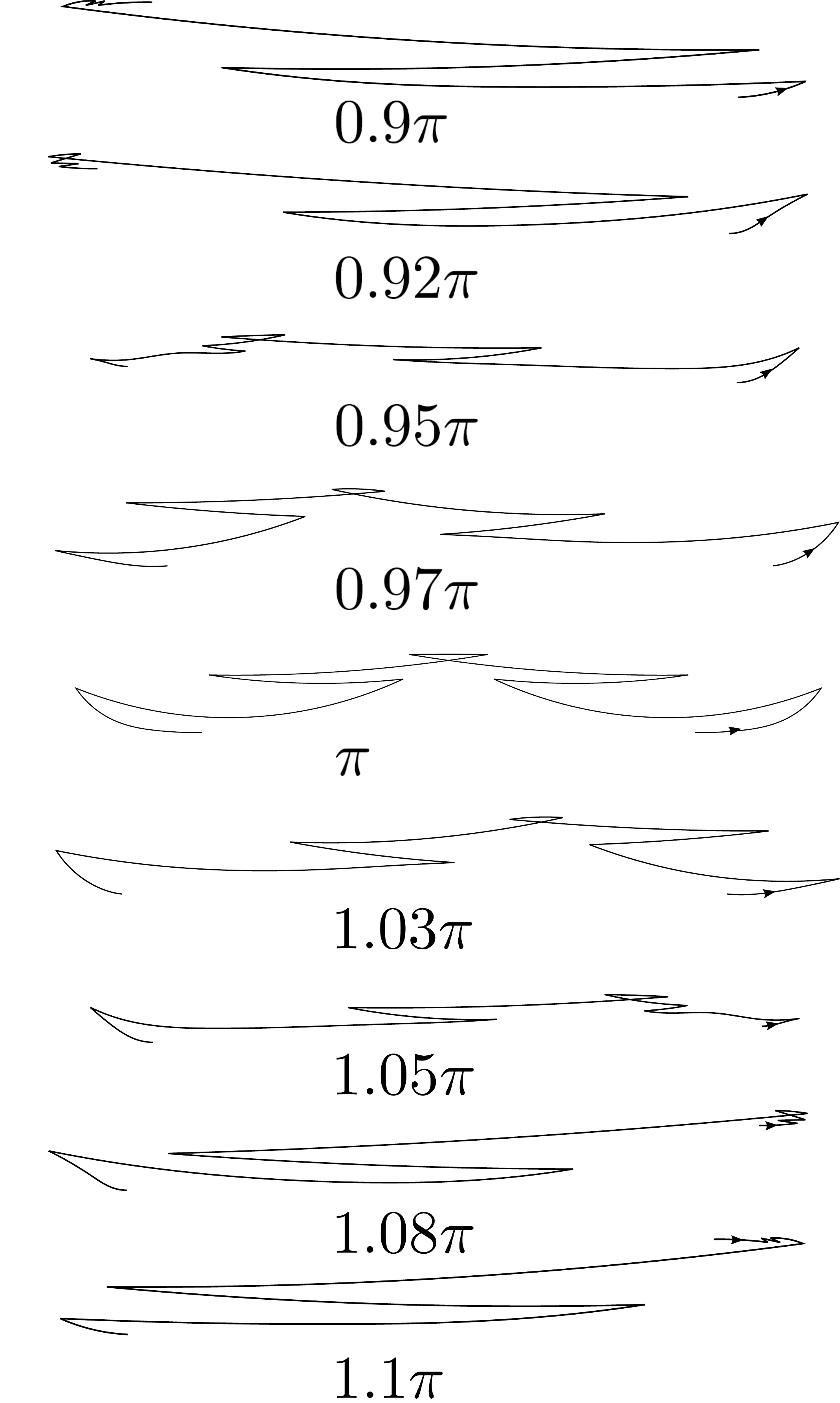}
	\centering
	\caption{Interpolation between $0.9\pi$ and $1.1\pi$ in Figure \ref{reloj}.}\label{InterpolationFina}
\end{figure}

The construction is clearly continuous in the angle, therefore it can be done in a canonical way. This readily implies that it works parametrically. Moreover, the approximation coincides with the formal derivative at the beginning and the end of the interval $I_i$, thus the constructed approximations are smooth. 

We claim that the $\gamma^{k,1}$ are embeddings. The segments $\gamma^{k,1}_{|I_i}$ are clearly embedded: just by visual inspection and since they just possess $R$ $II$ front moves. However, intersections may show up between the end of $\gamma^{k,1}_{|I_i}$ and the beginning of $\gamma^{k,1}_{|I_{i+1}}$. By sufficiently increasing $N$, we assume that the angle of the two consecutive segments is almost the same. In other words, we are in the same position in the clock figure (Figure \ref{reloj}). There are no tangencies between the two consecutive segments if the angle is away from $\pi$ radians since the front curve of the union is clearly embedded. In the neighborhood of $\pi$ radians, we need to be much more careful and depict the precise movie of the two consecutive segments. The problem is that there are intersections in front projection between the end of the first approximating curve and the beginning of the next one. In Figure \ref{pegado}, $\gamma^{k,1}_{|I_i}$ is depicted in red and $\gamma^{k,1}_{|I_{i+1}}$ in black.  The key point is to make sure that the big--sized stabilization boundary, i.e the beginning of the approximating curve for angle smaller than $\pi$ (respectively the end for angle bigger than $\pi$) has the shape of a scimitar with a fat blade. The fatness of the blade allows to fit inside its bell the whole end of the previous approximating curve (consisting of very small stabilizations) and the convexity of the bottom part of the scimitar blade allows to cross with $R$  $II$ moves to the exterior part of the scimitar. Once there, the small stabilizations are running away from the danger area. It is left to check the crossing of the first stabilization of the second curve with the last stabilization of the first one. They do not create tangencies thanks to the convexity of the scimitar shape.
 
\begin{figure}[h!]
	\includegraphics[scale=0.27]{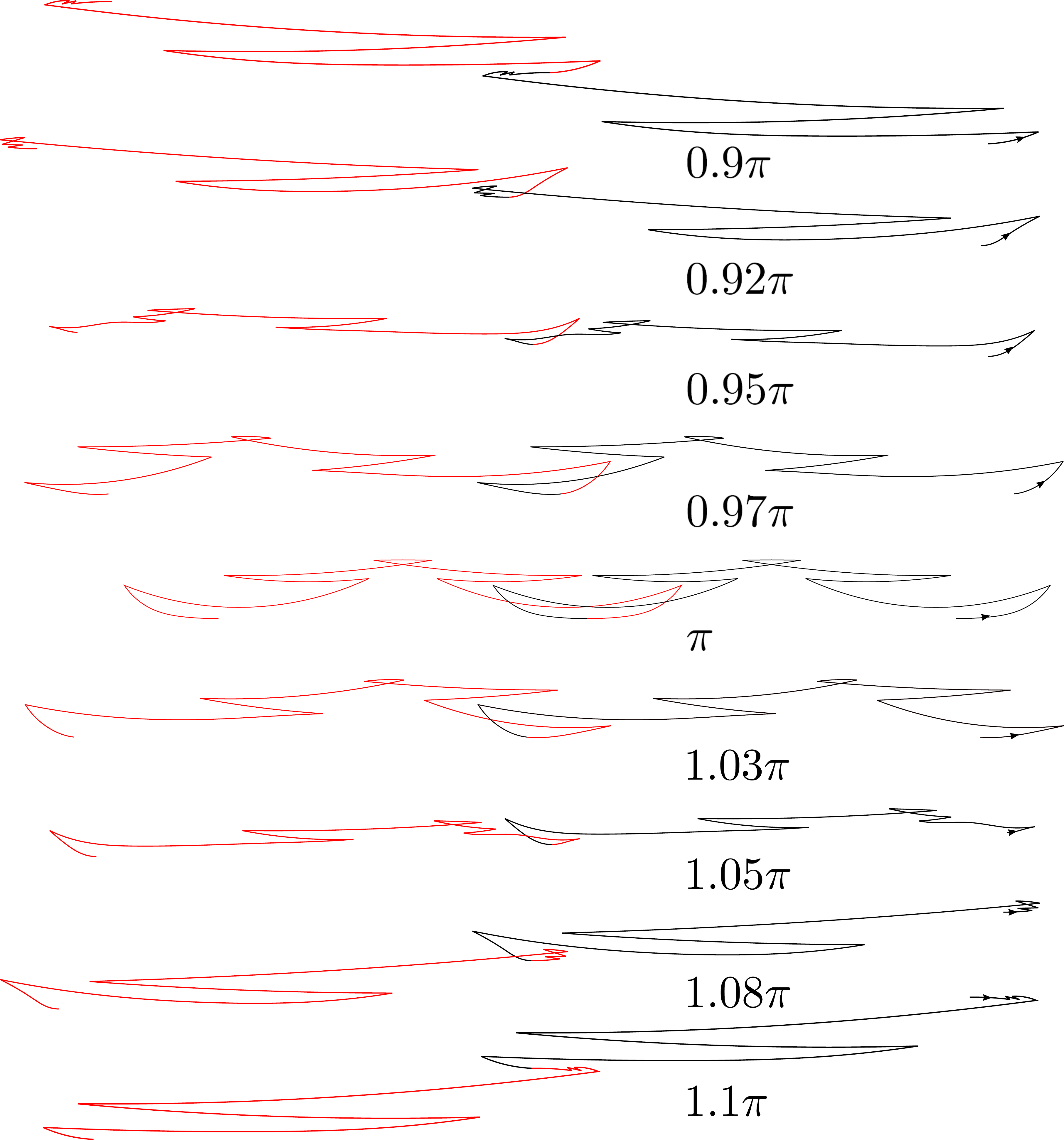}
	\centering
	\caption{Checking that the gluing process does not create tangency points.}\label{pegado}
\end{figure}

In the general case, apply Corollary \ref{CuspSameTime}, to assume that the cusp points of each $\gamma^k$, $k\in C$, are at times $c_1,\ldots,c_m$. It is very simple to deform each $\gamma^k$, $k\in K$, at small neighborhoods of times $c_1,\ldots, c_m$, making the family Legendrian and having a cusp point over each $c_1,\ldots, c_m$. Finally, add $c_0=0$ and $c_{m+1}=1$ and apply the first case over the family $\gamma^{k}_{|[c_i ,c_{i+1}]}$, $i=0,\ldots,m$, to conclude the proof.
\end{proof}

We are able to prove the following:

\begin{proof}[Proof of Lemma \ref{HorizontalAproximationOfLoops}.]
	Let $\gamma$ be a horizontal embedding such that $\gamma_G$ is a Legendrian embedding and let $\gamma^s\in\Emb(\NS^1,\R^4)$, $s\in[0,1]$, be any loop based at $\gamma$, ie $\gamma^0=\gamma^1=\gamma$. Without loss of generality we may assume that the Geiges projection of $\gamma^s$ is embedded except for a finite sequence of times $0<s_1<\cdots<s_N<1$ and the strict immersions $\gamma^{s_i}_{G}$ has a generic self--intersection. By Corollary \ref{cor:AproxHorizontalRot0} there exists a loop $\hat{\gamma}^s$ of horizontal embeddings with $\hat{\gamma}^0=\hat{\gamma}^1=\gamma$ homotopic to $\gamma^s$. Note that the total zero area condition for a closed horizontal curve can be easily achieved with Lemma \ref{DeltaBarquillo} since the map $A:[0,1]\rightarrow\R,s\mapsto\int (\hat{\gamma}^s)^*zdx$; is continuous. Moreover the approximation procedure fixes the cusps of $\gamma$ and, thus, the rotation number of this loop is zero so the result follows. 
\end{proof}

\subsection{Cancellation by pairs of zeroes in the same connected component.} \label{sub:same}
For a generic $2$--disk of Legendrians $\D$, the intersection of the disk with the space of strict immersions is a set of curves with crossings. They come from a set of curves in $\SSLegImm(\R^3)$ by projection. We declare a connected component of the set $\D \cap \SLegImm(\R^3)$ to be a connected component of the set seen in $\SSLegImm(\R^3)$. In other words, the two branches of a crossing (possibly) are in different connected components.

Let $\varphi:\D\rightarrow\HorImm(\R^4),z\mapsto\gamma^z$, be a disk of horizontal immersions satisfying:
\begin{itemize}
	\item $\varphi(\partial\D)\subseteq\Hor(\R^4)$,
	\item $\pi_G \circ \varphi:\D\rightarrow \LegImm(\R^3)$ is in generic position, 
	\item $(\pi_G\circ\varphi)^{-1}(\SLegImm(\R^3))$ is a finite union of connected components $C_1,\ldots,C_n$,
	\item $C_1$ contains, at least, two zeroes of the area function.
\end{itemize}
Note that Lemma \ref{lem:ThomTransversalityDiskOfHorizontalImmersions} implies that there is no loss of generality in thesew assumptions. Let $\bar{C}\subseteq C_1$ be a connected arc containing exactly two zeroes of the area function. Parametrize the curve $(\pi_G\circ\varphi)(\bar{C})$ as $\{\gamma^{s}_{G}:s\in[0,1]\}$ such that $\gamma^{\frac{1}{4}}_G$ and $\gamma^{\frac{3}{4}}_G$ are the two zeroes of the area function. Moreover, we can assume that $\{\gamma^{s}_G:s\in[0,\frac{1}{4}]\cup[\frac{3}{4},1]\}$,  is smooth (does not contain any cusp) in $\LegImm(\R^3)$. 

\begin{figure}[h]
	\includegraphics[scale=0.25]{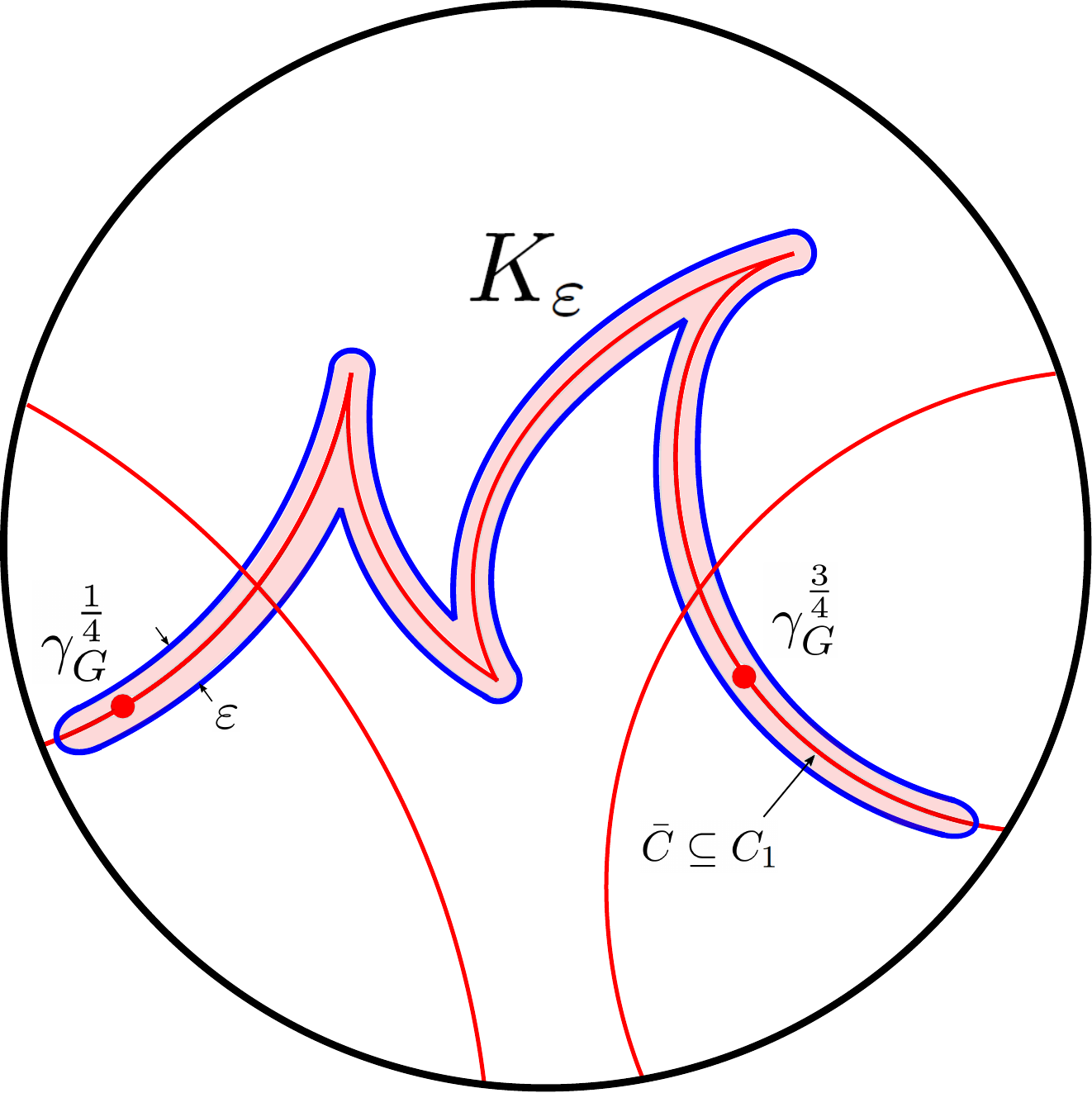}
	\centering
	\caption{Compact neighborhood $K_\varepsilon$ of $\bar{C}$ in $\D$.\label{entorno}}
\end{figure}

Consider $\varepsilon>0$ small enough  and the closed $\varepsilon$--neighborhood $K_\varepsilon$ of $\bar{C}$ in $\D$ (see Figure \ref{entorno}), satisfying that $K_{\varepsilon}$ does not contain any other zero of the area function for a different connected component and such that $K_\varepsilon \cap \partial \D = \emptyset$.
We have a well--defined continuous map $\bar{\Gamma}:\bar{C}\rightarrow\NS^1\times\NS^1, s\mapsto (t_{0}^{s},t_{1}^{s})$, where the components of the ordered pair $(t_{0}^{s},t_{1}^{s})$ are the tangency times of the front projection of $\gamma^{s}_G (t)=(x^s (t),z^s (t),w^s (t))^\intercal$, $s\in[0,1]$, and the order of the pair is determined by the following rule: $(t_{0}^{0},t_{1}^{0})$ are chosen with the usual rule to define the area function $\varepsilon_A$ and the pair of points $(t_{0}^{s},t_{1}^{s})$ are extended in the unique continuous way. In order to apply Lemma \ref{DeltaBarquillo} to cancel the zeroes $\gamma^{\frac{1}{4}}_G$ and $\gamma^{\frac{3}{4}}_G$ it is necessary to find an extension $\Gamma:K_\varepsilon\rightarrow\NS^1\times\NS^1,k\mapsto(t_{0}^{k},t_{1}^{k})$, of $\bar{\Gamma}$. Just take any retraction $r:K_\varepsilon\rightarrow \bar{C}$ and define $\Gamma:=r^* \bar{\Gamma}=\bar{\Gamma}\circ r$. Write $K=K_\varepsilon$.

We state the main result of this section.

\begin{proposition}\label{zeroesSameComponent}
	There exists a $1$--parametric family $\varphi_u :\D\rightarrow\HorImm(\R^4)$, $u\in[0,1]$, satisfying:
	\begin{itemize}
		\item [(i)] $\varphi_0 =\varphi$,
		\item [(ii)] $\varphi_u (z)=\varphi(z)$, $z\in\D\backslash K$,
		\item[(iii)] $(\pi_G \circ \varphi_u)^{-1}(\SLegImm(\R^3))=(\pi_G \circ \varphi)^{-1}(\SLegImm(\R^3))$,
		\item[(iv)] $\varepsilon_A(\pi_G(\varphi_1 (z)))\neq0$, $z\in\bar{C}\subseteq K$.
	\end{itemize}
\end{proposition}
\begin{proof}
Assume that the area function over $\gamma^{s}_G$ is positive in $s\in[0,\frac{1}{4})$.  Define \[\tilde{\varepsilon}_A (s)=\tilde{\varepsilon}_A(\tilde{\gamma}^{s}_G)=\int_{t_{0}^{s}}^{t_{1}^{s}} z^s (t) (x^s)'(t)dt.\] Observe that $\tilde{\varepsilon}_A$ is just the continuous version of the area function over $\bar{C}$ and satisfies $|\varepsilon_A|=|\tilde{\varepsilon}_A|$. The function $\tilde{\varepsilon}_A:\bar{C}\rightarrow\R$ is depicted in Figure \ref{grafoareafunction}.

\begin{figure}[h]
	\includegraphics[scale=0.14]{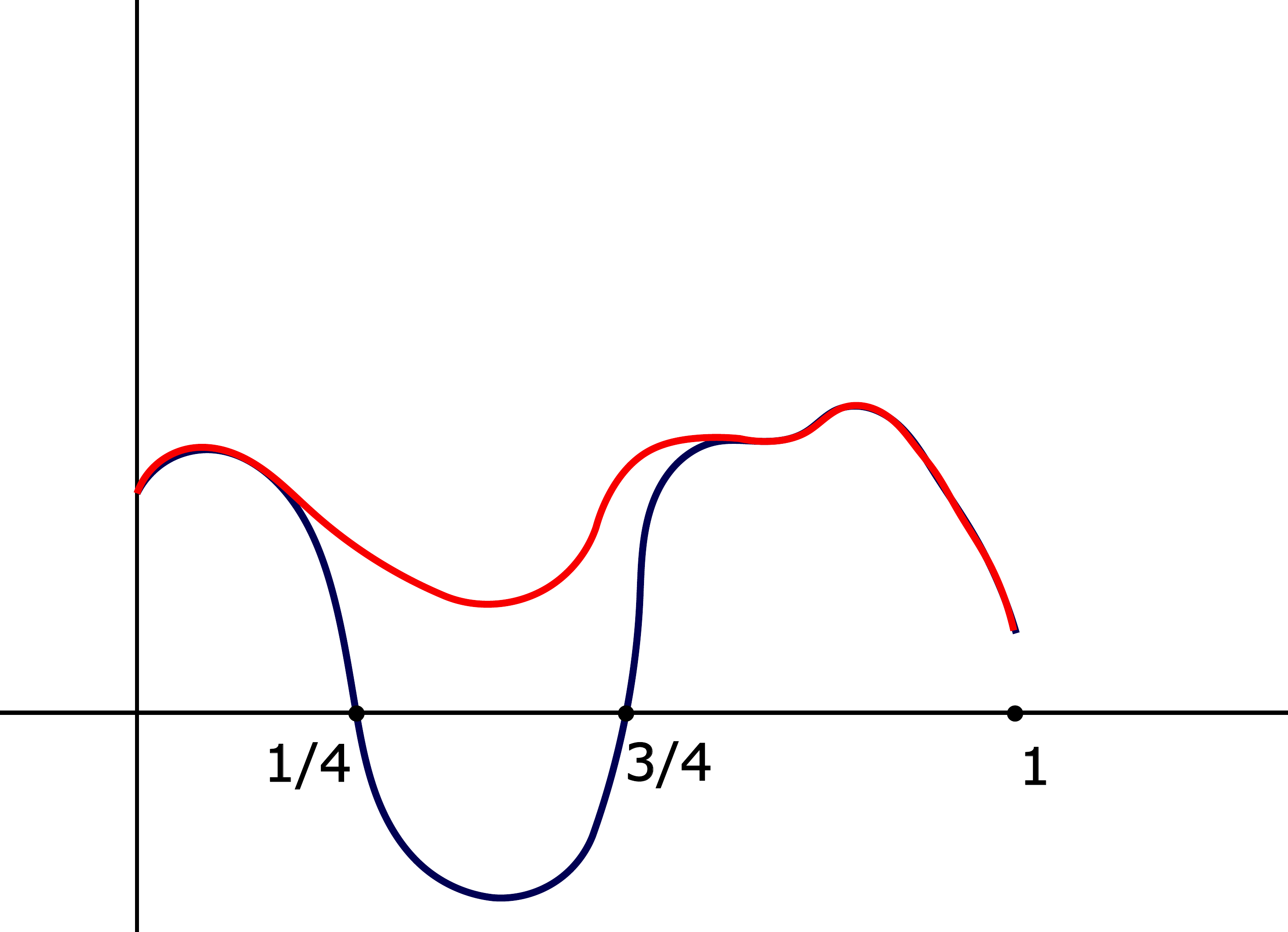}
	\centering
	\caption{Possible graph of $\tilde{\varepsilon}_A$ before (blue) and after (red) adding area using Lemma \ref{DeltaBarquillo}.\label{grafoareafunction}}
\end{figure}

Define $\varphi_u=\varphi$ over $\D\backslash K$. The problem is clearly to define $\varphi_u$ over $K$. Proceed as follows: In order to \em add area \em choose two points $p^k,n^k\in\NS^1$ in the domain of $\gamma^{k}_G$, $k\in K$. We require these points to have the following properties 
\begin{itemize}
	\item  $p^k$ and $n^k$ are regular for the front projection of $\gamma^{k}_{G}$,
	\item if $\gamma^{k}_{G}$ has more tangency points, different from $(t_{0}^{k},t_{1}^{k})$, the area function over these points remains unchanged after adding area over $p^k$ and $n^k$.
\end{itemize}
For a given $\gamma^{k}_G$, choose two regular points $p^k$ and $n^k$ very close to $t_{0}^{k}$, one to the left and another one to the right of $t_{0}^{k}$. These two points satisfy the properties above. By compactness, there exists a finite open covering $\{ U^i\}$ of $K$ and an election $(p^i,n^i)\in\NS^1\times\NS^1$ compatible with each $\gamma^{k}_G$, $k\in U^i$. Take a partition of unity $\{\phi^i\}$ subordinated to $\{U^i\}$ and a continuous function $\eta:K\rightarrow[0,1]$ such that
\begin{itemize}
	\item $\eta\equiv0$ in $\partial K$ and 
	\item $\eta\equiv1$ in $(\pi_G\circ\varphi)^{-1}(\{\gamma^{s}_G:s\in[\frac{1}{4},\frac{3}{4}]\})\cap K$. 
\end{itemize}
	Finally, fix some positive number $A>\max\{|\tilde{\varepsilon}_A (\tilde{\gamma}^{s}_G)|:s\in[0,1]\}$. Define $\varphi_u (k)$, $k\in K$, as the lift of
\[ \gamma^{k}_G\#_i \delta_{p^i} (u\phi^i (k)\eta(k)A)\#_i \delta_{n^i} (-u\phi^i (k)\eta(k)A),\]
to $\HorImm(\R^4)$. The family $\varphi_u$ satisfies the required properties.
\end{proof}

\section{The core of the proof}\label{sec:Core}

We explain now how for a given capping disk of horizontal immersions we are able to cancel zeroes of the area function by deforming the disk relative to the boundary. In particular, we prove that any pair of zeroes can be cancelled out without creating new zeroes: there will be no assumptions about the zeroes. This completes the proof of the main Theorem.

\subsection{Changing the diagram of curves of strict Legendrian immersions.} \label{sub:changedia}
The goal is to create a $2$--square ($2$--disk) in the space of $\LegImm(\R^3)$ such that the bottom of the square is a path $\gamma^s$ in the disk $\D_G$ whose endpoints are two strict Legendrian immersions $\gamma^0$ and $\gamma^1$. Moreover, there exists a finite set of values of the parameter  $0< s_1 < \cdots < s_\nu <1$ such that $\gamma^{s_j}$ is also a strict immersed Legendrian. The rest of the points $\gamma^s$, $s\in(0,1)\setminus \{ s_1, \ldots, s_{\nu} \}$ are Legendrian embeddings. The other three sides of the square will produce a path of strict Legendrian immersions that has the property that the self--intersections points $\hat\gamma^s(t_0)=\hat\gamma^s(t_1)$ form a continuous family. In other words, we are constructing a curve of Legendrians that makes the two initial Legendrians live in the same connected component of strict immersions. This is the key point and is the content of Proposition \ref{CreatingLegendrianTangency}. Afterwards, we play with the result: given an embedded curve in the disk whose ends live in two different connected components in the disk, we find a square transverse to the disk $\D_G$ that creates on the other three sides a curve of strict immersions joining the two strict Legendrian immersions. Then, we take a small open neighborhood of it removing the interior in order to produce a new disk (see Figure \ref{edificio}) in which the curves of strict Legendrian immersions have changed their topology (see Figure \ref{Surgery})\footnote{The picture looks like a bypass picture in contact topology, though geometrically it has nothing to do with a bypass. Formally, we play the same game: our initial disk plays the formal role of the convex surface, the strict immersions play the formal role of the dividing sets. The bypass changes the strict immersions/dividing set combinatorial structure.}. Proposition \ref{CreatingSmoothTangency} explains how to build the disk in the smooth category. Finally, we just use Theorem \ref{LegendrianAproximation} to conclude in the Legendrian setting.

Denote by $\SLegImmG (\R^3)\subseteq\SLegImm(\R^3)$ the open stratum of simple transverse self--intersections.

\subsubsection{Legendrian fingers.}

\begin{definition}
	Fix $\gamma\in\Emb(\NS^1,\R^3)$ and $b\in\NS^1$. A \em finger bone \em is an embedded path $\beta:[0,1]\rightarrow\R^3$ such that $\beta(0)=\gamma(b)$ and the intersection is transverse. Moreover, $\beta(0,1]\cap\gamma(\NS^1)$ is either one point or the empty set. In the former case the intersection must be transverse.
\end{definition}

\begin{definition} 
	Fix $\gamma\in\Emb(\NS^1,\R^3)$ and $b\in\NS^1$. A \em finger germ \em is a pair $(\beta,v)$ satisfying
		\begin{itemize}
			\item [(a)] $\beta:[0,1]\rightarrow\R^3$ is a finger bone such that $\beta(0)=\gamma(b)$,
			\item[(b)] $v\in\beta^* T\R^3$ satisfies
			\begin{itemize}
				\item [(i)] $v(t)$ and $\beta'(t)$ are linearly independent,
				\item[(ii)] $v(0)=\gamma'(b)$,
			\end{itemize}
		\end{itemize}
\end{definition}

Without loss of generality assume that $\gamma$ is geodesic on a very small neighborhood of $b$. For $\varepsilon>0$ small enough, fix an embedding $Q: [-\varepsilon, \varepsilon] \times [0,1] \to \R^3$ defined as $Q(x,z)= exp_{\beta(z)}( x \cdot v(z))$. Choose a bump function $\chi:[-\varepsilon,\varepsilon] \rightarrow[0,1]$ satisfying
\begin{itemize}
	\item {} $\chi$ is an even function,
	\item {} $\chi$ is increasing in $[-\varepsilon,0]$
	\item{} $\chi(-\varepsilon)=0$, $\chi^{(j)}(-\varepsilon)=0$ for all $j\in\Z^+$ and
	\item{} $\chi(0)=1$ and $\chi^{(j)} (0)=0$ for all $j\in\Z^+$.
\end{itemize}

\begin{definition}
	The \em finger deformation \em $\gamma\#(\beta,v)\in\Imm(\NS^1,\R^3)$ associated to the finger germ $(\beta,v)$ is given by the curve
$$ \gamma \# (\beta,v) (t) = \left\{ \begin{array}{cc} \gamma(t) & t\not\in (b-\varepsilon, b +\varepsilon), \\
Q(t-b, \chi(t-b))& t\in [b-\varepsilon, b +\varepsilon]. \end{array}\right. $$	
\end{definition}

This is a well--defined operation that works for compact families of curves and finger germs, provided that $\varepsilon>0$ is chosen small enough and uniform for the whole family.

 Note that $\gamma$ and $\gamma\#(\beta,v)$ are homotopic in $\Imm(\NS^1,\R^3)$. Denote by $(\beta_u,v_u)=(\beta_{|[0,u]},v_{|[0,u]})$ the obvious homotopy of germs of fingers. We say that $\gamma\#(\beta_u,v_u)$ is obtained from $\gamma$ \em adding the finger germ \em $(\beta_u,v_u)$. 

\begin{figure}[h]
	\includegraphics[scale=0.2]{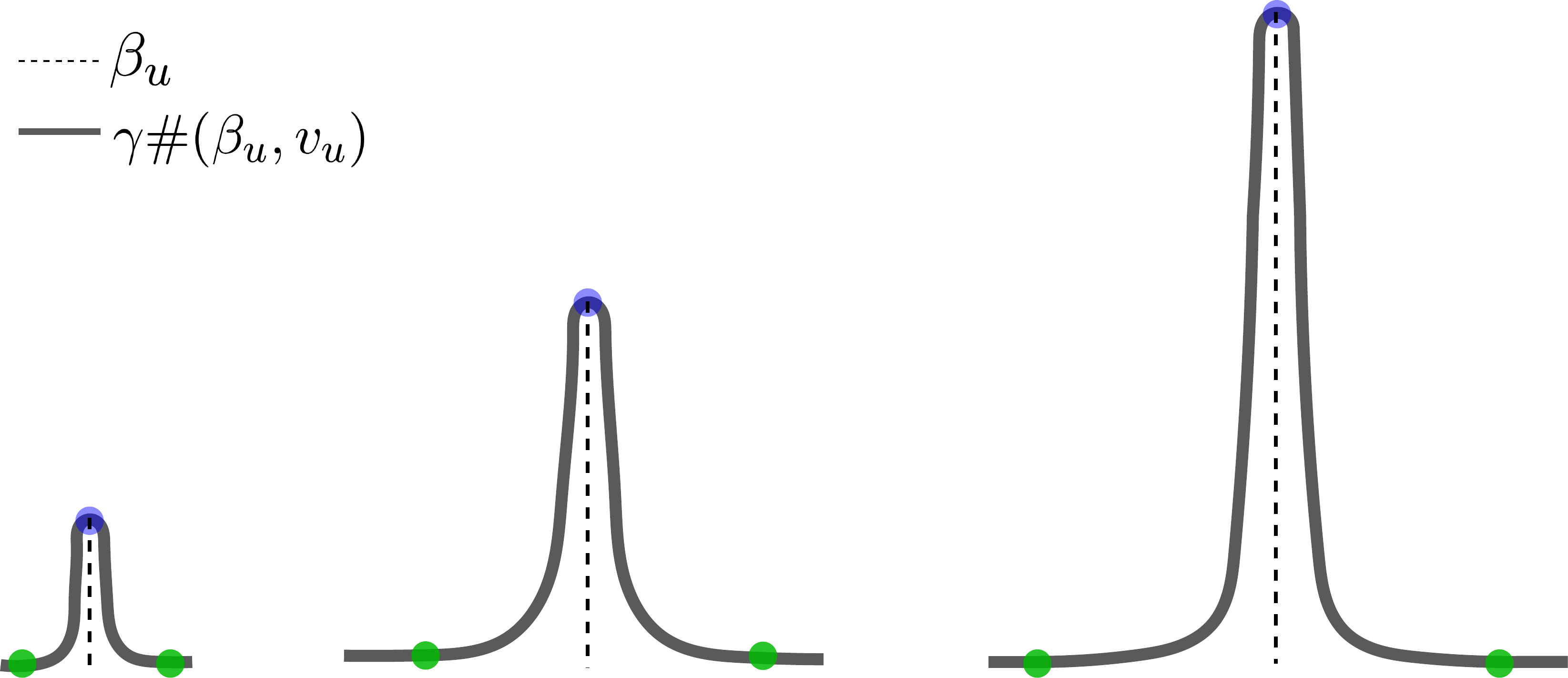}
	\centering
	\caption{Finger addition process $\gamma\#(\beta_u,v_u), u\in[0,1]$.}\label{CreacionFinger}
\end{figure}

\begin{figure}[h]
	\includegraphics[scale=0.7]{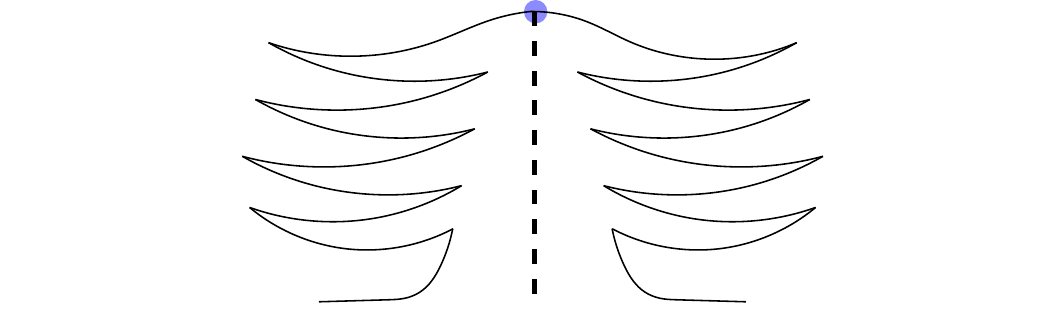} \label{FIngerTower}
	\centering
	\caption{Legendrian Finger construction process. The key tool is to create a huge {\em tower} of double stabilizations.}\label{CreacionFinger}
\end{figure}

\begin{definition}
Fix $\gamma\in\Leg(\R^3)$ and $(\beta,v)$ a finger germ such that
\begin{itemize}
	\item [(a)] $\beta$ is a transverse embedding, ie $\beta'(t)\notin\xi_{\beta(t)}$ for each $t\in[0,1]$, 
	\item[(b)]  $v$ is Legendrian, ie $v(t)\in\xi_{\beta(t)}$ for each $t\in[0,1]$.
\end{itemize}
We say that $(\beta,v)$ is a \em pre--Legendrian finger germ\em.
\end{definition}

Let us explain how to add the \em Legendrian finger \em associated to a pre--Legendrian finger germ. We are working on the standard $(\R^3(x,y,z),\xi=\ker(dz-ydx))$ with the trivialization of $\xi$ given by the Legendrian framing $\langle \partial_y \rangle$. Take a contactomorphism  $\Phi:(C_\varepsilon,\xi)\rightarrow (\D^2(0;\varepsilon)\times[0,1](x,y,z),\xi=\ker(dz-ydx))$. The domain $C_\varepsilon$, that we will call a cylindrical shape, is a sufficiently small tubular neighborhood of $\beta$ and fixing $\varepsilon>0$ small enough, we have that 
\begin{itemize}
	\item $\Phi \circ \beta (t)= (0,0,t)$,
	\item $\Phi_* v=\partial_x,$
	\item $\Phi \circ \gamma(t) = (t-b,0,0)$, for $|t-b|\leq \varepsilon$.
\end{itemize}
The existence of this contactomorphism is a direct consequence of the \em Tubular Neighborhood Theorem for Contact Submanifolds \em (see \cite[Theorem 2.5.15]{GeigesCont}).  From now on we work in these coordinates. The \em front projection \em in $(\D^2(0;\varepsilon)\times[0,1](x,y,z),\xi=\ker(dz-ydx))$ is given by the projection in the $XZ$--plane. The front projection of the original Legendrian in a neighborhood of the base of the finger $\gamma_{|\Op(\{b\})}$ is the horizontal axis in front coordinates. The front of the transverse finger bone $\beta$ is the  vertical axis, see  Figure \ref{CreacionFinger}. To approximate the finger bone by a Legendrian finger just add $N=N(\varepsilon)\in\Z_{\geq0}$ DSs to the Legendrian $\gamma$ over $\gamma(b)$ as in Figure \ref{CreacionFinger} in such a way that the derivative of the new Legendrian $\gamma\#(\beta,v)$ at time $b$ is $v(1)$. 

This operation is well defined up to Legendrian immersions homotopy and works for compact families. Note that this follows just taking the radius $\varepsilon$ of the tubular neighborhood and the number $N$ of DSs in the finger addition process constant in the whole family. This means that given a family $(\gamma^k,(\beta^k,v^k),C_{\varepsilon}^{k})$, $k\in K$, we obtain a continuous family of Legendrian finger deformations $\gamma^{k}\#(\beta^k,v^k)$. In particular, if there is some $k_0\in K$ such that the pre--Legendrian finger germ $(\beta^{k_0},v^{k_0})$ is trivial, ie $$(\beta^{k_0},v^{k_0})\equiv(\gamma^{k_0}(b),(\gamma^{k_0})'(b)),$$ and $N\neq 0$, then the Legendrians  $\gamma^{k_0}$ and  $\gamma^{k_0}\#(\beta^{k_0},v^{k_0})$ are different. In other words, to interpolate between a trivial pre--Legendrian finger germ and a non--trivial one we need to add DSs to the trivial finger. Another important observation is that the DSs operation is well--defined for horizontal embeddings (see Lemma \ref{DeltaSeta}), this is relevant since we are going to work with the Geiges projection of families of horizontal curves.

\begin{remark}
Another possible definition of Legendrian finger is by using Reidemeister I moves instead of DSs. Note that, in both cases, the transverse finger bone is approximated by \em half of the Legendrian finger \em which is a very stabilized Legendrian curve. Moreover, the definition using Reidemeister I moves does not work well for our purposes because the two \em halves \em of the finger are linked and in the later argument is very important to not have linking between the two halves.
\end{remark}

\begin{figure}[h]
	\includegraphics[scale=0.5]{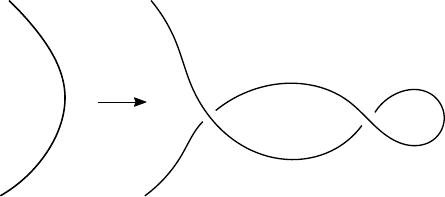}
	\centering
	\caption{Stabilization of a transverse curve.}\label{fig:tstab}
\end{figure}

A stabilization of a transverse curve is an operation that changes  the curve  in such a way that in front projection a small arc is replaced by a $2\pi$ rotation as in Figure \ref{fig:tstab}. The stabilization does not depend on the chosen small segment (see Etnyre \cite{Etnyre}). So we can speak of a $k$ stabilization of a transverse curve $\beta_0$ that we will denote by $\beta_k$.

The next Lemma will be useful in the next subsection and ensures that we can interpolate between a small finger with finger bone a small Reeb chord and another one with finger bone a very stabilized Reeb chord.

\begin{lemma}\label{lem:ReebChordStabilized}
	Let $\gamma\in\SLegImmG(\R^3)$ be a Legendrian immersion with a self--intersection at times $t_0,t_1\in \NS^1$ and $\tilde{\gamma}\in\Leg(\R^3)$ a Legendrian embedding such that
	\begin{itemize}
		\item [(i)] There is a small Reeb chord $\beta:[0,1]\rightarrow\R^3$ between $\tilde{\gamma}(b)$ and $\tilde{\gamma}(e)$ for some times $b,e\in\NS^1$;
		\item[(ii)] $\gamma=\tilde{\gamma}\#(\beta,v)$ where $v=\partial_x+y\partial_z$.\footnote{We may assume that $\tilde{\gamma}'(b)=v_{\gamma(b)}$.}
	\end{itemize}
Then, for any $k\in\Z^+$ there exists a $1$--parametric $\gamma^u\in\SLegImmG(\R^3)$, $u\in[0,1]$, with a self--intersection at times $t_0,t_1\in\NS^1$, such that $\gamma^0=\gamma$ and $\gamma^1=\tilde{\gamma}\#(\beta_k,v)$. Moreover, the family is constructed from $\gamma$ just by adding a sequence of DSs to the initial finger, ie the number $N$ of DSs in the definition of Legendrian finger can be chosen in such a way that $\tilde{\gamma}\#(\beta,v)=\tilde{\gamma}\#(\beta_k,v)$. 
\end{lemma}
\begin{proof}
Since $\beta_k$ is $C^0$--close to $\beta$ we may assume that the finger associated to $(\beta_k,v)$ is contained in the tubular neighborhood $C_\varepsilon$ of $\beta$. Moreover, $(C_\varepsilon,\xi)$ is contactomorphic to $(\R^3,\xi=\ker(dz-ydx))$ by Eliashberg's classification of contact structures in $\R^3$ \cite{EliashbergR3}. In this sense the right side of the fingers associated to $(\beta,v)$ and $(\beta_k,v)$, denoted by $R$ and $R_k$ respectively, can be understood as positively stabilized \em long \em Legendrian unknots. In the same way the left side, denoted as $L$ and $L_k$, as negatively stabilized long Legendrian unknots. Thus, there exist two non--negative numbers $M,M_k\in\Z_{\geq0}$ such that 
\begin{itemize}
	\item $S_{+}^{M}(R)$ is Legendrian isotopic to $S_{+}^{M_k}(R_k)$,
	\item $S_{-}^{M}(R)$ is Legendrian isotopic to $S_{-}^{M_k}(R_k)$;
\end{itemize}
where $S_+$ ($S_-$) denotes a positive (negative) stabilization of a Legendrian curve. This follows from the classification of Legendrian unknots given in \cite{EliashbergFraser}. The Lemma follows from the fact that the number of positive and negative stabilizations are the same.
\end{proof}

\begin{remark}
	Another possible proof of the previous Lemma is to use Theorem \ref{LegendrianAproximation} to each half of the finger and, then, slide the DSs over the finger nails.
\end{remark}

A Legendrianization of a transverse knot/curve (with fixed framings at the ends) is an induced Legendrian knot/curve that is created by drawing the characteristic foliation of the boundary of a small solid torus fixed around the transverse knot for a radius such that the slope is rational and pick one integral curve. The transverse push off of a Legendrianization is always transverse isotopic to the original one (see \cite{Etnyre}).

The next lemma ensures the existence of $2$--parametric families of pre--Legendrian finger germs, it generalizes the  $1$--parametric result explained in Lemma 3.3.3 in \cite{GeigesCont}. It is the key lemma that we will use to produce controlled deformations on our families of curves.

\begin{lemma}\label{TransverseFingers}
	Fix $\beta^{s,u}\in\Emb([0,1],\R^3)$, $(s,u)\in[0,1]\times[0,1]$, assume that $\beta^{s,0}$ is transverse and $\beta^{s,u} _{|\Op(\{0,1\})}$ is also transverse, then there exist a constant $k>0$ and a family of transverse embeddings $\tilde{\beta}^{s,u}:[0,1]\rightarrow\R^3$, $C^0$--close to $\beta^{s,u}$, such that $\tilde{\beta}^{s,0}=\beta_k^{s,0}$ and $\tilde\beta^{s,u} _{|\Op(\{0,1\})} = \beta^{s,u} _{|\Op(\{0,1\})}$.
\end{lemma}

 \begin{proof}
Perform a Legendrianization of $\beta^{s,0}$ that we call $\gamma^{s,0}$, that extends to a family $\gamma^{s,u}$ that is smoothly isotopic to $\beta^{s,u}$. Further asume that we Legendrianize the small segments $\beta^{s,u} _{| t \in \mathcal {O}p(0) \cup \mathcal {O}p(1)}$. Abusing  notation, we still call $\gamma^{s,u}$ a perturbation of the original family in order to make it Legendrian on a very small neighborhood of a time $t_0$.

Now, apply Theorem \ref{LegendrianAproximation} to the two segments $\gamma^{s,u}_{|t\in[0, t_0]}$  and $\gamma^{s,u}_{|t\in[t_0,1]}$. This produces a family $\hat\gamma^{s,u}$ that is Legendrian. The process introduces a high number of DSs on the family $\gamma^{s,0}$. We denote by $k$ the number of introduced DSs,.

Construct a new family  $\tilde{\gamma}^{s,0}$ that it is just built by performing $k$ DS at time $t_0$ on $\gamma^{s,0}$. It is well known that DSs can be moved around a Legendrian (they are Legendrian homotopic \cite{Etnyre}).  This, in particular implies, that $\tilde\gamma^{s,0}$ and $\hat\gamma^{s,0}$ are homotopic through Legendrian embeddings $\gamma^{s,0}_r$, $r\in[0,1]$, ie $\gamma^{s,0}_0=\tilde\gamma^{s,0}$ and $\gamma^{s,0}_1=\hat\gamma^{s,0}$. So, we have proven that 
$$ \overline\gamma^{s,u} = \left\{ \begin{array}{cc}
\gamma^{s,0}_{2u} &, u \in [0, \frac12] \\
\hat\gamma^{s,2u-1} &, u\in [\frac12 ,1]
\end{array} \right.
$$
is a family of Legendrians approximating $\gamma^{s,u}$ such that $\tilde{\gamma}^{s,0}= \overline\gamma^{s,0}$.

Perform a transverse push off to obtain $\tilde\beta^{s,u}$ a family of transverse knots. Since the transverse push off of $k$ DSs of a Legendrian knot corresponds to $k$ stabilizations and since the transverse push off of a Legendrianization coincides with the original transverse knot, we have shown that 
$\tilde\beta^{s,0}= \beta^{s,0}_k$.
 \end{proof} 

The \em unknotting number \em of an embedding is finite (see Adams \cite{Adams}). We can think of the unknotting number as the number of fingers that we need to add to the knot in order to make it the unknot. Moreover the necessary number of fingers in order to unlink two embeddings is also finite. We summarize this information as 

\begin{proposition}\label{UnknottingSmooth}
	Fix $\gamma,\eta\in\Emb(\NS^1,\R^3)$ disjoint embeddings. There exists a parametrized unknot $\tilde{\gamma}\in\Emb(\NS^1,\R^3)$, obtained from $\gamma$ just by adding a finite number of finger germs, which is unlinked from $\eta$. 
\end{proposition}

We are interested in applying this result to the two \em resolutions \em of a strict immersion with exactly one (generic) self--intersection. More precisely, let $\gamma\in\Imm(\NS^1,\R^3)$ be an immersion with exactly one self--intersection at times $(t_0,t_1)\in\NS^1\times\NS^1$. Assume without loss of generality that there exists $\varepsilon>0$ small enough such that $\gamma((t_0-\varepsilon,t_0+\varepsilon)\cup(t_1-\varepsilon,t_1+\varepsilon))$ is contained in the affine plane $\gamma(t_0)+\Span\{\gamma'(t_0),\gamma'(t_1)\}$. Take coordinates $(x,y)\in\R^2$ in the plane $\gamma(t_0)+\Span\{\gamma'(t_0),\gamma'(t_1)\}$ such that $\gamma(t_0+r)=(r,0)$ and $\gamma(t_1+r)=(0,r)$ for $r\in(-\varepsilon,\varepsilon)$. Finally, fix any smooth non increasing function $\chi:[-\varepsilon,\varepsilon]\rightarrow[0,1]$ such that $\chi(-\varepsilon)=1$ and $\chi(\varepsilon)=0$.

\begin{definition}
	The \em resolution \em of $\gamma$ over $[t_0,t_1]$ is the embedding $\gamma_{\Res(t_0,t_1)}\in\Emb(\NS^1,\R^3)$ defined as:
	
	\begin{equation*}
		\gamma_{\Res(t_0,t_1)}(t)=\begin{cases}
			\gamma(t) & \text{ if $t\in(t_0+\varepsilon,t_1-\varepsilon)$,} \\
			\chi(t-t_0)(0,t-t_0)+(1-\chi(t-t_0))(t-t_0,0)& \text{ if $t \in[t_0,t_0+\varepsilon),$}  \\
			\chi(t-t_1)(0,t-t_1)+(1-\chi(t-t_1))(t-t_1,0)& \text{ if $t \in(t_1-\varepsilon,t_1].$} 
			
		\end{cases}
	\end{equation*} 
where $\NS^1=[t_0,t_1]/\sim$.
\end{definition}

\begin{remark}
	The resolution 
$\gamma_{\Res(t_0,t_1)}$ is well--defined up to isotopy. Moreover, observe
	that $\gamma_{\Res(t_0,t_1)}(\NS^1)\cap\gamma_{\Res(t_1,t_0)}(\NS^1)=\emptyset$.
	
	In these terms, Proposition \ref{UnknottingSmooth} applied over $\gamma$ just says that, up to adding a finite number
	of fingers to $\gamma$, we may assume that the two resolutions $\gamma_{\Res(t_0,t_1)}$ and $\gamma_{\Res(t_1,t_0)}$ are unlinked and one of them is an unknot. In fact, we can add fingers in order two obtain two unlinked unknots.
\end{remark}

In the same fashion we can define the resolution of a strict Legendrian immersion $\gamma\in\SLegImmG(\R^3)$. Assume that $\gamma(t_0)=\gamma(t_1)$. Take adapted coordinates such that the front $\gamma_F (t_0+r)=(r,0)$ and $\gamma_F (t_1+r)=(\pm r,r^2)$ for $r\in(-\varepsilon,\varepsilon)$. Fix a non increasing smooth function $\chi:[-\varepsilon,0]\rightarrow[0,1]$ such that $\chi(-\varepsilon)=1$ and $\chi(0)=0$. Finally, take a decreasing smooth function $\eta:[-\varepsilon,0]\rightarrow[0,\delta]$, where $\delta>0$ is very small, satisfying that $\eta(-\varepsilon)=\delta$ and $\eta^{j)}(0)=0$ for $j\in\Z^+\cup\{0\}$.

\begin{definition}\label{def:resolutions}
	The \em resolution \em of $\gamma$ over $[t_0,t_1]$ is the Legendrian embedding $\gamma_{\Res(t_0,t_1)}\in\Leg(\R^3)$ defined as follows:
	\begin{itemize}
		\item [(i)] If $\gamma_F(t_1+r)=(r,r^2)$ for $r\in(-\varepsilon,\varepsilon)$ then
		\begin{equation*}
			\gamma_{\Res(t_0,t_1),F}(t)=\begin{cases}
				\gamma_F(t) & \text{ if $t\in[t_0,t_1-\varepsilon)$,} \\
				(r,\chi(r)r^2+(1-\chi(r))\eta(r))& \text{ if $t=t_1+r\in(t_1-\varepsilon,t_1]$,} 
			\end{cases}
		\end{equation*} 
		\item [(ii)] If $\gamma_F(t_1+r)=(-r,r^2)$ for $r\in(-\varepsilon,\varepsilon)$ then
	\begin{equation*}
		\gamma_{\Res(t_0,t_1),F}(t)=\begin{cases}
			\gamma_F(t) & \text{ if $t\in[t_0,t_1-\varepsilon)$,} \\
			(-r,\chi(r)r^2+(1-\chi(r))r^{\frac{3}{2}})& \text{ if $t=t_1+r\in(t_1-\varepsilon,t_1]$,} 
		\end{cases}
	\end{equation*} 
	\end{itemize}
where $\NS^1=[t_0,t_1]/\sim$.
\end{definition}

	\begin{figure}[h]
		\centering
		\subfloat[Resolution  (i).]{
			\label{fig:resolution1}
			\includegraphics[width=0.4\textwidth]{resolution2.pdf}}\phantom{ aaaaa  }
		\subfloat[Resolution (ii).]{
			\label{fig:resolution2}
			\includegraphics[width=0.43\textwidth]{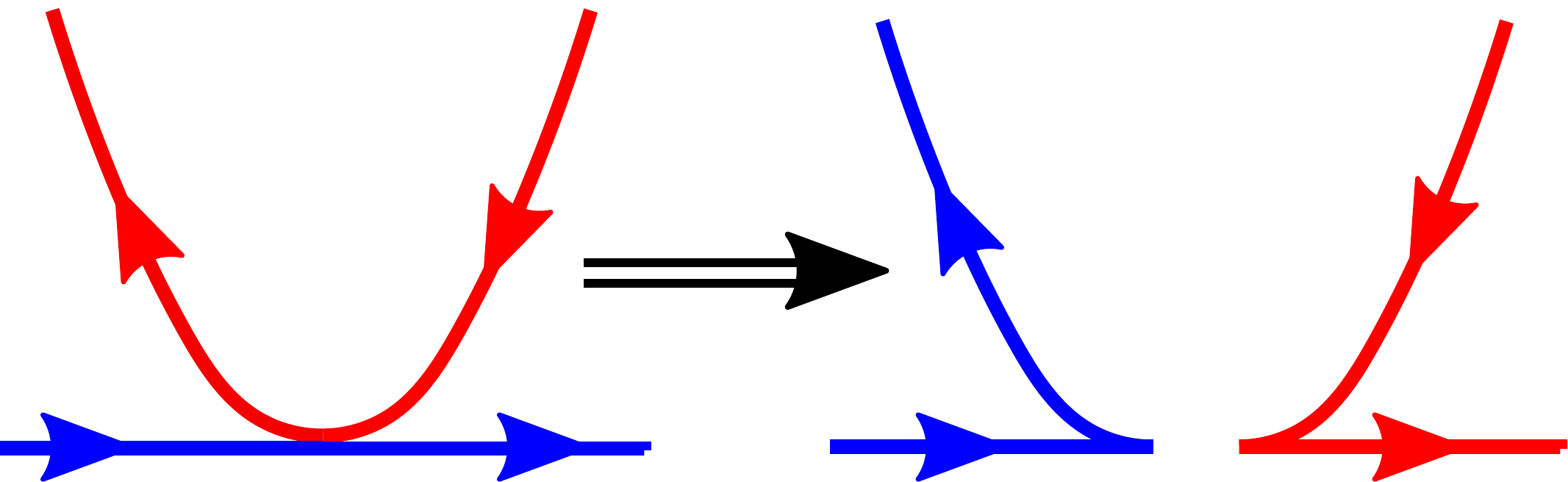}}
		\caption{Resolutions in Definition \ref{def:resolutions}. $\gamma_{\Res(t_0, t_1)}$ is represented in red.}
	\end{figure}

We have the following result for $1$--parametric families of smooth curves.

\begin{proposition}\label{CreatingSmoothTangency}
Fix $\gamma^s \in\Imm(\NS^1,\R^3)$, $s\in[0,1]$, such that 
\begin{itemize}
	\item [(i)] $\gamma^s$, $s\in(0,1)\setminus \{s_1,\ldots,s_\nu, 1\}$, for $\nu\geq 0$, where $0<s_1<\cdots<s_\nu<1$, is an embedding,
	
	\item[(ii)] $\gamma^i$, $i\in\{0,1\}$, has exactly one (generic) self--intersection at times $(t_{0},t_{1})\in\NS^1\times\NS^1$ and $\gamma^{s_j}$, $j\in\{s_1, \ldots, s_\nu \}$, has exactly one (generic) self--intersection at times $(\tau_{s_j}^0,\tau_{s_j}^1)\in\NS^1\times\NS^1$. The times $t_0,t_1,\tau_{s_j}^{0},\tau_{s_j}^{1}$, $j\in\{1,\ldots,\nu\}$, are pairwise different.
\end{itemize}

Then, there exist $k$ $1$--parametric families $\{\beta^{s}_{1},\ldots,\beta^{s}_{k}\}$ of finger germs and a family $\gamma^{s,u}\in\Imm(\NS^1,\R^3)$, $(s,u)\in[0,1]\times[0,1]$, such that $\gamma^{s,0}=\gamma^s$ and $\gamma^{s,1}=\gamma^s\#\beta^{s}_{1}\#\cdots\#\beta^{s}_{k}$, ie the $1$--parametric family $\gamma^{s,1}$ is obtained from $\gamma^s$ by adding in consecutive times the finite sequence of fingers, satisfying
\begin{itemize}
	\item[(i)] $\gamma^{s,1}$, $s\in[0,1]$, has one self--intersection at times $(t_0,t_1)\in\NS^1\times\NS^1$. This self-intersection is unique for $s\in[0,1]\backslash\{s_1,\ldots,s_\nu\}$.
	
	\item [(ii)] $\gamma^{i,u}$, $i\in\{0,1\}$, has one self--intersection at times $(t_0,t_1)\in\NS^1\times\NS^1$ and $\gamma^{s_j,u}$, $j\in\{1,\ldots,\nu\}$, has one self--intersection at times $(\tau_{s_j}^{0},\tau_{s_j}^{1})\in\NS^1\times\NS^1$.
	
	\item[(iii)] There exists a finite sequence $0<u_1<\ldots<u_N <1$ such that $$\gamma^{s,u_j}\in\Imm(\NS^1,\R^3)\backslash\Emb(\NS^1,\R^3)$$ for each $j\in\{1,\ldots,N\}$ and $s\in(0,1)$,
	
	\item[(iv)] $\gamma^{s,u}\in\Emb(\NS^1,\R^3)$ for $u\neq 1, u_1,\ldots,u_N$ and $s\in(0,1)\backslash\{s_1,\ldots,s_\nu\}$. 
\end{itemize}

Moreover, if there is $c\in(t_0,t_1)$ such that $\gamma^{s}_{|\Op(\{c\})}$ is constant then the construction can be done relative to $\Op(\{\gamma^{s}(c)\})$; ie $\gamma^{s,u}_{|\Op(\{c\})}$ is constant.
\end{proposition}

Note that we start with a segment $\gamma^s$ and we create a square $\gamma^{s,u}$ in such a way that one of the edges is the original segment and the other three edges constitute a curve of self--intersections. Therefore, this is the key result in the smooth category and we will adapt it to the Legendrian case in Proposition \ref{CreatingLegendrianTangency}.

\begin{proof}[Proof of Proposition \ref{CreatingSmoothTangency}]. \\

{\bf \em Case 1: $\{ s_1, \ldots, s_{\nu} \} = \emptyset$. }		
		
		The proof is divided in three steps.
	
\textit{Step I. Unknot $\gamma^{0}_{\Res(t_0,t_1)}$:} For $\gamma^0$ and  $\gamma^1$, we declare $t_0$ to be the point in the lower branch in front projection.\footnote{Being precise we should write $t_0(0)$ and $t_0 (1)$ since the intersection times of $\gamma^{0}$ and $\gamma^{1}$ could be different. We avoid this to simplify the notation.} Apply Proposition \ref{UnknottingSmooth} and create a sequence of finger germs $$(\beta^{0}_{1},v^{0}_{1}),\ldots,(\beta^{0}_{K},v^{0}_{k})$$ to unknot and unlink $\gamma^{0}_{\Res(t_0,t_1)}$ from $\gamma^{0}_{\Res(t_1,t_0)}$. Let $(b^{0}_{1},e^{0}_{1})\in\NS^1\times\NS^1$ be the base point and the end point of the finger $(\beta^{0}_{1},v^{0}_{1})$; ie the finger germ $(\beta^{0}_{1},v^{0}_{1})$ joins $\gamma^0 (b^{0}_{1})$ with $\gamma^0 (e^{0}_{1})$. We extend this finger germ to a whole family of finger germs for the family $\gamma^s$. In order to do this, choose any continuous map $F:[0,1]\rightarrow\NS^1\times\NS^1$, $s\mapsto (b^{s}_{1},e^{s}_{1})$; such that $b^{s}_{1}\neq e^{s}_{1}$ and $b^{1}_{1},e^{1}_{1}\neq t_i$, $i\in\{0,1\}$, this provides a continuous family $(\beta^{s}_{1},v^{s}_{1})$ of finger germs joining $\gamma^s(b^{s}_{1})$ with $\gamma^s (e^{s}_{1})$. This produces a family $\gamma^{s,u}$, with $s\in[0,1]$ and $u\in[0,\frac{1}{3K}]$, where $\gamma^{s,\frac{1}{3K}}=\gamma^{s}\#\beta^{s}_{1}$.

Proceed in the same way with the rest of the finger germs to construct a family $\gamma^{s,u}$, $s\in[0,1]$ and $u\in[0,\frac{1}{3}]$, satisfying
	\begin{itemize}
		\item[(i)] $\gamma^{s,0}=\gamma^s$,
		
		\item[(ii)] $\gamma^{s,\frac{1}{3}}=\gamma^s \# \beta^{s}_{1}\#\cdots\#\beta^{s}_{K}$,
		
		\item [(iii)] $\gamma^{i,u}$, $i\in\{0,1\}$, has one self--intersection at times $(t_0,t_1)\in\NS^1\times\NS^1$,
		
		\item[(iv)] There exists a finite sequence $0<u_1<\cdots<u_K<\frac{1}{3}$ such that $$\gamma^{s,u_j}\in\Imm(\NS^1,\R^3)\backslash\Emb(\NS^1,\R^3)$$ for each $j\in\{1,\ldots,K\}$ and $s\in(0,1)$,
		
		\item[(v)] $\gamma^{s,u}\in\Emb(\NS^1,\R^3)$ for $u\neq u_1,\ldots,u_K$ and $s\in(0,1)$,
		
		\item[(vi)] $\gamma^{0,\frac{1}{3}}_{\Res(t_0,t_1)}$ is a parametrized unknot and is unlinked from $\gamma^{0,\frac{1}{3}}_{\Res(t_1,t_0)}$.
		
	\end{itemize}

\textit{Step II. Unknot $\gamma^{1,\frac{1}{3}}_{\Res(t_0,t_1)}$:} Do the same as in Step I but be careful in order to add families of fingers which restricted to $s=0$ end up in $\gamma^{0,\frac{1}{3}}_{\Res(t_1,t_0)}$ and to not relink this resolution from the other resolution $\gamma^{0,\frac{1}{3}}_{\Res(t_0,t_1)}$. More precisely, let $(\beta^{1}_{K+1},v^{1}_{K+1}),\ldots,(\beta^{1}_{N},v^{1}_{N})$ be the sequence of finger germs obtained from Proposition \ref{UnknottingSmooth} to unknot and unlink $\gamma^{1,\frac{1}{3}}_{\Res(t_0,t_1)}$ from $\gamma^{1,\frac{1}{3}}_{\Res(t_1,t_0)}$. As in the previous step, begin with $(b^{1}_{K+1},e^{1}_{K+1})\in\NS^1\times\NS^1$ the base and the end points of the finger $(\beta^{1}_{K+1},v^{1}_{K+1})$ and extend this finger to the whole family via a map $F:[0,1]\rightarrow\NS^1\times\NS^1$, $s\mapsto (b^{s}_{K+1},e^{s}_{K+1})$. We need to make sure that the finger germ $(\beta^{1}_{K+1},v^{1}_{K+1})$ is extended to a $1$--parametric family $(\beta^{s}_{K+1},v^{s}_{K+1})$ such that 
\begin{itemize}
	\item {} $\beta^{0}_{K+1} (0)=\gamma^{0,\frac{1}{3}}(b^{0}_{K+1}),\beta^{0}_{K+1} (1)=\gamma^{0,\frac{1}{3}}(e^{0}_{K+1})\in\gamma^{0,\frac{1}{3}}(\NS^1\backslash(t_0,t_1))$; ie $b^{0}_{K+1},e^{0}_{K+1}\in\NS^1\backslash(t_0,t_1)$, and
	\item {} $\beta^{0}_{K+1}$ is unlinked from $\gamma^{0,\frac{1}{3}}_{\Res(t_0,t_1)}$.
\end{itemize}
In other words, denote by $\D$ an embedded $2$--disk whose boundary is given by $\gamma^{0,\frac{1}{3}}_{\Res(t_0,t_1)}$. Then, the condition $\beta^{0}_{K+1} [0,1]\cap\D=\emptyset$ is sufficient to keep the unlinking property. We say that a $1$--parametric family of finger germs that satisfies this condition is \em valid\em. Hence, we must show that there always exists a \em valid \em $1$--parametric family of finger germs which extends $(\beta^{1}_{K+1},v^{1}_{K+1})$. 

Assume that the family $(\beta^{s}_{K+1},v^{s}_{K+1})$ does not satisfy this condition, ie $\beta^{0}_{K+1}[0,1]\cap\D=\{p_1,\ldots,p_k\}$. Let $\delta>0$ be small enough and let $\tilde{\D}\supsetneq\D$ be a small deformation of $\D$ such that $\tilde{\D}\cap\gamma^{\delta,\frac{1}{3}}(\NS^1)=\gamma^{\delta,\frac{1}{3}}[t_0,t_1]\subsetneq\partial\tilde{\D}$. Take any $\phi\in\Diff(\tilde{\D},\partial\tilde{\D})$, isotopic to the identity, such that $\phi(p_i)\in\tilde{\D}\backslash\D$, $i\in\{1,\ldots,k\}$, and extend it to a diffeomorphism $\Phi$ of $\R^3$ with compact support close enough to $\tilde{\D}$ isotopic to the identity. Let $\varepsilon>0$ be small enough and let $\Phi_s$, $s\in[0,1]$, be any isotopy between $\Phi_0=\Phi$ and $\Phi_1=\Id$, such that $\Phi_s=\Id$ for $s\geq\varepsilon$.

\begin{figure}[h]
	\includegraphics[scale=0.13]{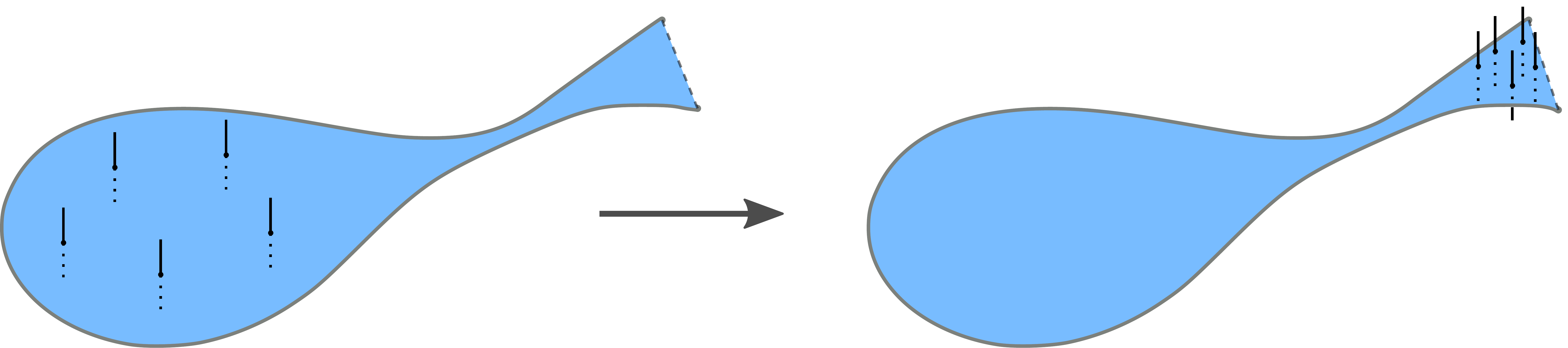}
	\centering
	\caption{Family of finger bones before ($\beta^{s}_{K+1}$) and after ($\tilde{\beta}^{s}_{K+1}=\Phi\circ\beta^{s}_{K+1}$) applying the described isotopy.\label{isotopiaBudney}}
\end{figure}

 The family of finger germs $(\tilde{\beta}^{s}_{K+1},\tilde{v}^{s}_{K+1})=(\Phi_s)_* (\beta^{s}_{K+1},v^{s}_{K+1})=(\Phi_s\circ\beta^{s}_{K+1}, d\Phi_s\circ v^{s}_{K+1})$ satisfies the condition $\tilde{\beta}^{0}_{K+1}\cap\D=\emptyset$ and extends $(\beta^{1}_{K+1},v^{1}_{K+1})$ so we find a \em valid \em family of finger germs which extends $(\beta^{1}_{K+1},v^{1}_{K+1})$, see Figure \ref{isotopiaBudney}.
 
 In conclusion, we may assume that the family of finger germs $(\beta^{s}_{K+1},v^{s}_{K+1})$ which extends $(\beta^{1}_{K+1},v^{1}_{K+1})$ is valid. Use this family to construct $\gamma^{s,u}$, $s\in[0,1]$, $u\in[0,\frac{1}{3}+\frac{1}{3(N-K)}]$, such that $\gamma^{s,\frac{1}{3}+\frac{1}{3(N-K)}}=\gamma^{s,\frac{1}{3}}\#\beta^{s}_{K+1}$.

Proceed in the same way with the rest of finger germs ($N-K$ in total), ie extending them to \em valid \em families of finger germs, to obtain a $2$--parametric family $\gamma^{s,u}$ with $s\in[0,1]$ and $u\in[0,\frac{2}{3}]$, such that 
\begin{itemize}
	\item[(i)] $\gamma^{s,0}=\gamma^s$,
	
	\item[(ii)] $\gamma^{s,\frac{2}{3}}=\gamma^{s}\#\beta^{s}_{1}\#\cdots\#\beta^{s}_{N}$,
	
	\item [(iii)] $\gamma^{i,u}$, $i\in\{0,1\}$, has one self--intersection at times $(t_0,t_1)\in\NS^1\times\NS^1$,
	
	\item[(iv)] There exists a finite sequence $0<u_1<\cdots<u_N<\frac{2}{3}$ such that
	$$\gamma^{s,u_j}\in\Imm(\NS^1,\R^3)\backslash\Emb(\NS^1,\R^3),$$
	for each $j\in\{1,\ldots,N\}$ and $s\in(0,1)$,
	
	\item[(v)] $\gamma^{s,u}\in\Emb(\NS^1,\R^3)$ for $u\neq u_1,\ldots,u_N$ and $s\in(0,1)$,
	
	\item[(vi)] The resolution $\gamma^{i,\frac{2}{3}}_{\Res(t_0,t_1)}$, $i\in\{0,1\}$, is a parametrized unknot and is unlinked from the other resolution $\gamma^{i,\frac{2}{3}}_{\Res(t_1,t_0)}$.

\end{itemize}

\textit{Step III. Create a $1$--parametric family of self--intersections using $\gamma^{s,\frac{2}{3}}$:} 

We will further extend the family $\gamma^{s,u}$, $u\in [0, \frac{2}3]$ to a new family $\gamma^{s,u}$, $u\in [0,1]$.

There exists a disk which bounds $\gamma^{i, \frac23}_{\Res(t_0,t_1)}$. Interpret this disk as a band $B_i :[0,1]\times[0,1]\rightarrow\R^3$, $i\in\{0,1\}$, satisfying:
\begin{itemize}
\item[(i)] $B_i (1,h)=\gamma^{i,\frac23} (t_0)=\gamma^{i, \frac23} (t_1)$,
\item[(ii)] $B_i (0,h)=\gamma^{i, \frac23} (\frac{t_0+t_1}{2})$,
\item[(iii)] $B_i (t,1)=\gamma^{i, \frac23} (\frac{t_0+t_1}{2}+t\frac{t_1-t_0}{2})$,
\item [(iv)] $B_i (t,0)=\gamma^{i, \frac23} (\frac{t_0+t_1}{2}-t\frac{t_1-t_0}{2})$,
\item[(v)] $B_i ((0,1)\times(0,1))\cap\gamma^{i,\frac23}(\NS^1)=\emptyset$.
\end{itemize}

Then, there exists a $1$--parametric family of bands $B_s:[0,1]\times[0,1]\rightarrow\R^3$ interpolating between $B_0$ and $B_1$ satisfying the properties (ii),(iii),(iv),(v) and 
\begin{itemize}
	\item [(i)'] $B_s (1,0)=\gamma^{s, \frac23} (t_0), B_s (1,1)= \gamma^{s, \frac23} (t_1)$ and $B_s (\{1\}\times(0,1))\cap\gamma^{s, \frac23} (\NS^1)=\emptyset$.
\end{itemize}

This is just built by continuously collapsing the two initial bands into small bands $\tilde{B}^i_{\delta}: [0,1] \times [0,1] \to \R^3$ defined as $\tilde{B}^i_{\delta}(t,h)=B^i(\delta \cdot t, h)$. From there, it is simple to interpolate between $\tilde{B}^0_{\delta}$ and $\tilde{B}^1_{\delta}$ for $\delta>0$ small enough.

In particular, $(\beta^{s}_{N+1}(-),v^{s}_{N+1}(-))=(B_s (1,-), \frac{\partial B_s}{\partial t}(1,-))$ defines a continuous family of finger germs for $\gamma^{s, \frac23} $. Use this family to create the required self--intersections. Ie define $\gamma^{s,u}$, $s\in[0,1]$, $u\in[0,1]$, such that 
\begin{itemize}
	\item[(i)] $\gamma^{s,0}=\gamma^s$
	\item [(ii)] $\gamma^{s,1}=\gamma^{s,\frac{2}{3}}\#\beta^{s}_{N+1}=\gamma^s\#\beta^{s}_{1}\#\cdots\#\beta^{s}_{N+1}$,
	\item [(iii)] $\gamma^{i,u}$, $i\in\{0,1\}$, has one self--intersection at times $(t_0,t_1)\in\NS^1\times\NS^1$,
	\item[(iv)] $\gamma^{s,1}$ has exactly self--intersection at times $(t_0,t_1)\in\NS^1\times\NS^1$,
	\item[(v)] There exists a finite sequence $0<u_1<\ldots<u_N <1$ such that $$\gamma^{s,u_j}\in\Imm(\NS^1,\R^3)\backslash\Emb(\NS^1,\R^3)$$ for each $j\in\{1,\ldots,N\}$ and $s\in(0,1)$,
	\item[(vi)] $\gamma^{s,u}\in\Emb(\NS^1,\R^3)$ for $u\neq 1, u_1,\ldots,u_N$ and $s\in(0,1)$.
	
\end{itemize}

This completes the argument to prove the first part of the Proposition. For the relative construction observe that since $\NS^1\backslash\{c\}$ is connected all the fingers germs $\beta^{s}_{j}$, $j\in\{1,\ldots,N+1\}$, can be chosen in such a way that $\beta^{s}_{j} (0),\beta^{s}_{j} (1)\neq \gamma^{s} (c)$, ie $b^{s}_{j},e^{s}_{j}\in\NS^1\backslash\{c\}$.

{\bf \em Case 2: $\{ s_1, \ldots, s_{\nu} \} \neq \emptyset$. }

Observe that in the proof of the first case all the deformations over the family of immersions $\gamma^{s,u}:\NS^1\rightarrow\R^3$ take place in a neighborhood of at most two points in the domain $z_0(s,u),z_1(s,u)\in\NS^1$. There is a finite number of steps at times $0<u_1<\cdots u_N\leq 1$, $N\geq 0$, such that $\gamma^{s,u}\in\Imm(\NS^1,\R^3)\backslash\Emb(\NS^1,\R^3)$, $s\in[0,1]$. We may assume in all the steps that the points $z_0(s_j,u_k),z_1(s_j,u_k),\tau_{s_j}^{0},\tau_{s_j}^{1}$, $(j,k)\in\{1,\ldots,\nu\}\times\{1,\ldots,N\}$, are pairwise different. Observe that a generic choice of the times $z_0(s,u),z_1(s,u)$ satisfies this assumption. Thus, this case follows immediately from the previous one.
\end{proof}

\subsubsection{Construction of the square.}
Theorem \ref{LegendrianAproximation} and Lemma \ref{TransverseFingers} allow us to state the last Proposition in the Legendrian setting, which provides the answer that we are looking for. To simplify the statement we introduce the following definition:

\begin{definition}
	Let $\gamma^s$, $s\in[0,1]$, be a path of Legendrian immersions such that
	\begin{itemize}
		\item[(i)] $\gamma^s$, $s\in[0,1]\setminus \{0, s_1,\ldots,s_\nu, 1\}$, for $\nu\geq 0$, where $0<s_1<\cdots<s_\nu<1$, is a Legendrian embedding,
		
		\item [(ii)] $\gamma^i\in\SLegImmG(\R^3)$, $i\in\{0,1\},$ has a self--intersection at times $(t_0,t_1)\in\NS^1\times\NS^1$; and $\gamma^{s_j}\in\SLegImmG(\R^3)$, $j\in\{1, \ldots, \nu \},$ has a self--intersection at times $(\tau_{s_j}^0,\tau_{s_j}^1)\in\NS^1\times\NS^1$.
	\end{itemize}
A homotopy $\gamma^{s,u}$, $(s,u)\in[0,1]\times[0,1]$, of paths of Legendrian immersions between $\gamma^{s,0}=\gamma^s$ and $\gamma^{s,1}$ is \em admissible \em if there exists a finite sequence $0<u_1<\ldots<u_N\leq1$ such that 
\begin{itemize}
	\item [(i)] $\gamma^{s,u_j}\in\SLegImmG(\R^3)$ for $j\in\{1,\ldots,N\}$ and $s\in(0,1)$,
	\item[(ii)] $\gamma^{s,u}\in\Leg(\R^3)$ for $u\neq u_1,\ldots,u_N$, and $s\in(0,1)\setminus \{s_1, \ldots, s_\nu \}$,	
	\item[(iii)] $\gamma^{i,u}\in\SLegImm(\R^3)$ has a self--intersection at times $(t_0,t_1)\in\NS^1\times\NS^1$ for $i\in\{0,1\}$ and $u\in[0,1]$,
	\item [(iv)] $\gamma^{s_j,u}\in\SLegImmG(\R^3)$, $j\in\{1, \ldots, \nu \},$ has a self--intersection at times $(\tau_{s_j}^0,\tau_{s_j}^1)\in\NS^1\times\NS^1$ and $u\in [0,1]$. Moreover, the intersection times of $\gamma^{s_j,u_k}$, for $k\in\{1,\ldots,N\}$, do not coincide, ie the real numbers $t_0, t_1, \tau_{s_j}^0, \tau_{s_j}^1$ are pairwise different.
	\item[(v)] Moreover, the segments $\gamma^{i,u}\in\SLegImm(\R^3)$ has self--intersections of type $\Sigma^{1,0}$ except at an even finite number of values $u\in (0,1)$ in which the self--intersections will be of type $\Sigma^{1,1,0}$.
\end{itemize}
\end{definition}  

Let us detail some co--orientation issues in $\Sigma^{1,0}$. There is a standard orientation of the normal bundle of $\Sigma^{1,0}\backslash\Sigma^{1,1,0}$. It is defined at $\gamma \in \Sigma^{1,0}\backslash\Sigma^{1,1,0}$ as follows. Denote $\gamma(t_0)= \gamma(t_1)$ the self--intersection point. Let us order the two branches. We declare that the branch through $\gamma(t_0)$ goes first and the branch through $\gamma(t_1)$ second if and only if $\langle \gamma'(t_0),\gamma'(t_1)\rangle $ is a positive basis of $\xi_{\gamma(t_0)}$. We want to define the normal vector to $\Sigma^{1,0}\backslash\Sigma^{1,1,0}$ at $\gamma$ that is a field in $T_{\gamma} \LegImm(\R^3)$, ie a vector field on $\gamma^* T\R^3$. We define it as a vector field supported around $\gamma(t_1)$ such that it extends to a contact vector field in a small ball $\tilde{B}$ centered at $\gamma(t_1)\in \R^3$. This contact vector field has as associated Hamiltonian function a cut-off function $\chi:\tilde{B}\rightarrow[0,1]$ satisfying that $\chi(0)=1$, $d\chi(0)=0$  and $\chi|_{\Op(\partial\tilde{B})}=0$. In other words, we are pushing through the Reeb vector field the second branch.

This definition does not extend continuously to $\Sigma^{1,1,0}$. However, it provides a co--orientation for $\Sigma^{1,0}\backslash\Sigma^{1,1,0}$. We call this co--orientation the {\em standard} one. Let $\gamma^s$ the path of Legendrian immersions such that $\gamma^i\in \Sigma^{1,0}$, $i\in\{0,1\}$. The vectors $v_i= (-1)^{i-1} \frac{d\gamma^s}{ds}|_{s=i} \in T_{\gamma^i}\LegImm(\R^3)$, $i\in\{0,1\}$,  determine two orientations of the normal bundle of $\Sigma^{1,0}\backslash\Sigma^{1,1,0}$ at $\gamma^0$ and $\gamma^1$. We need the co--orientation defined by these two vectors to coincide with the standard one.

In these terms, the main result of the section is the following:

\begin{proposition}\label{CreatingLegendrianTangency}
Let $\gamma^s$, $s\in[0,1]$, be a path of Legendrian immersions such that
\begin{itemize}
	\item[(i)] $\gamma^s$, $s\in(0,1)\setminus \{s_1, \ldots, s_\nu \}$, is a Legendrian embedding,
	\item [(ii)] $\gamma^i\in\SLegImmG(\R^3)$, $i\in\{0,1\}$, has one self--intersection at times $(t_0,t_1)\in\NS^1\times\NS^1$ and $\gamma^s\in\SLegImmG(\R^3)$, $i\in\{s_1, \ldots, s_\nu \}$, has a self--intersection at times $(\tau_{s_j}^0,\tau_{s_j}^1)\in\NS^1\times\NS^1$,
	\item[(iii)] moreover, $v_i= (-1)^{i-1} \frac{d\gamma_t}{dt}|_{t=i}$ defines the standard co--orientation at $\gamma^i$, $i\in\{0,1\}$.
\end{itemize}
Then, there exists an admissible homotopy $\gamma^{s,u}\in\LegImm(\R^3)$, $(s,u)\in[0,1]\times[0,1]$, with $\gamma^{s,0}=\gamma^s$ such that $\gamma^{s,1}\in\SLegImmG(\R^3)$ has one self--intersection at times $(t_0,t_1)\in\NS^1\times\NS^1$.
\end{proposition}

We need a sequence of simple lemmas. We assume them to be true and we will prove them later.

The following lemma is fairly general in homotopy theory saying that if you have a retraction $K\rightarrow L$ then, by precomposing with the retraction, any map on $L$ extends to $K$. For our purpose $K=[0,1]\times[0,1]$ and $L\subseteq K$, such that $K$ deformation retracts to $L$, and we have a map $L\rightarrow \LegImm(\R^3)$. This allows us to extend local deformations in a path of Legendrian immersions to produce a complete square.

\begin{lemma}\label{LemaCuadrado}
	Let $0<\delta<s_1$ be a positive number. Fix $\gamma^{s,u}\in\LegImm(\R^3)$, $(s,u)\in L=[0,\delta]\times[0,1]\cup[0,1]\times[0,\delta]$. Assume that 
	\begin{itemize}
		\item [(i)] $\gamma^{i,u}\in\SLegImm(\R^3)$, $i\in\{0,1\}$,
		\item[(ii)] $\gamma^{s,u}\in\Leg(\R^3)$ for $s\in(0,1)\backslash\{s_1,\ldots,s_\nu\}$.
	\end{itemize}
Then, there exists a family $\tilde{\gamma}^{s,u}\in\LegImm(\R^3)$, $(s,u)\in[0,1]\times[0,1]$, satisfying
\begin{itemize}
	\item [(i)] $\tilde{\gamma}^{s,u}=\gamma^{s,u}$ if $(s,u)\in L$,
	\item[(ii)] $\tilde{\gamma}^{s,u}\in\Leg(\R^3)$ if $s\in(0,1)\backslash\{s_1,\ldots,s_\nu\}$,
	\item[(iii)] $\tilde{\gamma}^{1,u}=\gamma^{1,\delta}$ and $\tilde{\gamma}^{s_j,u}=\gamma^{s_j,\delta}$, $j\in\{1,\ldots,\nu\}$, if $u\in[\delta,1]$.
\end{itemize}
\end{lemma}

The following Lemma allows us to assume, up to reparametrization, that there is a fixed cusp in any family of Legendrian immersions parametrized by the segment (we need to add parametric Reidemeister type 1 Legendrian moves, see the proof for the details). 

\begin{lemma}\label{LemaCuspides}
	Let $\gamma^s$ be a path of Legendrian immersions such that 
\begin{itemize}
	\item[(i)] $\gamma^s$, $s\in(0,1)\setminus \{s_1, \ldots, s_\nu \}$, is a Legendrian embedding,
	\item [(ii)] $\gamma^i\in\SLegImmG(\R^3)$, $i\in\{0,1\}$, has one self--intersection at times $(t_0,t_1)\in\NS^1\times\NS^1$; and $\gamma^s\in\SLegImmG(\R^3)$, $s\in\{s_1, \ldots, s_\nu \}$, has a self--intersection at times $(\tau_{s_j}^0,\tau_{s_j}^1)\in\NS^1\times\NS^1$.	
\end{itemize}

There exists an admissible homotopy $\gamma^{s,u}$, $(s,u)\in[0,1]\times[0,1]$, with $\gamma^{s,0}=\gamma$ satisfying
\begin{itemize}
	\item[(i)] $\gamma^{s,u}\in\Leg(\R^3)$ for $s\in(0,1)\backslash\{s_1,\ldots,s_\nu\}$,
	\item[(ii)] there exists a sequence of times $c_1,\ldots,c_{2k}\in\NS^1$ such that all the cusps of $\gamma^{s,1}$ are at times $c_i$.
\end{itemize}
\end{lemma}

The following lemma plays the role of Steps I and II in the proof of Proposition \ref{CreatingSmoothTangency} and solves our problem modulus the creation of the interpolating bands. 

\begin{lemma}\label{LemaUnknottingLegendriano}
	Let $\gamma^s$ be a path of Legendrian immersions such that 
	\begin{itemize}
		\item[(i)] $\gamma^s$, $s\in(0,1)\setminus \{s_1, \ldots, s_\nu \}$, is a Legendrian embedding,
	\item [(ii)] $\gamma^i\in\SLegImmG(\R^3)$, $i\in\{0,1\}$, has one self--intersection at times $(t_0,t_1)\in\NS^1\times\NS^1$; and $\gamma^s\in\SLegImmG(\R^3)$, $s\in\{s_1, \ldots, s_\nu \}$, has a self--intersection at times $(\tau_{s_j}^0,\tau_{s_j}^1)\in\NS^1\times\NS^1$.	
	\end{itemize}

	There exists an admissible homotopy $\gamma^{s,u}$, $(s,u)\in[0,1]\times[0,1]$, with $\gamma^{s,0}=\gamma^s$ such that
	\begin{itemize}
		\item[(i)] The self--intersection times $(t_0(i,u),t_1(i,u))$ of $\gamma^{i,u}$, $i\in\{0,1\}$, satisfy that $t_1(i,1)<\tau_{s_j}^0,\tau_{s_j}^1<t_0 (i,0)$, for $j\in\{1,\ldots,\nu\}$.
		\item[(ii)] $\gamma^{s,1}\in\Leg(\R^3)$ for $s\in(0,1)\setminus \{s_1, \ldots, s_\nu \}$,
	    \item[(iii)] $\gamma^{i,1}_{\Res(t_0,t_1)}$, $i\in\{0,1\}$, is a Legendrian (stabilized) unknot and is unlinked from $\gamma^{i,1}_{\Res(t_1,t_0)}$.
    \end{itemize}	
\end{lemma}

We will need to extend a strict Legendrian immersion into a path in which the self--intersection points cross a positive (negative) stabilization. This will be key to change the formal invariants associated to the resolutions $\gamma_{\Res(t_0, t_1)}$ and $\gamma_{\Res(t_1, t_0)}$: since we will add/substract a stabilization to each of the resolutions. The price that you pay is that the new path has two cusp points ($\Sigma^{1,1,0}$). We are being precise here: \em moving the self--intersection point through two cusps in the front projection creates two cusps in the induced curve in the moduli of strict Legendrian immersions.\em 

\begin{lemma} \label{lem:cuspcreation}
Assume that there is a $2$--disk $\varphi:\D\rightarrow\LegImm(\R^3)$ such that it posseses  just one connected curve $C_0$, parametrized as $\gamma^s$, $s\in[0,1]$, of strict immersions that is a properly embedded segment on the disk. Denote by $t_0$ and $t_1$ the self--intersection times.
Moreover, assume that $\gamma^s \in C_0 \subset \SLegImmG(\R^3)$ possesses a positive DS at time $t_0+ \delta$ for $\delta>0$ small enough. Then:
\begin{enumerate}
\item There is a local $C^0$ deformation of the disk  $\varphi_t:\D\rightarrow\LegImm(\R^3)$, $t\in[0,1]$, relative to the boundary, such that for $t=1$, we have that the curve $C_0$ is deformed into a curve $C_1$ with $4$ cusps: points of type $\Sigma^{1,1,0}$ (a DS in the moduli). Moreover, $d_{C^0}(\varphi,\varphi^t)=O(\delta)$.
\item Parametrize the curve  $C_1$ by $\tilde{\gamma}^s$, $s\in [0,1]$. Assume that $\tilde{\gamma}^{1/2}$ lies between the second and the third cusp in the moduli, and denote its self-intersection times $t_0$ and $t_1$, then $(\tilde{\gamma}^{1/2})_{\Res(t_0, t_1)}^1$ inherits a negative stabilization and $(\tilde{\gamma}^{1/2})_{\Res(t_1, t_0)}^{1}$ inherits a positive stabilization with respect to the resolutions associated to $\gamma^{1/2}$.
\end{enumerate}
\end{lemma}

The {\em standard} normal vector for a strict Legendrian immersion of type $\Sigma^{1,0}$ is preserved when we cross over a cusp (a point of type in $\Sigma^{1,1,0}$) in the moduli space (see Figure \ref{DeformationModuli}), so we can easily assume that the orientations are preserved in the previous construction. Ie the {\em standard} normal vector is the same as the vector field normal to the curve in the $1$--parametric family constructed by the rule of avoiding tangencies with the curve of strict immersions.

\begin{figure}[h]
	\centering
	\subfloat[Change in the topology of the disk in the moduli space via $\varphi_t$.]{
		\label{fig:moduliDS}
		\includegraphics[width=0.5\textwidth]{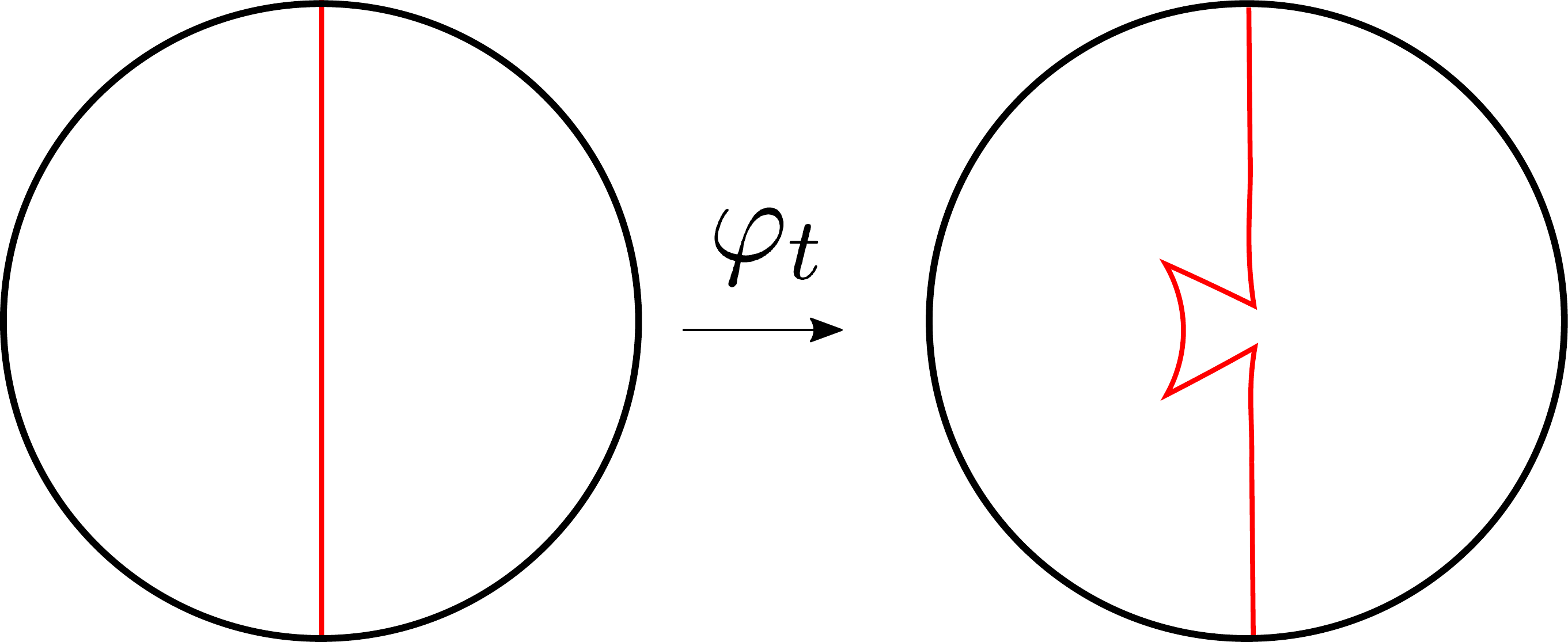}}\phantom{ aa }
	\subfloat[Birth of four cusps in the moduli space.]{
		\label{fig:3dimensionalDS}
		\includegraphics[width=0.3\textwidth]{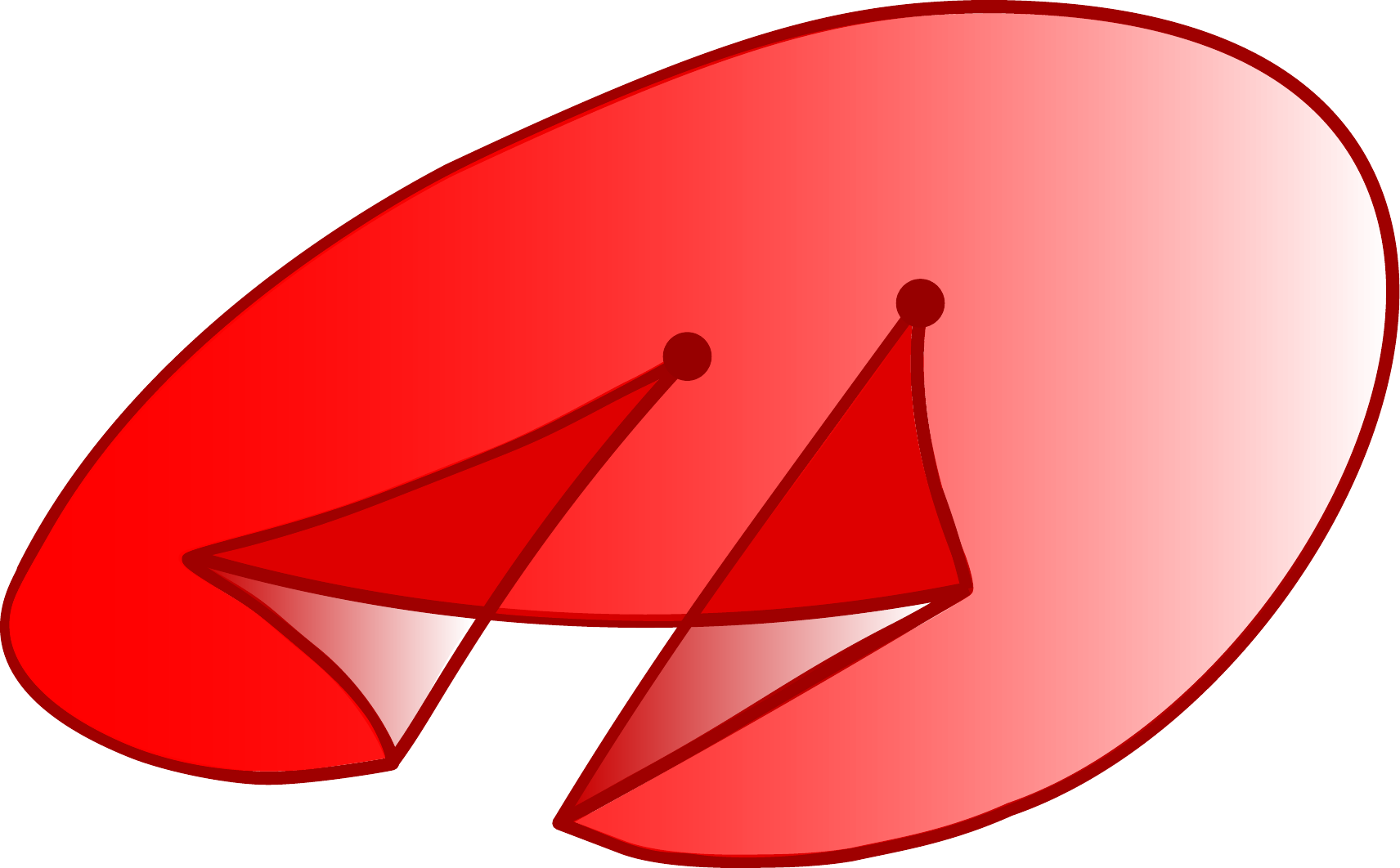}}
	\caption{The $1$--parametric family of disks $\varphi_t:\D\rightarrow\LegImm(\R^3)$, $t\in[0,1]$, described in Lemma \ref{lem:cuspcreation}.}
\end{figure}

The next lemma provides a way to match the formal invariants of the resolutions of two generic strict Legendrian immersions. This is a necessary condition to join two generic strict Legendrian immersions by a path of generic strict Legendrian immersions. In other words we want the ends of the path to be stabilized Legendrian unknots, with the same formal invariants. Therefore, by the classification theorem of Legendrian unknots \cite{EliashbergFraser}, they are the same Legendrian knot. Thanks to that fact, it makes sense to mimic the smooth proof and look for a $1$--parametric family of fingers interpolating between both of them.

\begin{lemma}\label{LemaInvariantesFormales}
		Let $\gamma^s$ be a path of Legendrian immersions such that 
	\begin{itemize}
		\item[(i)] $\gamma^s$, $s\in(0,1)\setminus \{s_1, \ldots, s_\nu \}$, is a Legendrian embedding,
	\item [(ii)] $\gamma^i\in\SLegImmG(\R^3)$, $i\in\{0,1\}$, has one self--intersection at times $(t_0,t_1)\in\NS^1\times\NS^1$; and $\gamma^s\in\SLegImmG(\R^3)$, $s\in\{s_1, \ldots, s_\nu \}$, has a self--intersection at times $(\tau_{s_j}^0,\tau_{s_j}^1)\in\NS^1\times\NS^1$.	
	\end{itemize}
	There exists an admissible homotopy $\gamma^{s,u}$, $(s,u)\in[0,1]\times[0,1]$, with $\gamma^{s,0}=\gamma^s$ satisfying
	\begin{itemize}
			\item[(i)] $\gamma^{s,1}\in\Leg(\R^3)$, for $s\in(0,1)\setminus \{s_1, \ldots, s_\nu \}$,
		    \item [(ii)] $\tb(\gamma^{0,1}_{\Res(t_0,t_1)})=\tb(\gamma^{1,1}_{\Res(t_0,t_1)})$ and $\Rot(\gamma^{0,1}_{\Res(t_0,t_1)})=\Rot(\gamma^{1,1}_{\Res(t_0,t_1)})$ .
	\end{itemize}
\end{lemma}

The next Lemma provides the last step in the proof, we show how to mimic the band construction in the Legendrian setting. It allows us to connect two strict Legendrian immersions by a path of strict Legendrian immersions. Thus, to prove Proposition \ref{CreatingLegendrianTangency} we just need to deform our initial path, via an admissible homotopy, to a path which satisfies the hypothesis of this lemma.

\begin{lemma}\label{LemaBandaLegendrianoMismosInvariantes}
	Let $\gamma^s$, $s\in[0,1]$, be a path of Legendrian immersions such that 
	\begin{itemize}
		\item[(i)] $\gamma^s$, $s\in(0,1)\setminus \{s_1, \ldots, s_\nu \}$, is a Legendrian embedding,
	\item [(ii)] $\gamma^i\in\SLegImmG(\R^3)$, $i\in\{0,1\}$, has one self--intersection at times $(t_0,t_1)\in\NS^1\times\NS^1$; and $\gamma^s\in\SLegImmG(\R^3)$, $s\in\{s_1, \ldots, s_\nu \}$ has a self--intersection at times $(\tau_{s_j}^0,\tau_{s_j}^1)\in\NS^1\times\NS^1$. Moreover,  $t_1<\tau_{s_j}^0,\tau_{s_j}^1<t_0 $, for $j\in\{1,\ldots,\nu\}$.
		\item[(iii)] There exists a time $c_1\in(t_0,t_1)\subseteq\NS^1$ such that $\gamma^s$ has a cusp at time $c_1$ and $\gamma^{s}_{|\Op(\{c_1\})}$ is constant,
		\item[(iv)] $\gamma^{i}_{\Res(t_0,t_1)}$, $i\in\{0,1\}$, is a Legendrian (stabilized) unknot and is unlinked from $\gamma^{i}_{\Res(t_1,t_0)}$,
		\item[(v)] $\gamma^{0}_{\Res(t_0,t_1)}$ and $\gamma^{1}_{\Res(t_0,t_1)}$ are formally isotopic.
	\end{itemize}

Then, there exists an admissible homotopy $\gamma^{s,u}$, $(s,u)\in[0,1]\times[0,1]$, with $\gamma^{s,0}=\gamma^s$ satisfying
\begin{itemize}
	\item[(i)] $\gamma^{s,u}_{|\Op(\{c_1\})}$ is constant,
	\item[(ii)] $\gamma^{s,u}\in\SLegImm(\R^3)$, $(s,u)\in\{0,1\}\times[0,1]\cup[0,1]\times\{1\}$, has one self--intersection at times $(t_0,t_1)\in\NS^1\times\NS^1$.
\end{itemize}
\end{lemma}

\begin{proof}[Proof of Proposition \ref{CreatingLegendrianTangency}]
	The proof is divided in four steps.
	
	\textit{Step I. Fix the cusp times:} Apply Lemma \ref{LemaCuspides} to produce a family $\gamma^{s,u}$, $(s,u)\in[0,1]\times[0,\frac{1}{4}]$, such that 
		\begin{itemize}
		\item [(i)] $\gamma^{s,0}=\gamma^s$,
			\item[(ii)] $\gamma^{s,u}\in\Leg(\R^3)$ for $s\in(0,1)\backslash\{s_1,\ldots,s_\nu\}$,
		\item[(iii)] $\gamma^{i,u}\in\SLegImmG(\R^3)$ has a self--intersection at times $(t_0,t_1)\in\NS^1\times\NS^1$ for $i\in\{0,1\}$; and $\gamma^{s,u}\in\SLegImmG(\R^3)$, $s\in\{s_1, \ldots, s_\nu \}$, has a self--intersection at times $(\tau_{s_j}^0,\tau_{s_j}^1)\in\NS^1\times\NS^1$,
		\item[(iv)] there exists a sequence of times $c_1,\ldots,c_{2k}\in\NS^1$ such that, the cusps of $\gamma^{s,\frac{1}{4}}$ are at times $c_i$.
	\end{itemize}
	
	Assume that $c_1\in(t_0,t_1)$ and that $\gamma^{s,\frac{1}{4}}_{|\Op(\{c_1\})}$ is constant. This property will be preserved all over the construction. 
	
	\textit{Step II. Unknotting process:} 
	Apply Lemma \ref{LemaUnknottingLegendriano} to find an extended family $\gamma^{s,u}$, $(s,u)\in[0,1]\times[0,\frac{2}{4}]$, satisfying 	
	\begin{itemize}
		\item [(i)] $\gamma^{s,u}$ is an admissible homotopy, 
		\item[(ii)] $\gamma^{i,u}$, $i\in\{0,1\}$, has one self--intersection at times $(t_0(i,u),t_1(i,u))$ and $\gamma^{s,u}$, $s\in\{s_1, \ldots, s_\nu \}$, has a self--intersection at times $(\tau_{s_j}^0,\tau_{s_j}^1)$. Moreover,  $t_1(i,1)<\tau_{s_j}^0,\tau_{s_j}^1<t_0 (i,0)$, for $(i,j)\in\{0,1\}\times\{1,\ldots,\nu\}$.
		\item[(iii)] $\gamma^{s,u}_{|\Op(\{c_1\})}$ is constant for $u\geq\frac14$,
		\item[(iv)] The resolution $\gamma^{i,\frac24}_{\Res(t_0,t_1)}$, $i\in\{0,1\}$, is an (stabilized Legendrian) unknot unlinked from the other one.
	\end{itemize}
	
	\textit{Step III. Match the formal invariants of the resolutions $\gamma^{0,\frac24}_{\Res(t_0,t_1)}$ and $\gamma^{1,\frac24}_{\Res(t_0,t_1)}$:} Use Lemma \ref{LemaInvariantesFormales} to extend our family to $\gamma^{s,u}$, $(s,u)\in[0,1]\times[0,\frac34]$, such that 
		\begin{itemize}
		\item [(i)] $\gamma^{s,u}$ is an admissible homotopy, 
		\item[(ii)] $\gamma^{i,u}$, $i\in\{0,1\}$, has one self--intersection at times $(t_0(i,u),t_1(i,u))$ and $\gamma^{s,u}$, $s\in\{s_1, \ldots, s_\nu \}$, has a self--intersection at times $(\tau_{s_j}^0,\tau_{s_j}^1)$. Moreover,  $t_1(i,1)<\tau_{s_j}^0,\tau_{s_j}^1<t_0 (i,0)$, for $(i,j)\in\{0,1\}\times\{1,\ldots,\nu\}$.
		\item[(iii)] $\gamma^{s,u}_{|\Op(\{c_1\})}$ is constant for $u\geq\frac{1}4$,
		\item[(iv)] The resolution $\gamma^{i,\frac34}_{\Res(t_0,t_1)}$, $i\in\{0,1\}$, is an unknot unlinked from the other one,
		\item[(v)] $\gamma^{0,\frac34}_{\Res(t_0,t_1)}$ and $\gamma^{1,\frac34}_{\Res(t_0,t_1)}$ are formally isotopic Legendrian embeddings (stabilized unknots).
	\end{itemize}
	
    \textit{Step IV. Create a $1$--parametric family of self--intersections:} Finally, observe that $\gamma^{s,\frac34}$ satisfies the hypothesis of Lemma \ref{LemaBandaLegendrianoMismosInvariantes}. So, apply this Lemma to extend our family to $\gamma^{s,u}$, $(s,u)\in[0,1]\times[0,1]$. The family $\gamma^{s,u}$ satisfies the required properties.
\end{proof}

\subsubsection{Proofs of the Lemmas \ref{LemaCuadrado}, \ref{LemaCuspides}, \ref{LemaUnknottingLegendriano}, \ref{lem:cuspcreation}, \ref{LemaInvariantesFormales} and \ref{LemaBandaLegendrianoMismosInvariantes}} \textbf{ \\Case 1: $\{ s_1, \ldots, s_{\nu} \}= \emptyset$.}

We assume that the set of points $\{ s_1, \ldots, s_{\nu} \}$ is empty. We explain afterwards the reason why the same proofs  work for the general case.

\begin{proof}[Proof of Lemma \ref{LemaCuadrado}]
	Take a retraction $r:[0,1]\times[0,1]\rightarrow L$ such that: \begin{itemize}
		\item[(i)] $r(1,u)=(1,\delta)$ for all $u\in[\delta,1]$.
		\item[(ii)] $r\left((0,1)\times [0,1]\right)\subseteq(0,1)\times [0,1].$
	\end{itemize}
	Then define $\tilde{\gamma}^{s,u}=\gamma^{r(s,u)}$.
\end{proof}

\begin{proof} [Proof of Lemma \ref{LemaCuspides}]
	Lemma \ref{EliminationRI} and Corollary \ref{CuspSameTime} also work for $\gamma^s$.
\end{proof}

\begin{proof} [Proof of Lemma \ref{LemaUnknottingLegendriano}]
	The first two steps in the proof of Proposition \ref{CreatingSmoothTangency} readily imply, via Lemma \ref{TransverseFingers}, the result.
\end{proof}

\begin{proof} [Proof of Lemma \ref{lem:cuspcreation}]
Draw in the Lagrangian projection a picture with two intersecting branches: one of the branches with slope $-1$  and orientation from right to left and the other one horizontal with a stabilization on it. The stabilization lies on the right of the intersection point. We continuously move the intersection point from left to right, making sure that the area of the intersection point is zero. So we get an actual intersection in the Legendrian above. Visual inspection shows that if the track intersection point has area zero all the other intersections in the Lagrangian projection have non zero area, therefore they are not actual intersections in the Legendrian. There are $2$ times $t_0 < \tau_0 < \tau_1 <t_1$ in which the branch becomes tangent. Those times correspond to cusps in the moduli space: curves of type $\Sigma^{1,1,0}$. See Figure \ref{DobleCusp} for a description of half of the deformation (crossing of the first stabilization). Figure \ref{fig:3dimensionalDS} for a description of the $2$-dimensional family as a subset of a $3$-dimensional ball of Legendrian immersions.
\begin{figure}[h] 
	\includegraphics[scale=0.55]{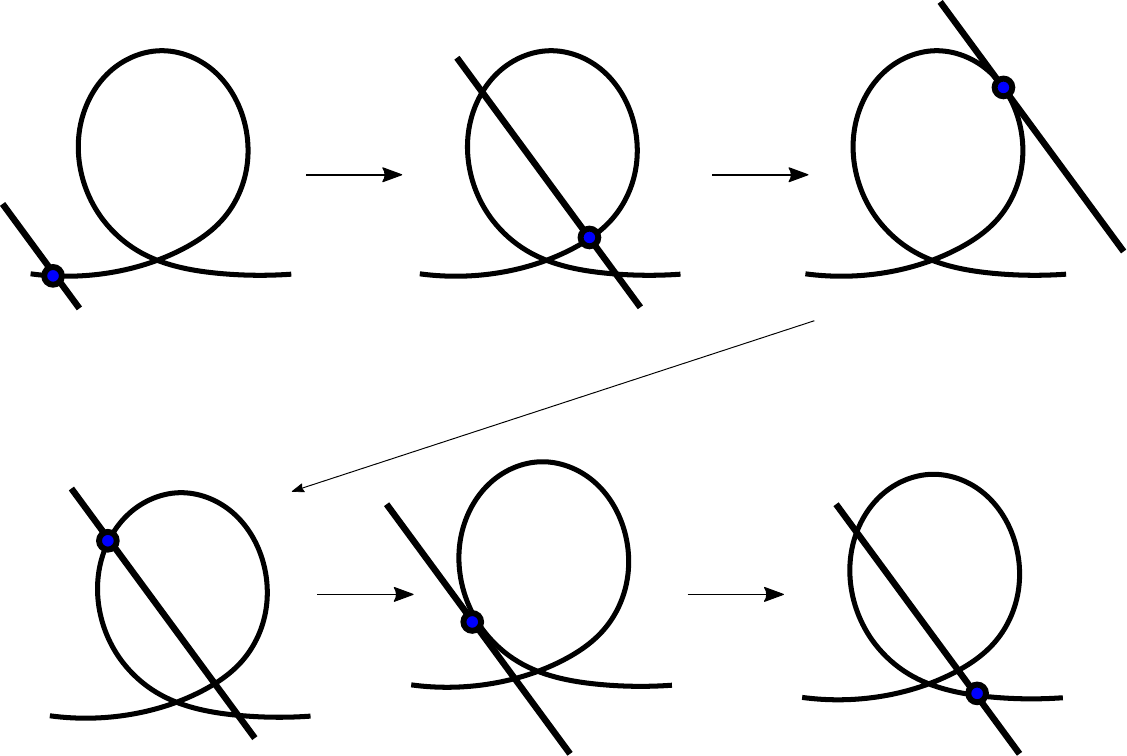}
	\centering
	\caption{Description of the movie in Lagrangian projection.}\label{DobleCusp}
\end{figure}

The last statement follows if we are able to extend the standard orientation over the cusps. We want to create a vector field  $v_s \in T_{\gamma_s} \LegImm(\R^3)$ transverse to the curve that defines the orientation in the moduli space given by the Figure \ref{DeformationModuli}. Flowing along this vector, we obtain the required band. First, observe that away from the $2$ cusp points this just corresponds to the standard orientation. Denote by $\gamma_{s}|_{s\in [0,1/2]}$ the part of the curve that crosses the first two cusps in the moduli and fix $\gamma_{1/6}$ and $\gamma_{2/6}$ to be those cusps. For $s\in [0, \varepsilon]$, with $\varepsilon>0$ small enough, we interpolate between the standard normal vector for $s=0$ and for $s=\varepsilon$  the vector defined by a Hamiltonian function constant equal to $1$ in the image of a segment joining two points $t_{0}^{'},t_{1}^{'}\in \NS^1$ and supported on an arbitrarily small neighborhood of it. We require $t_0< t_{0}^{'} < t_1 < t_{1}^{'}  < t_0$ . This last condition allows to homotope those two vectors through transverse vectors. The vector $v_{\varepsilon}$ is clearly equivalent the same as the one obtained by adding a positive bump (area) to the lagrangian projection at $t_{0}^{'}$ and a negative bump at $t_{1}^{'}$.  Thus, the deformation has support in $[t'_{0}, t'_{1}]$. The algorithm to extend the deformation to the cusp points of the moduli space is depicted in Figure \ref{PeliculaBumps} and it corresponds to change the bump function from $t_{i}^{'}$ to $t_{i}^{'}+2 \pi$, $i\in\{0,1\}$, creating a family $t_{i,s}^{'}$, $s\in [0,1/2]$, such that $t_{i,s}^{'}$ is strictly increasing with $s$and such that  $t_{0,1/6}^{'}=t_0$and $t_{0,2/6}^{'}=t_1$  Likewise, $t_{1,1/6}^{'}=t_0$ and $t_{1,2/6}^{'}=t_1$. See Figure \ref{PeliculaBumps}. This clearly produces a continuous family of the deformations that removes the self--intersections. See Figure \ref{DeformationModuli} to visualize the deformation as a tangent vector in the moduli space. 

\begin{figure}[h] 
	\includegraphics[scale=0.55]{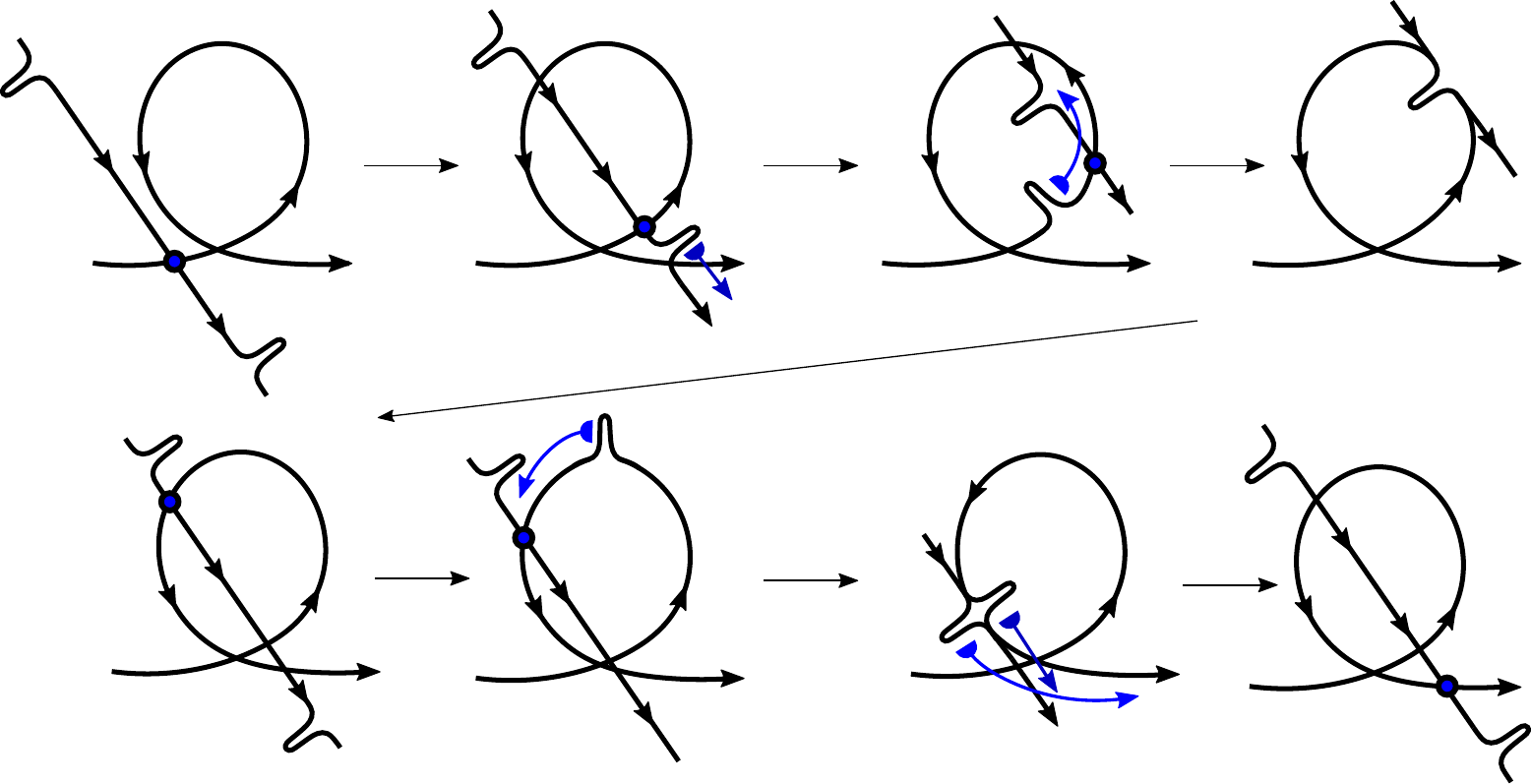}
	\centering
	\caption{Algorithm in Lagrangian projection to extend the deformation to the cusp points of the moduli space.}\label{PeliculaBumps}
\end{figure}

\begin{figure}[h] 
	\includegraphics[scale=1]{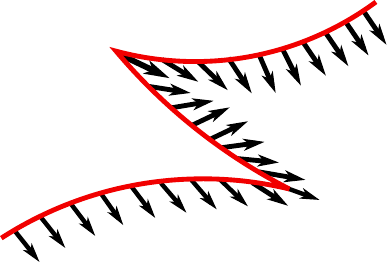}
	\centering
	\caption{Deformation in the moduli space.}\label{DeformationModuli}
\end{figure}

\end{proof}

\begin{proof}[Proof of Lemma \ref{LemaInvariantesFormales}]
	Define  $\gamma^{s,u}$, $(s,u)\in[0,1]\times[0,1]$, as follows:
	\begin{itemize}
		\item [(i)] If the parity of $\tb(\gamma^{0}_{\Res(t_0,t_1)})$ and $\tb(\gamma^{1}_{\Res(t_0,t_1)})$ coincide define $\gamma^{s,u}=\gamma^s$ for $u\in[0,\frac{1}{3}]$. In other case add one DS in $\gamma^{s}$ such that for $s=0,1$ it is inserted on the segment $(t_0,t_1)$ for $u\in[0,\frac{1}{6}]$ . Then, apply Lemma \ref{lem:cuspcreation} to create a band $\hat \gamma^{s',v}$ $s'\in [0,1]$, $v\in [-\varepsilon, \varepsilon]$, such that $\hat\gamma^{0,v}= \gamma^{v, \frac16}$ for $v\in[0, \varepsilon]$. Thus we define $\gamma^{s,u}= \hat\gamma^{6(u- \frac16), s}$, $s\in[0,\varepsilon]$, $u \in [\frac16 ,\frac13]$. Finally, apply Lemma \ref{LemaCuadrado} to extend $\gamma^{s,u}$ for all $s\in [0,1]$, $u\in [\frac16, \frac13]$. Thus, we have made sure that $\tb(\gamma^{0,\frac13}_{\Res(t_0,t_1)})-\tb(\gamma^{1, \frac13}_{\Res(t_0,t_1)})=2k$, $k\in \Z$. Recall from \cite[Section 2.7]{Etnyre} that the Thuston--Bennequin invariant decreases by 1 after a stabilization and that the rotation number increases (decreases) by 1 whenever the stabilization is positive (negative).

		\item[(ii)] If $k>0$, for each $i\in\{1,\ldots,k\}$, take a path $B_i:[0,1]\rightarrow\NS^1,s\mapsto B_i (s)$, such that $B_i (0)\in(t_0,t_1)$ and $B_i (1)\in\NS^1\backslash[t_0,t_1]$. Define recursively $\gamma^{s,u}$ in $[0,1]\times[\frac{1}{3},\frac{1}{3}+\frac{i}{3k}]$ creating a DS on $\gamma^{s,\frac{1}{3}+\frac{i-1}{3k}}$ at time $B_i (s)$.
		
		Observe that a DS decreases the Thurston--Bennequin invariant by $2$. Thus, we obtain 
\begin{equation} \label{f:cuadre} \tb(\gamma^{0,\frac{2}{3}}_{\Res(t_0,t_1)})=\tb(\gamma^{1,\frac{2}{3}}_{\Res(t_0,t_1)}).
\end{equation}
		
		Do the symmetric sequence of moves in case $k<0$ to obtain again (\ref{f:cuadre}).
		
		\item[(iii)] For any $\gamma\in\Leg(\R^3)$ we have $\tb(\gamma)+\Rot(\gamma)\in2\Z+1$ (see \cite[Proposition 3.5.23]{GeigesCont}). Thus, since $\tb(\gamma^{0,\frac{2}{3}}_{\Res(t_0,t_1)})=\tb(\gamma^{1,\frac{2}{3}}_{\Res(t_0,t_1)})$ we conclude that 
\begin{equation} \label{eq:rotcuadre}	
\Rot(\gamma^{0,\frac{2}{3}}_{\Res(t_0,t_1)})-\Rot(\gamma^{1,\frac{2}{3}}_{\Res(t_0,t_1)})=2m.
\end{equation}
		Assume that $m\geq0$, the other case is totally analogous. Add a DS in $(t_0,t_1)$ and consecutively apply Lemma \ref{lem:cuspcreation} and Lemma \ref{LemaCuadrado}.This adds a square $\gamma^{s, u}$, $s\in [0,1]$, $u \in[\frac23, \frac23 +\frac1{6m}]$. Now add a DS in $(t_1, t_0)$ and consecutively apply Lemma \ref{lem:cuspcreation} and Lemma \ref{LemaCuadrado}. This adds a second square. In the new top $\gamma^{s, \frac23 +\frac1{3m}}$, the difference (\ref{eq:rotcuadre}) has decreased by $2$. Repeat $m-1$ more times to make $m=0$ in that equation.

	It follows that $\gamma^{0,1}_{\Res(t_0,t_1)}$ and $\gamma^{1,1}_{\Res(t_0,t_1)}$ have the same formal invariants.
	\end{itemize}
\end{proof}

\begin{proof} [Proof of Lemma \ref{LemaBandaLegendrianoMismosInvariantes}]
	
The proof is divided in several steps.

\textit{Step 0: Preparation.}

By hypothesis, the standard co--orientations are matched, so we may assume that the self--intersection times $t_0$ and $t_1$ satisfy that $\{(\gamma^i)'(t_0), (\gamma^i)'(t_1)\}$, $i\in\{0,1\}$; is a positive basis of $\xi$ at the point of self--intersection $\gamma^i(t_0)=\gamma^i(t_1)$, and that there exists $\varepsilon>0$ such that 
\begin{itemize}
	\item $\gamma^{\varepsilon}$ is obtained from $\gamma^{0}$ by pushing $\gamma^{0}(\Op(\{t_1\}))$ in the Reeb direction for time $\varepsilon$,
	\item $\gamma^{1-\varepsilon}$ is obtained from $\gamma^{1}$ by pushing $\gamma^{1}(\Op(\{t_1\}))$ in the Reeb direction for time $\varepsilon$.
\end{itemize}

We may also think that $\gamma^i=\hat{\gamma}^{i}\#(\beta^{i}_{0},v_0)$ where $\hat{\gamma}^i$ is a Legendrian embedding and $\beta^{i}_{0}$ is a very small Reeb chord between $\hat{\gamma}^i(t_0)$ and $\hat{\gamma}^i(t_1)$. Moreover, the hypothesis about the co--orientations allows us to further assume that 
\begin{itemize}
	\item [(A.1)] $\gamma^s=\hat{\gamma}^{0}\#(\beta^{0}_{0}((1-\frac s\varepsilon)t),v_0((1-\frac s\varepsilon)t))$, for $s\in[0,\varepsilon]$;
	\item [(A.1)'] $\gamma^{1-s}=\hat{\gamma}^{1}\#(\beta^{1}_{0}((1-\frac s\varepsilon)t),v_0((1-\frac s\varepsilon)t))$, for $s\in[0,\varepsilon]$;
	\item [(A.2)] $\gamma^s\equiv\hat{\gamma}^{0}$ is constant for $s\in[\varepsilon,2\varepsilon]$,
	\item [(A.2)'] $\gamma^{1-s}\equiv\hat{\gamma}^{1}$ is constant for $s\in[\varepsilon,2\varepsilon]$.
\end{itemize}

Assume that the fixed cusp point of the curves $\gamma^s$ is placed at time $c_1=\frac{t_0+t_1}{2}$ in the whole family. Choose a capping disk for $\D_i$ for $\gamma^{i}_{\Res(t_0,t_1)}$, $i\in\{0,1\}$. Choose small neighborhoods $U_i$ of each disk. By \cite{EliashbergR3}, $U_1$ is contactomorphic to the standard contact $\R^3$ and there is a contactomorphism $\Phi:U_0\rightarrow U_1$ such that $\Phi(\gamma^{0}(t_0))=\gamma^{1}(t_0)$ and $\Phi(\gamma^{0}(c_1))=\gamma^{1}(c_1)$. By \cite{EliashbergFraser} we can further assume that $\gamma^{1}_{\Res(t_0,t_1)}=\Phi\circ\gamma^{0}_{\Res(t_0,t_1)}$. In other words, the construction that we are going to do is symmetric (there is a symmetry $s\rightarrow 1-s$ for $s\in[0,2\varepsilon]$). Ie all the constructions that we are going to do near $s=0$ immediately translate to $s=1$ by means of the contactomorphism $\Phi$. We do not detail the construction in a neighborhood of $s=1$, just declaring it to be the symmetric construction to the one that we provide for $s=0$ 

Use the disk $\D_0$ to construct a band $B_0:[0,1]^2\rightarrow\R^3$ with the following properties:
	\begin{itemize}
		\item $B_0(0,v)=\beta^{0}_{0}(v)$,
		\item $B_0(1,v)=\hat{\gamma}^0(c_1)$,
		\item $B_0(\sigma,0)=\hat{\gamma}^0((1-\sigma)t_0+\sigma c_1)$,
		\item $B_0(\sigma,1)=\hat{\gamma}^0((1-\sigma)t_1+\sigma c_1)$,
		\item $B_0((0,1)^2)\cap\hat{\gamma}^0(\NS^1)=\emptyset$,
		\item for each $\sigma\in(0,1]$ the map $B_{0}^{\sigma}:[0,1]\rightarrow\R^3,y\mapsto B_0(\sigma,v)$; is an embedding.		
	\end{itemize}

Let $\delta>0$ be small enough in such a way that $\gamma^{s}_{[c_1-\delta,c_1+\delta]}$ is constant (independent of $s$). Consider the linear interpolation $\chi:[0,1]\rightarrow[t_0,c_1-\delta],s\mapsto (1-s)t_0+s(c_1-\delta)$. The family $\tilde{\beta}^{s}_{0}=B^{\chi(s)}_{0}$ defines a family of smooth finger bones joining $\hat{\gamma}^{0}(\chi(s))$ with $\hat{\gamma}^{0}(t_1+t_0-\chi(s))$. Note that $\tilde{\beta}^{0}_{0}=\beta^{0}_{0}$. Use Lemma \ref{TransverseFingers} to homotope this family through finger bones to a family of transverse finger bones $\beta^{s}_{k}$. Observe that $\beta^{0}_{k}=\tilde{\beta}^{0}_k$ by the construction of the interpolation. There is no obstruction to find a family of vector fields $v_s$ in $(\beta^{s})^{*}\xi$ such that $v_0=v$ and $(\beta_s,v_s)$ is a pre--Legendrian finger germ for $\hat{\gamma}$.

Finally, by Lemma \ref{lem:ReebChordStabilized} there exists a non--negative integer $K\in\Z_{\geq0}$ in such a way that we can interpolate through strict Legendrian immersions between $\hat{\gamma}^{0}\#(\beta^{0}_{0},v_0)$ and $\hat{\gamma}^{0}\#(\beta^{0}_{k},v_0)$ just by adding a tower of $K$ DSs to the initial finger (by using Lemma \ref{lem:cuspcreation}). This number $K$ is going to be fixed for the rest of the proof.

\textit{Step 1: Adding DSs to the family.}

Define $\gamma^{s,u}$, $(s,u)\in[0,1]\times[0,\frac13]$; by adding a tower of $K$ DSs (see Figure \ref{CreacionFinger}) to $\gamma^{s,0}=\gamma^{s}$ based at the point $\gamma^s (T(s))$ where $T(s)$ satisfy 
\begin{itemize}
\item $T(s)\equiv t_0+\delta$ for $s\in[0,\varepsilon]$ and, symmetrically, for $s\in[1-\varepsilon,\varepsilon]$;
\item $T(s)$ linearly interpolates between $t_0+\delta$ and $c_1-\delta$ for $s\in[\varepsilon,2\varepsilon]$ and, symmetrically, for $s\in[1-2\varepsilon,1-\varepsilon]$\footnote{Note that the path of Legendrians is constant over $[\varepsilon,2\varepsilon]$ and over $[1-2\varepsilon,1-\varepsilon]$.};
\item $T(s)\equiv c_1-\delta$ for $s\in[2\varepsilon,1-2\varepsilon]$.
\end{itemize}

\textit{Step II: Crossing the first half of the tower of DSs.}

Use Lemma \ref{lem:cuspcreation} to define $\gamma^{s,u}$, $(s,u)\in[0,1]\times[\frac13,\frac23]$; just by crossing half of the tower of DSs through the self--intersection of $\gamma^{i,\frac13}$, $i\in\{0,1\};$ at time $t_0$. Moreover, do this in such a way that $\gamma^{s,u}=\gamma^{s,\frac13}$ for $(s,u)\in[\varepsilon,1-\varepsilon]\times[\frac13,\frac23]$.

\textit{Step III: Creating the self--intersections.}

To define $\gamma^{s,u}$, $(s,u)\in[0,1]\times[\frac23,1]$; proceed as follows:
\begin{itemize}
	\item Define $\gamma^{s,u}$, $s\in[0,\varepsilon]$, in such a way that $\gamma^{s,1}\in\SLegImmG(\R^3)$, $s\in[0,\varepsilon]$ interpolates between $\gamma^{0,\frac23}$ and $\hat{\gamma}^{0}\#(\beta^{0}_{k},v_0)$. This is a direct application of Lemma \ref{lem:ReebChordStabilized} and the fact that we have a tower of $K$ DSs to produce this interpolation. Define $\gamma^{s,u}$, $s\in[1-\varepsilon,1]$, symmetrically.
	\item Use the family of pre--Legendrian finger germs $(\beta^{s}_{k},v^s)$ to define $\gamma^{s,u}$, $s\in[\varepsilon,2\varepsilon]$. Do the symmetric construction to define $\gamma^{s,u}$ for $s\in[1-2\varepsilon,1-\varepsilon]$. Note that the self--intersection of $\gamma^{2\varepsilon,1}$ and, symmetrically, of $\gamma^{1-2\varepsilon,1}$ happens at times $c_1-\delta$ and $c_1+\delta$.
	\item Define $\gamma^{s,u}$, $s\in[2\varepsilon,1-2\varepsilon]$, just by using the pre--Legendrian finger germ $(\beta^{1}_{k},v^{1})$ to create a self--intersection at times $c_1-\delta$ and $c_1+\delta$. Note that this is trivial because the family is constant near the cusp point.
\end{itemize}
This concludes the construction of the required family.
\end{proof}

{\bf \large Case 2: $\{ s_1, \ldots, s_{\nu} \}\neq \emptyset$}.

Reasoning as in the smooth case, we conclude. Let us provide the details. In all the previous proofs we have made deformations in the moduli space $\LegImm(\R^3)$. These deformations over a point $\gamma^{s,u}:\NS^1\rightarrow\R^3$ of the moduli space change the type of Legendrian by creating tangencies at a finite number of values of the parameter $0<u_1<\cdots u_N\leq 1$, $N\geq 0$, ie  $\gamma^{s,u_j}\in\SLegImm(\R^3)$, $s\in[0,1]$ . The tangencies are created identifying at most two points in the domain $z_0(s,u_j),z_1(s,u_j)\in\NS^1$. We assume in all the steps that the points $z_0(s_j,u_k), z_1(s_j,u_k), \tau_{s_j}^{0}, \tau_{s_j}^{1}$, $(j,k)\in\{1,\ldots,\nu\}\times\{1,\ldots,N\}$, are pairwise different. This implies that there is no interference between the deformation and the pair of points $\tau_{s_j}^{0}, \tau_{s_j}^{1}$. Observe that a generic choice of the times $z_0(s,u),z_1(s,u)$ satisfies this assumption. So, we need to change nothing if we add the set $\{s_1,\ldots,s_\nu\}$ in the proofs, except for the fact that we need to make the choices of times (mainly for the finger creation process) generic in the previous sense.

\subsubsection{Index $1$ surgeries.}

\begin{proposition}\label{TwoComponents}
		Let $\D(x,y)\subseteq\R^2$ be the disk of radius $2$ and let $\varphi:\D(x,y)\rightarrow\Hor(\R^4)$ be a map into the horizontal embeddings space such that
		\begin{itemize}
			\item [(i)] the configuration of $\pi_G\circ \varphi$ is given by $\nu+2$ vertical segments, ie \[(\pi_G\circ\varphi)^{-1}(\SLegImm(\R^3))=(\{x=0\}\cup\{x=s_1\}\cup\cdots\cup\{x=1\})\cap\D(x,y),\]
			for $\nu\geq0$ and $0<s_1<\cdots<s_\nu<1$;
			\item [(ii)] the area function has the same sign over the first and the last segments.
		\end{itemize}
	Then, there exists a $1$--parametric family of disks $\varphi_t:\D\rightarrow\Hor(\R^4)$, $t\in[0,1]$, satisfying that
	\begin{itemize}
		\item [(i)] $\varphi_0=\varphi$,
		\item[(ii)] $\varphi_{t |\partial\D}=\varphi_{0 |\partial\D}$,
		\item [(iii)] The configuration of $\pi_G\circ\varphi_1$ is given by the right diagram of Figure \ref{Surgery}.
	\end{itemize}
The changes in the configuration of curves of strict immersions in $\pi_G\circ \varphi_t$ moving from $\pi_G\circ \varphi_0$ to $\pi_G\circ\varphi_1$ correspond to an ambient index one surgery whose attaching points are the middle points $(0,0)$ and $(1,0)$ of the two vertical segments $\{x=0\}$ and $\{x=1\}$.  The attached disk is given by $[0,1]\times[-\varepsilon,\varepsilon]$, with $\varepsilon>0$ small enough, and the surgery replaces the pair of $1$--disks $\{0,1 \}\times[-\varepsilon,\varepsilon]$ by the pair of $1$--disks $[0,1]\times\{\pm\varepsilon\}$. The new $1$--disks, that are curves of strict immersions, have an even number of cusp points and they do not add new zeroes to their the area functions. The other segments $\{x=s_1\},\ldots,\{x=s_\nu\}$ remain untouched in the diagram. Moreover, a finite number of closed circular components of type $S_n=\{(x,y)\in\D:(x-\frac12)^2+y^2=r_{n}^{2}\}$, with $0<r_1<\cdots<r_n<c$. These new closed components do not have cusp points and have positive area function everywhere.
\begin{figure}[h] 
	\includegraphics[scale=0.3]{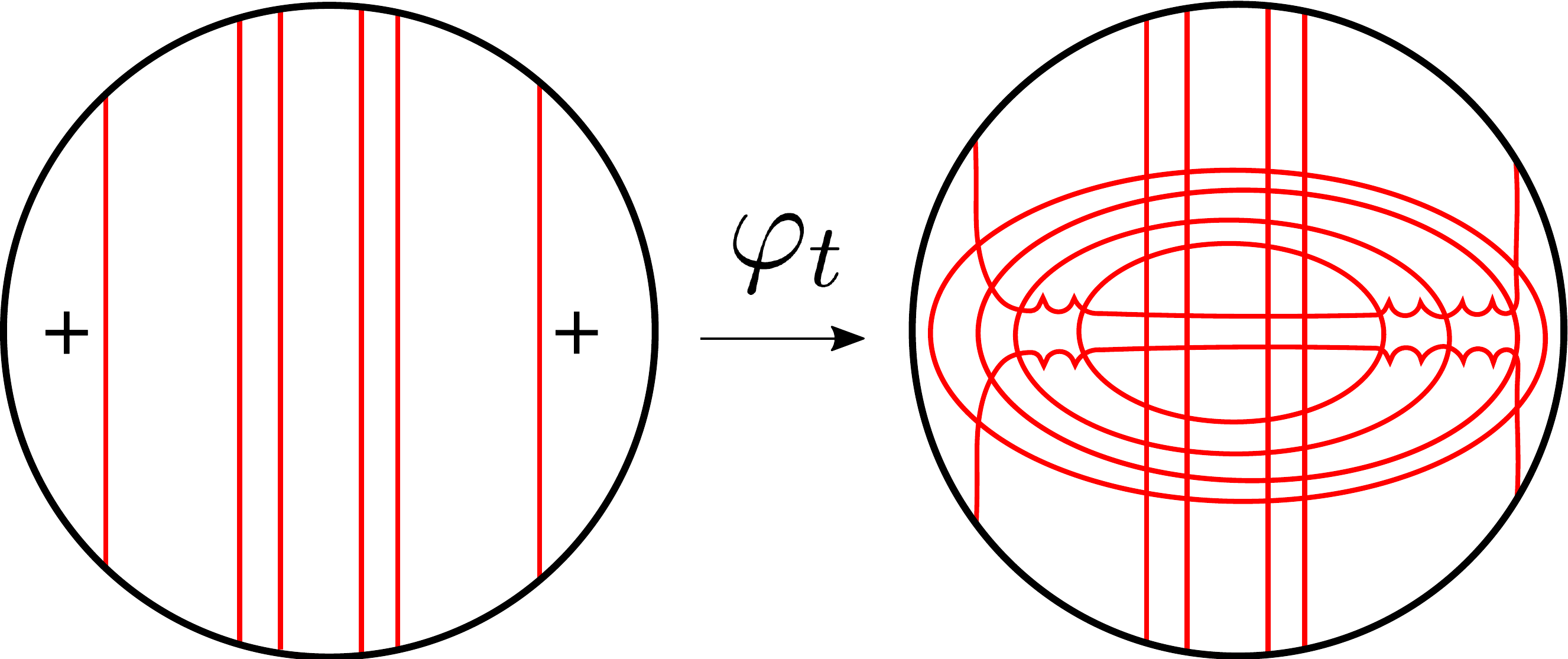}
	\centering
	\caption{Description of the homotopy $\varphi_t$. Observe that the cusp points appear in pairs, therefore no zeroes of the area function appear.\label{Surgery}}
\end{figure}
\end{proposition}
\begin{proof}
Denote the two curves of strict Legendrian immersions in $\D_G=\pi_G\circ\varphi(\D)$ by $C_0$ and $C_1$. Let $i\in\{0,1\}$ and fix $\gamma^i\in C_i$. Choose a path $\gamma^s$, $s\in[-1,2]$, transversal to the strict immersions curves at  $C_i$ such that $C_i\cap\{\gamma^s:s\in[-1,2]\}=\{\gamma^i\}$. 
A neighborhood $K=\{\gamma^{s}_{r}: s\in[-1,2], r\in[-1,1]\}\subseteq\D$ of the path $\gamma^{s}$ may be assumed to satisfy $\gamma^{s}_{r}\equiv\gamma^{s}$ for $s\in[-1,2]$.

Apply Proposition \ref{CreatingLegendrianTangency} to find a disk $E=\{\gamma^{s,u}:s,u\in[0,1]\}$ such that
\begin{itemize}
	\item {}$\gamma^{s,0}=\gamma^{s}$, $s\in[0,1]$.
	\item {} $C=\{\gamma^{s,u}\in\SLegImm(\R^3):(s,u)\in\{0,1\}\times[0,1]\cup[0,1]\times\{1\}\}$ is a continuous curve of strict Legendrian immersions, with self-intersection times $t_{0}^{s,u},t_{1}^{s,u}\in\NS^
1$, $t_{0}^{s,u}\neq t_{1}^{s,u}$.
\item {} $\mathcal{S}_j=\{\gamma^{s_j,u}:u\in[0,1]$, $j\in\{1,\ldots,\nu\}$, is a curve of strict Legendrian immersions,
	\item {}  There is a finite number of curves $\mathcal{C}_n$ of strict Legendrian immersions defined by \[ \mathcal{C}_n=\{\gamma^{s,u}:(s,u)\in[0,1]\times\{u_n\}\},\]
	where $0<u_1<\cdots<u_n<1$. The self-intersection times of $\gamma^{s,u_n}$ are $\tau_{0}^{s,u_n},\tau_{1}^{s,u_n}\in\NS^1$, with $\tau_{0}^{s,u_n}\neq \tau_{1}^{s,u_n}$. These curves do not have cusp points.
	\item {} $\tau_{i}^{s,u_n}\neq t_{j}^{s,u_n}$, for each $i,j\in\{0,1\}$ and $(s,u_n)\in C\cap\mathcal{C}_n$.	
\end{itemize}
Apply Lemma \ref{DeltaBarquillo} to lift $E$ into a disk of horizontal immersions. Now, apply again Lemma \ref{DeltaBarquillo} over $C$ in $t_{0}^{s,u}-\varepsilon$ and $t_{0}^{s,u}+\varepsilon$, for $\varepsilon>0$ small enough in order to increase the area function at will. Do likewise over the self--intersection points in $\mathcal{C}_n$. Note that there is no interaction between the two deformations because the self--intersection times are different. 

By construction the outward normal vector field $X$ to $C=\partial E$ coincides with the standard orientation. Push $C$ in the direction of $X$ to obtain an extended disk $\tilde{E}=\{\gamma^{s,u}:s\in[-1,2], u\in[0,2]\}\supseteq E$ in $\LegImm(\R^3)$, which lifts to $\Hor(\R^4)$, such that $\gamma^{s,0}=\gamma^{s}$ for each $s\in[-1, 2]$.  We define the map
$$\psi: [-1,2]\times[0,2]\times [-1,1]\to \Hor(\R^4): (s,u,r) \mapsto \gamma^{s,u}_r= \gamma^{s,u}.$$

See Figure \ref{edificio}, to visualize the piecewise linear disk  that we are about to define in $[-1,2]\times[0,2]\times [-1,1] \subseteq \R^3$. First, we define the base (a closed piecewise linear curve in $\R^2$): 
$$ \Delta_B = \{ (s,r) \in [-1,2]\times [-1,1]: |r|\leq g(s) \}, $$
for
$$g(s) = \left\{ \begin{array}{cc} 1+s, &\text{if }  s \in [-1,0], \\ 1, &\text{if }  s\in [0,1],  \\ 2-s, &\text{if }  s\in [1,2]. \end{array} \right. $$
And now the piecewise linear graph in $\R^3$ as 
$$ \Delta = \left\{ (s,u,r) \in [-1,2]\times [1,2] \times [-1,1]: (s,r) \in \Delta_B,  u= \left\{ \begin{array}{cc} 2+s-|r|, & \text{if } s \in[-1,0], \\ 2-|r|, & \text{if } s \in [0,1], \\ 3-s-|r|, & \text{if } s \in [1,2].  \end{array} \right. \right\}. $$
\begin{figure}[h] 
	\includegraphics[scale=0.55]{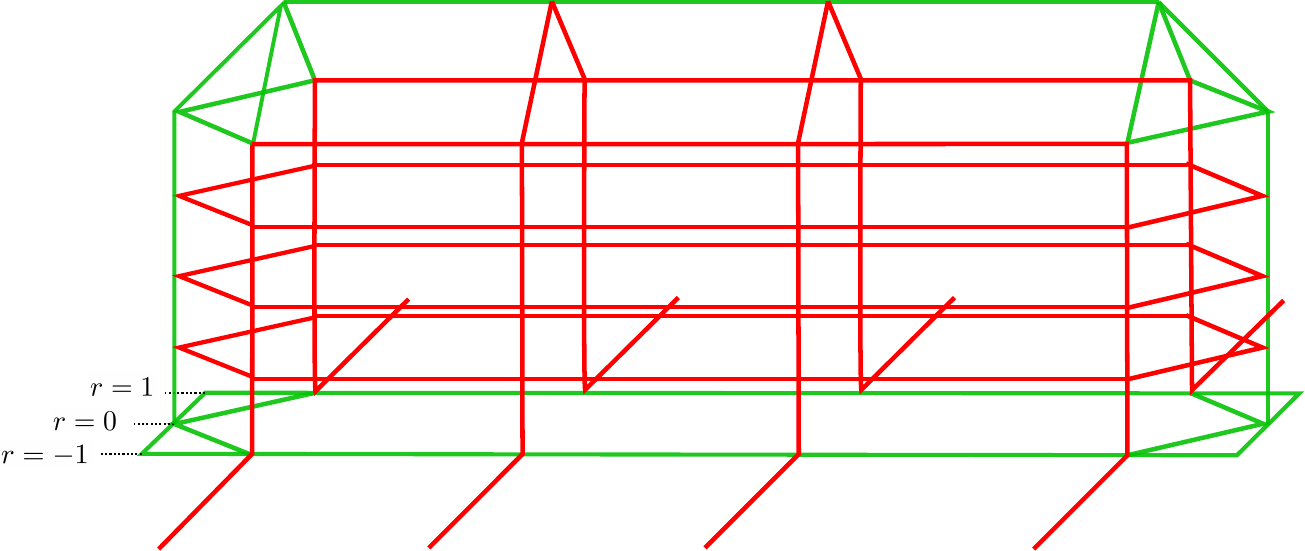}
	\centering
	\caption{Construction of $\tilde{\D}$.}\label{edificio}
\end{figure}
We define a new (piecewise linear) disk $\hat\Delta =( \partial (\Delta_B) \times [0,1]) \cup \Delta \subseteq [-1,2]\times [0,2] \times [-1,1]$. Observe that the boundary of $\hat\Delta$ coincides with $\partial\Delta_B$, therefore we can define a new disk  $\tilde{\D}= (\D \setminus \Delta_B) \cup \hat\Delta$ gluing along the boundaries. We can define a new
morphism
\begin{eqnarray*}
\varphi_1: \tilde{\D} & \to & \Hor(\R^4) \\
p & \mapsto & \left\{ \begin{array}{cc} \varphi_0(p), & p \in (\D \setminus \Delta_B), \\ \psi(p), & p \in \hat\Delta. \end{array} \right. 
\end{eqnarray*}
It is left to check that $\varphi_0$ and $\varphi_1$ are homotopic. This is performed by building any homotopy between the two disks $\Delta_B$ and $\hat{\Delta}$ relative to the boundary $\partial \Delta_B$ inside $[-1,2]\times [0,2] \times [-1,1]$.
\end{proof}

\subsection{Globalization.} \label{sub:global}
We use the theory developed in the last two Subsections to prove the main theorem of Section \ref{h--Principle}. We first adapt a previous construction.

\begin{lemma} \label{lem:hellcircle} [Corollary of Lemma \ref{lem:cuspcreation}]

	Let $\varphi:\D(x,y)\rightarrow\HorImm(\R^4)$ be a disk of horizontal immersions such that:
	\begin{itemize}
		\item [(i)] $\varphi(\partial \D)\subseteq\Hor(\R^4)$,
		\item [(ii)] $(\pi_G\circ\varphi)^{-1}(\SLegImm(\R^3))=C$ is a connected curve with a unique zero of the area function $\beta\in C$.
	\end{itemize}
Then, there exists a $1$--parametric family of disks $\varphi^{t}:\D\rightarrow\HorImm(\R^4)$, $t\in[0,1]$, satisfying
\begin{itemize}
	\item [(i)] $\varphi^0=\varphi$,
	\item [(ii)] $\varphi^{t}_{|\partial\D}=\varphi^{0}_{|\partial\D}$,
	\item [(iii)] there is a $C^0$ deformation $C_t$, $t\in[0,1]$, of $C=C_0$ and a family of connected curves $S_t$ such that \begin{itemize}
		\item  $(\pi_G\circ\varphi^t)^{-1}(\SLegImm(\R^3))=C_t\cup S_t$,
		\item  $\beta\in C_t$ is the unique zero of the area function in $\varphi^{t}$ and the curves $S_t$ encircle $\beta$,
		\item $C_t=C_0$ for $t\in[0,1/2)$,
		\item $C_1$ is a DS of $C_{1/2}$ in the moduli (see Lemma \ref{lem:cuspcreation}) and $\beta$ lies between the first and the second cusp.
		\item $S_t=\emptyset$ for $t\in[0,1/2)$,
		\item $S_t=\{x^2+y^2=(\frac{1}{2}-t)^2\}$ for $t\in[1/2,3/4]$ and
		\item $S_t=S_{3/4}$ for $t\in[3/4,1]$.
	\end{itemize}
\end{itemize}
See Figure \ref{fig:hellcirc}. 
\end{lemma}
\begin{proof}
	\begin{figure}[h] 
		\includegraphics[scale=0.23]{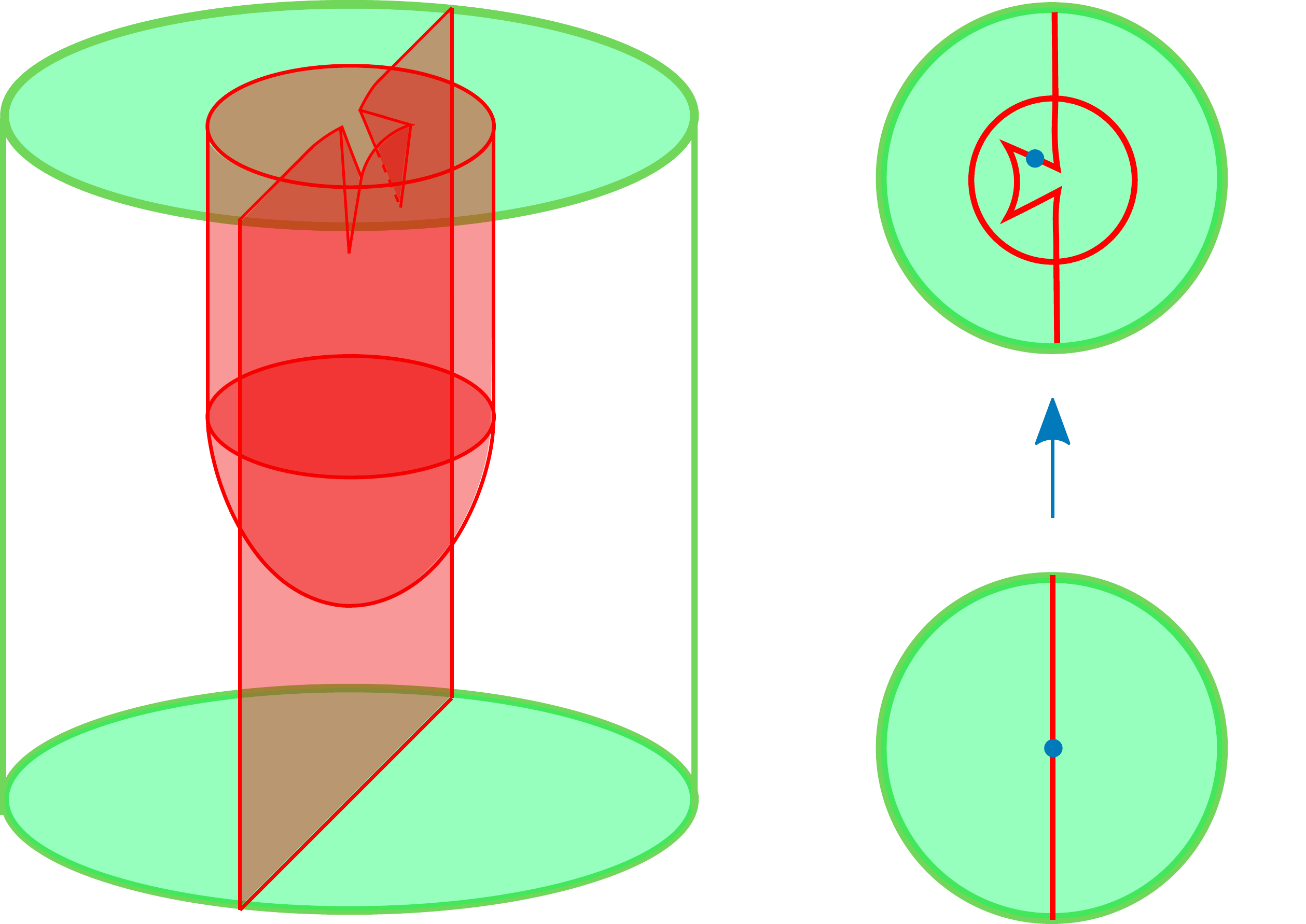}
		\centering
		\caption{Schematic description of Lemma \ref{lem:hellcircle}. The blue point in the middle of the circle represents $\beta$ and the red middle circle in the cylinder the self--intersection (with non zero area) created in the DS.}\label{fig:hellcirc}
	\end{figure}
Since the problem is local we may assume that the disk is small and the self--intersection times for the Legendrian projections in the curve $C$ are constant and equal to $t_0$ and $t_1$. Choose $t_2\in \NS^1$ a time value different from $t_0$ and $t_1$. Use polar coordinates $\D=\D(r,\theta)$ and define parametrically for each horizontal immersion $\varphi(r,\theta)$ a DS that is denoted $B = \D \times [0,1/2] \to  \HorImm(\R^4)$, where $B(r, \theta,h)$, $h\in [0,1/2]$, is the DS of $\varphi(r,\theta)$. 

Now we apply Lemma \ref{lem:cuspcreation} in order to further deform and create $4$ cusp points in the moduli: just sliding one branch ($t_0$) over the other  ($t_1$). Apply Lemma \ref{DeltaBarquillo} to place the zero in the required position. This second deformation produces a family of disks $B(r, \theta,h)$ with $h\in [1/2,1]$.

Define
$$ L= \{ B(r,\theta, h): \, r=\frac12, h \in [0,1] \} \subseteq \im(B),$$
$$ T= \D\left(r \leq \frac12\right) \times \{ h=1 \} \subseteq \im(B),$$
that builds a piecewise linear disk  $\D' = L \cup T$ whose boundary coincides with the boundary of $\D(r \leq \frac12) \subseteq \D$. Replace $\D(r \leq \frac12)$ by $\D'$, ie define $\tilde\D = (\D \setminus \D(r \leq \frac12)) \cup \D'$. The associated diagram of strict immersions conforms the statement of the Lemma. See Figure \ref{fig:hellcirc}.
\end{proof}

\begin{proof}[Proof of Theorem \ref{AreaIsomorphism}.]
	Let $\gamma^\theta$ be a loop of horizontal embeddings which lies in $\pi_1(\Hor(\R^4))_0$, ie $\Rot_L (\gamma^\theta)=0$. Assume that $\Ar(\gamma^\theta)=0$, we must check that $\gamma^\theta$ is trivial.
	
	Let $\D$ be a disk in $\HorImm(\R^4)$ whose boundary is given by $\{\gamma^{\theta}\}$. Apply successively Proposition \ref{zeroesSameComponent} to assume that each curve of strict Legendrian immersions in $\D_G$ has at most one zero of the area function. Since $\Ar(\gamma^\theta)=0$ the number of zeroes is even. If $\D_G$ has not zeroes then we are done. Hence, assume that $\D_G$ has a positive number of zeroes and denote them  $\beta_{0},\ldots,\beta_{2n-1}$. 
	
	Let $C_k\subseteq\D_G$ be the curve of strict Legendrian immersions that contains $\beta_k$, $k\in\{0,\ldots,2n-1\}$. It is sufficient to show that we can do a surgery along any two curves $C_0$ and $C_1$ to produce a new curve of strict Legendrian immersions that contains the zeroes $\beta_0$ and $\beta_1$ without creating any new zero. Application of Proposition \ref{zeroesSameComponent} will provide the required result. 
	
	To formalize the previous discussion, note that $C_i$ is naturally oriented by the gradient vector  $g_i$ of the Area function at $\beta_i$. Together with the standard co-orientation of $\Sigma^{1,0}$ at $\beta_i$, denoted by $w_i$, it provides an oriented basis at $T_{\beta_i}\D$ namely  $\{ g_i, w_i \}$. We can assume that the orientations induced at the two points $\beta_0$ and $\beta_1$ are the same. If it is not the case, apply Lemma \ref{lem:hellcircle} in order to change the orientation of the curve $C_0$. This works because the Area function changes signs when crossing a cusp in the curve $C_0$ and thus the gradient of the area function over the zero gets multiplied by $-1$. 	
	
Select a pair of points $\alpha_i$ arbitrarily close to $\beta_i$, whose area function is positive. Fix an embedded curve $\alpha:[0,1]\to \D$ such that $\alpha(i)=\alpha_i$ and it is transverse to $C_i$. Moreover, we assume that $\frac{d\alpha_t}{dt}_{|t=0}$ defines the standard coorientation at $\alpha_0 \in \Sigma^{1,0}$. Respectively $\frac{d\alpha_t}{dt}_{|t=1}$ defines the opposite orientation to the standard one (see Figure \ref{fig:uturn} to learn how to deal with the wrong coorientation case).

	\begin{figure}[h] 
		\includegraphics[scale=0.35]{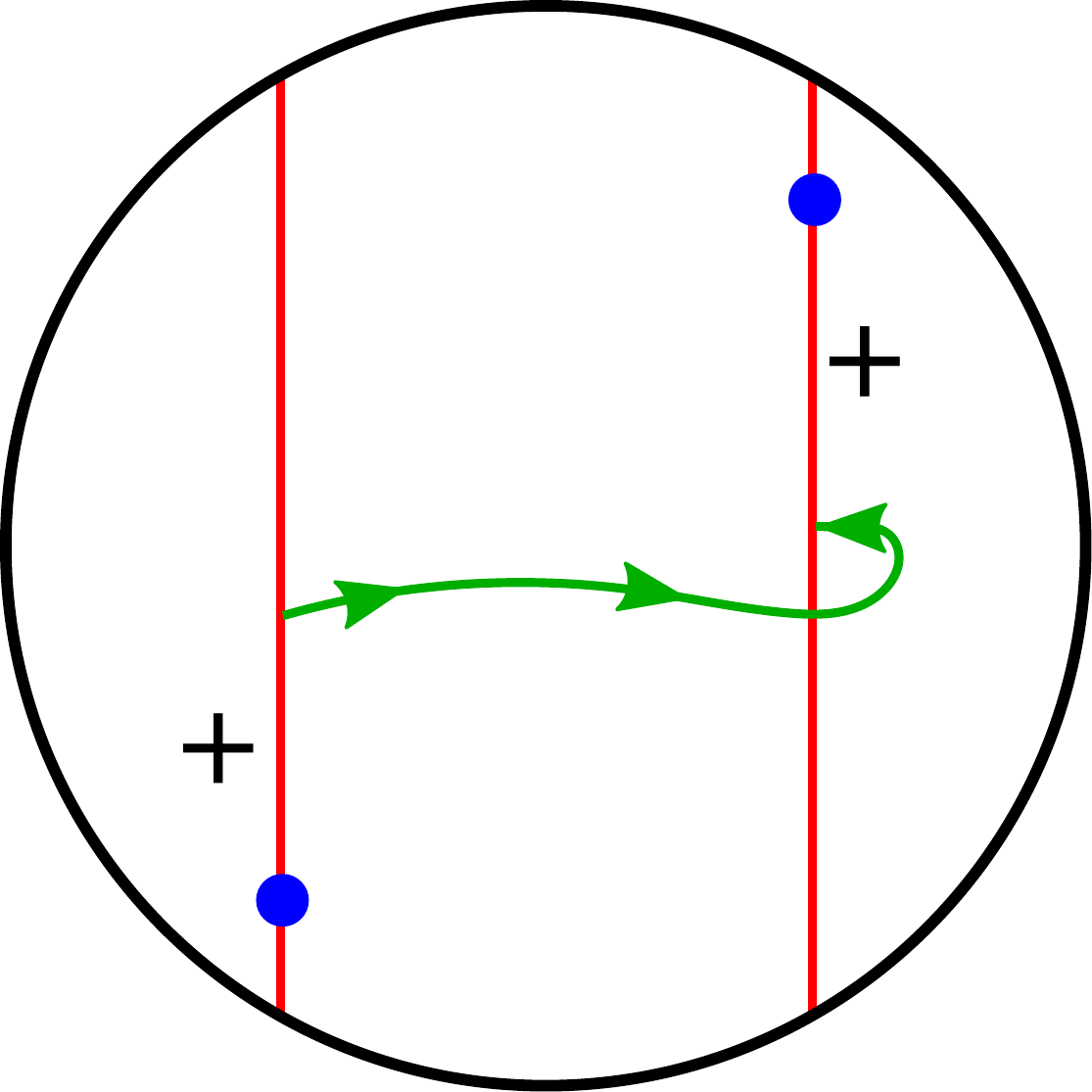}
		\centering
		\caption{Embedded path connecting $\alpha_0$ and $\alpha_1$ conforming to the coorientations rule.}\label{fig:uturn}
	\end{figure}

A small tubular neighborhood $\D_{\alpha}$ of the path satisfies the hypothesis of Proposition \ref{TwoComponents}. This implies that we can get a deformed disk with the same number of zeroes such that the deformation is performed compactly supported on the disk $\D_{\alpha}$ following the Diagram explained in Figure \ref{Surgery}. This new diagram makes the two zeroes $\beta_i$ be part of the same connected component of curves of strict immersions, see Figure \ref{m=0}. This is ensured by the compatibility of the orientations that we have enforced in the basis $\{ g_i, w_i \}$ of $T_{\beta_i}\D$. Apply Proposition \ref{zeroesSameComponent} to cancel the two zeroes.
	
\begin{figure}[h] 
		\includegraphics[scale=0.05]{lastfig1.pdf}
		\centering
		\caption{Schematic description of the neighborhood of $\alpha(t)$. The cusps between the concentric circles are not drawn to make the picture easier to understand.}\label{m=0}
	\end{figure}
\end{proof}


\begin{thebibliography}{9999}
	\bibitem{Ada} J. Adachi. \emph{Classification of horizontal loops in standard Engel space}. Int. Math. Res. Not. 2007; Vol. 2007: article ID rnm008, 29 pages, doi:10.1093/imrn/rnm008.
	
	\bibitem{Adams} C. C. Adams. The knot book. An elementary introduction to the mathematical theory of knots. American Mathematical Society, Providence, RI, 2004. xiv+307 pp. 
	
	\bibitem{Arnold} V. I. Arnold, S. M. Gusein-Zade, A. N. Varchenko. Singularities of Differentiable Maps. Volume I. Birkh\"{a}user, Boston, Basel, Sttutgart. 1985.
	
	\bibitem{Bennequin} D. Bennequin. \emph{Entrelacements et \'{e}quations de Pfaff.} Third Schnepfenried geometry conference, Vol. 1 (Schnepfenried, 1982), 87--161, 
	Ast\'{e}risque, 107-108, Soc. Math. France, Paris, 1983. 
	
	\bibitem{Budney} R. Budney. \emph{A family of embedding spaces.} Geometry \& Topology 13 (2008): 41--83.
	
	\bibitem{CasalsdelPino} R. Casals, A. del Pino. \emph{Classification of Engel knots}. arXiv:1710.11034.
	
	\bibitem{Chekanov} Y. Chekanov. \emph{Differential algebra of Legendrian links.} Invent. Math. 150 (2002), no. 3, 441--483.
	
	\bibitem{DingGeiges} F. Ding, H. Geiges. \emph{The diffeotopy group of $\NS^1\times\NS^2$ via contact topology.} Compos. Math. 146 (2010), no. 4, 1096?1112. 
	
	\bibitem{EliashbergR3} Y. Eliashberg. \emph{Classification of contact structures on $\R^3$}. Internat. Math. Res. Notices 1993, no. 3, 87--91. 		
	
	\bibitem{EliashbergFraser} Y. Eliashberg, M. Fraser. \emph{Topologically trivial Legendrian knots.} J. Symplectic Geom. 7 (2009), no. 2, 77--127.
	
	\bibitem{EliashMisch} Y. Eliashberg, N. Mishachev, Introduction to the $h$--Principle. Graduated Studies in Mathematics, 48. AMS, Providence, RI, 2002. 
	
	\bibitem{Etnyre} J. Etnyre. \emph{Legendrian and transversal knots.} Handbook of knot theory, 105--185, Elsevier B. V., Amsterdam, 2005.
	
	\bibitem{EtnyreHonda} J. Etnyre and K. Honda. \emph{Knots and contact geometry I: Torus knots and the figure eight knot.} J. Symplectic
	Geom. 1 (2001), 63-120.
	
	\bibitem{FuchsTabachnikov} D. Fuchs, S. Tabachnikov. \emph{Invariants of Legendrian and transverse knots in the standard contact space.} Topology 36 (1997), no. 5, 1025--1053.
	
	
	\bibitem{GeigesCont} H. Geiges. An Introduction to Contact Topology. Cambr. Studies in Adv. Math. 109. Cambr. Univ. Press 2008.
	
	\bibitem{GeigesLoops} H. Geiges. \emph{Horizontal loops in Engel space.} Math. Ann. Oct 2008, Vol. 342, no. 2, 291--296.
	
	\bibitem{ham} R. S. Hamilton. \emph{The inverse function theorem of Nash and Moser.} Bull. Amer. Math. Soc. (N.S.) 7 (1982), no. 1, 65–222.
	
	\bibitem{HatcherSmale} A. Hatcher. \emph{A proof of the Smale conjecture, ${\rm Diff(\NS^3)}\simeq O(4)$.} 
	Ann. of Math. (2) 117 (1983), no. 3, 553--607. 
	
	\bibitem{HatcherBook} A. Hatcher,
	Algebraic topology. Cambridge University Press, Cambridge, 2002. xii+544 pp. ISBN: 0-521-79160-X; 0-521-79540-0 
	
	\bibitem{hatcher} A. Hatcher. \emph{Spaces of knots.} arXiv:math/9909095.
	
	\bibitem{SmaleHirsch} M.W. Hirsch. \emph{Immersions of manifolds.} Trans. Amer. Math. Soc. 93 1959 242--276.
	
	\bibitem{Igusa} K. Igusa. \emph{Higher singularities of smooth functions are unnecessary.} 
	Ann. of Math. (2) 119 (1984), no. 1, 1?58.
	
	\bibitem{kalman} T. K\'{a}lm\'{a}n. \emph{Contact homology and one parameter families of Legendrian knots.} Geometry \& Topology 9 (2005): 2013--2078.
	
	\bibitem{Lang} S. Lang, Fundamentals of differential geometry. Graduate Texts in Mathematics, 191. Springer-Verlag, New York, 1999. xviii+535 pp. ISBN: 0-387-98593-X 53-01 (58-01)
	
	\bibitem{Lenhard} L. Ng. \emph{Computable Legendrian invariants.} Topology 42 (2003), no. 1, 55--82.  
	
	\bibitem{Milnor}  J. Milnor. Morse theory. Based on lecture notes by M. Spivak and R. Wells. Annals of Mathematics Studies, No. 51 Princeton University Press, Princeton, N.J. 1963 vi+153 pp. 57.50 (53.72)
	
	\bibitem{PinoPresas} A. del Pino, F. Presas. \emph{Flexibility for tangent and transverse immersions in Engel manifolds.} arXiv:1609.09306.
	

	\bibitem{Sabloff} J. M. Sabloff, M. G. Sullivan. \emph{Families of Legendrian submanifolds via generating families.} Quantum Topol. 7 (2016), no. 4, 639--668.
	
	\bibitem{Szabo} P. Ozsv\'{a}th, Z. Szab\'{o}, D. Thurston. \emph{Legendrian knots, transverse knots and combinatorial Floer homology}. Geometry \& Topology 12 (2008), no. 2, 941--980.	
	
\end{thebibliography}
\end{document}